\documentclass[12pt]{amsart}
\usepackage{amssymb,latexsym, amsmath, amsxtra,psfig}

\newtheorem{lemma}{Lemma}

\newtheorem{theorem}{Theorem}

\newtheorem{rems}{Remarks}

\newcommand{\abs}[1]{\left| #1 \right|}
\newcommand{\br}[1]{\left( #1 \right)}
\newcommand{\cbr}[1]{\left\{ #1 \right\}}

\newcommand{\cj}[1]{\overline{#1}}       
\newcommand{\mellin}[2]{{M}\br{#1;#2}}
\newcommand{\bl}{\boldsymbol{\lambda}}
\newcommand{\Elambda}[1]{E_{\bl}\br{#1}}
\newcommand{\Ecjlambda}[1]{E_{\cj{\bl}}\br{#1}}
\newcommand{\Glambda}[2]{G_{\bl}\br{#1,#2}}
\newcommand{\Gcjlambda}[2]{G_{\cj{\bl}}\br{#1,#2}}
\newcommand{\Gammalambda}[1]{\Gamma_{\bl}\br{#1}}
\newcommand{\gammalambda}[1]{\gamma_{\bl}\br{#1}}
\newcommand{\C}{{\mathbb{C}}}

\newcommand{\Q}{{\mathbb{Q}}}
\newcommand{\R}{{\mathbb{R}}}
\newcommand{\Z}{{\mathbb{Z}}}
\newcommand{\UH}{{\mathbb{H}}}
\newcommand{\e}{\varepsilon}
\newcommand{\sgn}{\text{sgn}}

\newcommand{\lhs}{l.h.s.\ }
\newcommand{\rhs}{r.h.s.\ }

\newcommand{\gam}[1]{\gamma_{#1}^{(d)}}

\def\yr #1#2#3#4{(#1#2#3#4)}
\def\publ{\relax}
\def\publaddr{\relax}

\def\jour{\relax}
\def\pages{\relax}

\begin{document}

\title[Methods and experiments]{Computational methods and experiments in analytic number theory}
\author{Michael Rubinstein}
%\date{June 2004}
\maketitle

\section{Introduction}

We cover some useful techniques in computational aspects of analytic 
number theory, with specific emphasis on ideas relevant to the evaluation of $L$-functions. 
These techniques overlap considerably with basic methods from analytic number theory. On the
elementary side, summation by parts, Euler-Maclaurin summation, and Mobius inversion
play a prominent role. In the slightly less elementary sphere, we find tools from analysis, such 
as Poisson summation, generating function methods, Cauchy's residue theorem, asymptotic methods, 
and the fast Fourier transform. We then describe conjectures and experiments that connect 
number theory and random matrix theory.

\section{Basic methods}

\subsection{Summation by parts}

Summation by parts can be viewed as a discrete form of integration by parts.
Let $f$ be a function from $\Z^+$ to $\R$ or $\C$, and
$g$ a real or complex valued function of a real variable. Then
\begin{equation}
\label{eq:summation by parts}
   \sum_{1 \leq n \leq x} f(n) g(n) =
   \left(
       \sum_{1 \leq n \leq x} f(n)
   \right) g(x) -
   \int_1^x 
   \left(
       \sum_{1 \leq n \leq t} f(n)
   \right) g'(t) dt.
\end{equation}
Here we are assuming that 
$g'$ exists and is continous on $[1,x]$. One verifies this identity by 
writing the integral as $\int_1^2+\int_2^3+\ldots+\int_{\lfloor x \rfloor}^x$, noticing that the sum in each
integral is constant on each open interval, integrating, and telescoping.
Although our integral begins at $t=1$, 
it is sometimes convenient to start earlier, for example at $t=0$. This doesn't change the value of the
integral, the sum in the integrand being empty if $t<1$. Formula~(\ref{eq:summation by parts}) can
also be interpreted in terms of the Stieltjes integral.

A slightly more general form of partial summation is over a set 
$\left\{\lambda_1,\lambda_2,\ldots\right\}$
of increasing real numbers:
\begin{equation}
   \notag
   \sum_{\lambda_n \leq x} f(n) g(\lambda_n) =
   \left(
       \sum_{\lambda_n \leq x} f(n)
   \right) g(x) -
   \int_{\lambda_1}^x
   \left(
       \sum_{\lambda_n \leq t} f(n)
   \right) g'(t) dt.
\end{equation}

As an application, let $$\pi(x) = \sum_{p \leq x} 1$$ denote the number of primes less than or equal
to $x$, and $$\theta(x) = \sum_{p \leq x} \log p$$ denote the number primes
up to $x$ with each prime weighted by its logarithm. The famous equivalence between 
$\pi(x) \sim x/\log x$ and $\theta(x) \sim x$ can be verified using partial summation.
Write
\begin{equation}
    \notag
    \pi(x)= \sum_{p \leq x} \log p \frac{1}{\log p} = 
    \theta(x) \frac{1}{\log x}  
    + \int_2^x
      \theta(t) \frac{dt}{t(\log t)^2},
\end{equation}
from which it follows that if $\theta(x) \sim x $ then $\pi(x) \sim x/\log x$.
The converse follows from
\begin{equation}
    \notag
    \theta(x) = \sum_{p \leq x} 1\cdot\log p =
    \pi(x) \log x - \int_2^x \pi(t) \frac{dt}{t}.
\end{equation}

\subsection{Euler-Maclaurin summation}

A powerful application of partial summation occurs when the function $f(n)$ is identically 
equal to 1 and
the function $g(t)$ is many times differentiable. In that case, summation by parts 
specializes to the Euler Maclaurin formula which involves one summation by parts
with $f(n)=1$ followed by repeated integration by parts.
For $a,b \in {\mathbb{Z}}$,$a<b$, partial summation gives
\begin{equation}
    \notag
    \sum_{a < n \leq b} g(n) = (b-a) g(b) 
    - \int_a^b  (\lfloor t \rfloor -a) g'(t) dt
    = b g(b) - a g(a) 
    - \int_a^b  \lfloor t \rfloor g'(t) dt.
\end{equation}
Here, we have chosen to start the integral at $t=a$, rather than at $t=a+1$.
Writing $\lfloor t \rfloor = t -\{t\}$, with $\{t\}$ the fractional
of $t$ we get
\begin{equation}
    \notag
    \sum_{a < n \leq b} g(n)  
    = \int_a^b g(t) dt
    + \int_a^b \{t\} g'(t) dt.
\end{equation}
The second term on the r.h.s. should be viewed as the necessary correction
that arises from replacing the sum on the left with an integral. 

The next step is to write $\{t\}=1/2+ (\{t\}-1/2)$, the latter term having nicer properties than 
$\{t\}$, for example being odd and also having zero constant term in its Fourier
expansion. So
\begin{equation}
\label{eq:euler mac 1}
    \sum_{a < n \leq b} g(n)  
    = \int_a^b g(t) dt + \frac{1}{2}(g(b)-g(a))
    + \int_a^b (\{t\}-1/2) g'(t) dt.
\end{equation}
Integrating the second integral repeatedly by parts
leads naturally to the introduction of Bernoulli polynomials, named after
Jacob Bernoulli (1654-1705), who discovered them in connection to the problem of 
studying sums of positive integer powers of consecutive integers. 
During the 1730's Euler (1707-1783), who studied mathematics from Jacob's brother Johann (1667-1748),
developed the summation formula being described in connection with computing
reciprocals of powers and Euler's constant.

\subsubsection{Bernoulli Polynomials}

The Bernoulli polynomials are defined recursively by the following relations
\begin{eqnarray}
    \notag
    B_0(t) &=& 1 \\
    \notag
    B_k'(t) &=& k B_{k-1}(t), \quad k \geq 1 \\
    \notag
    \int_0^1 B_k(t) dt &=& 0 , \quad k \geq 1.
\end{eqnarray}
The second equation determines $B_k(t)$ recursively up to the constant
term, and the third equation fixes the constant.
The first few Bernoulli polynomials are listed in Table~\ref{tab:bernoulli poly}.
\begin{table}[h]
\center{
\begin{tabular}{|c|c|}
\hline
$k$ & $B_k(t)$ \\ \hline
0 & $1$ \\ \hline
1 & $t-1/2$ \\ \hline
2 & $t^2-t+1/6$ \\ \hline
3 & $t^3-3/2t^2+1/2t$ \\ \hline
4 & $t^4-2t^3+t^2-1/30$ \\ \hline
5 & $t^5-5/2t^4+5/3t^3-1/6t$ \\ \hline
\end{tabular}
}
\caption{The first few Bernoulli polynomials}
\label{tab:bernoulli poly}
\end{table}

Let $B_k=B_k(0)$ denote the constant term of $B_k(t)$. $B_k$ is called the
$k$-th Bernoulli number.
We state basic properties
of the Bernoulli polynomials. Expansion in terms of Bernoulli numbers:
\begin{equation}
    \notag
    B_k(t) = \sum_0^k {k\choose m} B_{k-m} t^m, \quad k \geq 0
\end{equation}
Generating function:
\begin{equation}
    \notag
    \frac{ze^{zt}}{e^z-1} = \sum_0^\infty B_k(t) z^k/k!, \quad |z| < 2 \pi
\end{equation}
Fourier series:
\begin{equation}
    \label{eq:Fourier 1}
    B_1(\{t\}) = -\frac{1}{\pi} \sum_1^\infty \frac{\sin(2\pi mt)}{m}, \quad t \not\in {\mathbb{Z}}
\end{equation}
\begin{equation}
    \label{eq:bernoulli Fourier series}
    B_k(\{t\}) = -k! \sum_{m\neq 0} \frac{e^{2\pi imt}}{(2\pi i m)^k}, \quad k \geq 2.
\end{equation}
Functional equation:
\begin{equation}
    \notag
    B_k(t) = (-1)^k B_k(1-t), \quad k \geq 0
\end{equation}
Difference equation:
\begin{equation}
    \label{eq:bernoulli diff}
    \frac{B_{k+1}(t+1)-B_{k+1}(t)}{k+1} = t^k, \quad k \geq 0
\end{equation}
Special values:
\begin{equation}
    \notag
    B_k(1) = \begin{cases}    
                 (-1)^k B_k(0), &\quad k \geq 0 \\
                 0, \quad &\text{$k$ odd, $k \geq 3$}  \\
                 1/2, &\quad k=1 
             \end{cases}
\end{equation}
i.e.
\begin{equation}
    \label{eq:bernoulli special values}
    B_k(1)=B_k(0), \quad \text{unless $k=1$}
\end{equation}
Recursion:
\begin{equation}
    \notag
    \sum_{m=0}^{k-1} {k\choose m} B_m = 0, \quad k \geq 2
\end{equation}
Equation~(\ref{eq:Fourier 1}) can be obtained directly.
The other formulae can be verified using the defining relations
and induction. 

Property~(\ref{eq:bernoulli Fourier series}) can be used to obtain a formula
for $\zeta(2m)$.
Let 
\begin{equation}
    \notag
    \zeta(s) = \sum_1^\infty n^{-s}, \quad \Re{s} > 1.
\end{equation}
Taking $t=0$, $k=2m$, even, in the Fourier
expansion of $B_k(\{t\})$ gives
\begin{equation}
    \notag
    B_{2m} = \frac{(-1)^{m+1} (2m)!}{(2\pi)^{2m}} 2 \zeta(2m)
\end{equation}
so that
\begin{equation}
    \notag
    \zeta(2m) = \frac{(-1)^{m+1} (2\pi)^{2m}}{2(2m)!} B_{2m},
\end{equation}
a formula discovered by Euler.
Because $\zeta(2m) \to 1$ as $m \to \infty$, we have
\begin{equation}
    \notag
    B_{2m} \sim 
   \frac{(-1)^{m+1} 2 (2m)!}{(2\pi)^{2m}} 
\end{equation}
as $m \to \infty$.

\subsubsection{Euler-Maclaurin continued}

Returning to (\ref{eq:euler mac 1}), we write
\begin{equation}
    \notag
    \int_a^b (\{t\}-1/2) g'(t) dt = 
    \int_a^b B_1(\{t\}) g'(t) dt.
\end{equation}
Breaking up the integral 
$\int_a^b = \int_a^{a+1} + \int_{a+1}^{a+2} + \ldots \int_{b-1}^b$,
integrating by parts, and noting that $B_2(1)=B_2(0)$, 
we get, assuming that $g^{(2)}$ exists and is continous on $[a,b]$, 
\begin{equation}
    \notag
    \frac{B_2}{2}(g'(b)-g'(a)) - \int_a^b \frac{B_2(\{t\})}{2} g^{(2)}(t) dt.
\end{equation}
Repeating, using $B_k(1) = B_k(0)$ if $k \geq 2$, leads to the Euler-Maclaurin summation
formula. Let $K$ be a positive integer. Assume that $g^{(K)}$ exists and
is continous on $[a,b]$. Then
\begin{eqnarray}
    \sum_{a < n \leq b} g(n)  
    = \int_a^b g(t) dt + 
    \sum_{k=1}^K \frac{(-1)^k B_k}{k!}(g^{(k-1)}(b)-g^{(k-1)}(a)) \notag \\
    \notag
    +\frac{(-1)^{K+1}}{K!}
    \int_a^b B_K(\{t\}) g^{(K)}(t) dt.
\end{eqnarray}

\subsubsection{Application: Sums of consecutive powers}

We apply Euler-Macluarin summation to obtain Bernoulli's formula
for sums of powers of consecutive integers. Let $r \geq 0$ be
an integer. Then
\begin{equation}
   \notag
    \sum_{n=1}^N n^r = \frac{B_{r+1}(N+1)-B_{r+1}(1)}{r+1}.
\end{equation}
We can verify this directly using property (\ref{eq:bernoulli diff}),
substituting $n=1,2,\ldots,N$, and telescoping. However, it is instructive to
apply the Euler-Maclaurin formula, which, once begun, carries through in an 
automatic fashion. In this example, we have $g(t)=t^r$. Notice that 
$g^{(r+1)}(t)=0$, and that
\begin{equation}
    g^{(m)}(N)-g^{(m)}(0) = \begin{cases}
                                r(r-1)\ldots(r-m+1) N^{r-m}, \quad m \leq r-1 \notag \\
                                0, \quad m \geq r. \notag 
                            \end{cases}
\end{equation}
If $m=0$ we set $r(r-1)\ldots(r-m+1)=1$. Then
\begin{eqnarray}
    \sum_{n=1}^N n^r 
    &=& \int_0^N t^r dt +
    \sum_{k=1}^r \frac{(-1)^k B_k}{k!} r(r-1)\ldots(r-k+2) N^{r-k+1} \notag \\
    &=& \int_0^N t^r dt +
    \sum_{k=1}^r \frac{(-1)^k B_k}{r-k+1} {r\choose k} N^{r-k+1} \notag \\
    &=& \int_0^N 
    \sum_{k=0}^r (-1)^k B_k {r\choose k} t^{r-k} dt
    = \int_0^N (-1)^r B_r(-t) dt = \int_0^N B_r(t+1) dt 
    \notag \\
    &=& (B_{r+1}(N+1)-B_{r+1}(1))/(r+1). \notag
\end{eqnarray}
If $r \geq 1$, the last line simplifies according to~(\ref{eq:bernoulli special values})
and equals
\begin{equation}
    \notag
    \frac{B_{r+1}(N+1)-B_{r+1}}{r+1}.
\end{equation}

\subsubsection{Application: $\zeta(s)$}
The Euler-Maclaurin formula can be used to obtain the analytic continuation of $\zeta(s)$ and 
also provides a useful expansion for its numeric evaluation.
Consider
\begin{equation}
    \notag
    \sum_1^N n^{-s} = 
    1 + \sum_2^N n^{-s} 
\end{equation}
with $\Re{s}>1$. We have started the sum at $n=2$ rather than
$n=1$ to avoid difficulties near $t=0$ below.
Applying Euler-Maclaurin summation, with $g(t)=t^{-s}$,
$g^{(m)}(t) = (-1)^m s (s+1)\ldots(s+m-1) t^{-s-m}$, we get
\begin{eqnarray}
    \sum_1^N n^{-s} = 
    1 + \int_1^N t^{-s} dt
    &-&\sum_{k=1}^K \frac{B_k}{k} {s+k-2 \choose k-1} (N^{-s-k+1}-1) \notag \\
    \notag
    &-&{s+K-1 \choose K}
    \int_1^N B_K(\{t\}) t^{-s-K} dt.
\end{eqnarray}
Evaluating the first integral,
taking the limit as $N \to \infty$, with $\Re{s}>1$, we get
\begin{equation}
    \label{eq:zeta euler mac}
    \zeta(s) = \frac{1}{s-1} + \frac{1}{2} 
    +\sum_2^K {s+k-2 \choose k-1} \frac{B_k}{k}
    -{s+K-1 \choose K} \int_1^\infty B_K(\{t\}) t^{-s-K} dt.
\end{equation}
While we started with $\Re{s}>1$, the r.h.s. is meromorphic for 
$\Re{s} > -K+1$, so gives the meromorphic continuation of $\zeta(s)$ in
this region, with the only pole being the simple pole at $s=1$.

Taking $s=2-K$, $K \geq 2$,
\begin{equation}
    \notag
    \zeta(2-K) = 1 - \frac{1}{K-1} \sum_{k=0}^{K-1} (-1)^k {K-1 \choose k} B_k
    = \frac{(-1)^K B_{K-1}}{K-1}.
\end{equation}
Thus,
\begin{eqnarray}
    \zeta(1-2m) &=& -B_{2m}/(2m), \quad m=1,2,3,\ldots \notag \\
    \zeta(-2m) &=& 0, \quad m=1,2,3,\ldots \notag \\
    \notag
    \zeta(0) &=& -1/2.
\end{eqnarray}
Applying the functional equation for $\zeta$ (see for example
Roger Heath-Brown's notes)
\begin{equation}
    \notag
    \pi^{-s/2} \Gamma(s/2) \zeta(s) =
    \pi^{-(1-s)/2} \Gamma((1-s)/2) \zeta(1-s) 
\end{equation}
and
\begin{equation}
    \notag
    \Gamma(1/2) = \pi^{1/2} = 
    \frac{-1}{2}
    \frac{-3}{2}
    \frac{-5}{2} \cdots
    \frac{-(2m-1)}{2}
    \Gamma(1/2-m)
\end{equation}
gives another proof of Euler's identity
\begin{equation}
    \notag
    \label{eq:zeta 2m}
    \zeta(2m) = \frac{(-1)^{m+1} (2\pi)^{2m}}{2(2m)!} B_{2m}, \quad m \geq 1.
\end{equation}

\subsubsection{Computing $\zeta(s)$ using Euler-Maclaurin summation}

Next we describe how to adapt the above to obtain a practical method for
numerically evaluating $\zeta(s)$.
From a computational perspective, the following works better than 
using~(\ref{eq:zeta euler mac}).
Let $N$ be a large positive integer, proportional
in size to $|s|$. We will make this more explicit shortly.
For $\Re{s}>1$, write
\begin{equation}
    \label{eq:zeta break up sum}
    \zeta(s) = 
     \sum_1^\infty n^{-s} = 
     \sum_1^N n^{-s} +
     \sum_{N+1}^\infty n^{-s}.
\end{equation}
The first sum on the r.h.s. is evaluated term by term, while the second sum is 
evaluated using Euler-Maclaurin summation
\begin{equation}
    \label{eq:zeta tail sum}
     \sum_{N+1}^\infty n^{-s} = 
    \frac{N^{1-s}}{s-1} 
    +\sum_1^K {s+k-2 \choose k-1} \frac{B_k}{k} N^{-s-k+1}
    -{s+K-1 \choose K} \int_N^\infty B_K(\{t\}) t^{-s-K} dt.
\end{equation}
As before, the r.h.s. above gives the meromorphic continuation of the l.h.s. to
$\Re{s} > -K+1$.
Breaking up the sum over $n$ in this fashion allows us to throw away the
integral on the r.h.s., and obtain sharp estimates for its neglected contribution.
First, from property (\ref{eq:bernoulli Fourier series}), 
\begin{equation}
    \notag
    |B_K(\{t\})| \leq \frac{K!}{(2\pi)^K} 2 \zeta(K).
\end{equation}
It is convenient to take $K=2K_0$, even, in which case we have
from~(\ref{eq:zeta 2m})
\begin{equation}
    \notag
    |B_{2K_0}(\{t\})| \leq B_{2K_0}.
\end{equation}
Therefore, for $s=\sigma+i\tau, \sigma > -2K_0+1$,
\begin{eqnarray}
    &&\left|
        {s+2K_0-1 \choose 2K_0} \int_N^\infty B_{2K_0}(\{t\}) t^{-s-2K_0} dt
    \right| \notag\\
    &\leq&  
    \left|
        {s+2K_0-1 \choose 2K_0} B_{2K_0}
    \right| 
    \frac{N^{-\sigma-2K_0+1}}{\sigma+2K_0-1} \notag
    \\
    &=&
    \frac{|s+2K_0-1|}{\sigma+2K_0-1} 
    \left|
        {s+2K_0-2 \choose 2K_0-1} \frac{B_{2K_0}}{2K_0}
    \right|
    N^{-\sigma-2K_0+1} \notag \\
   &=&
    \frac{|s+2K_0-1|}{\sigma+2K_0-1} |\text{last term taken}|. 
    \notag
\end{eqnarray}
A more precise estimate follows by comparison of $B_{2K_0}$ with
$\zeta(2K_0)$, and we have that the remainder is 
\begin{equation}
    \notag
    \leq \frac{\zeta(2K_0)}{\pi N^\sigma} 
    \frac{|s+2K_0-1|}{\sigma+2K_0-1}
    \prod_{j=0}^{2K_0-2} \frac{|s+j|}{2\pi N}.
\end{equation}
We start to win when $2\pi N$ is bigger than $|s|,|s+1|, \ldots, |s+2K_0-2|$.
There are two parameters which we need to choose: $K_0$ and $N$, and we also need to specify
the number of digits accuracy, Digits, we desire. For example, with $\sigma \geq 1/2$,
taking
\begin{equation}
    \notag
    2\pi N \geq 10 |s+2K_0-2|
\end{equation}
with
\begin{equation}
    \notag
    2K_0-1 > \text{Digits} + \frac{1}{2}\log_{10}(|s+2K_0-1|)
\end{equation}
achieves the desired accuracy. The main work involves the computation
of the sum $\sum_1^N n^{-s}$ consisting of $O(|s|)$ terms.
Later we will examine the Riemann-Siegel formula and its smoothed variants
which, for $\zeta(s)$, involves a main sum of $O(|s|^{1/2})$ terms. However, for 
high precision
evaluation of $\zeta(s)$, especially with $s$ closer to the real axis,
the Euler-Maclaurin formula remains an ideal method allowing for
sharp and rigorous error estimates and reasonable efficiency. 

In fact, we can turn the above scheme into a computation
involving $O(|s|^{1/2})$ operations but requiring $O((\text{Digits}+\log{|s|})\log{|s|})$
precision due to cancellation that occurs.
In~(\ref{eq:zeta break up sum}) choose $N\sim |10s/(2\pi)|^{1/2}$,
and assume that $\Re{s}\geq 1/2$.
Expand $B_K(\{t\})$ into
its Fourier series (\ref{eq:bernoulli Fourier series}). We only need $M=O(|s|^{1/2})$ 
terms of the Fourier expansion
to assure a contribution from the neglected terms smaller than the desired precision.
Each term contributes
\begin{equation}
    \label{eq:contr Fourier series term}
    K! {s+K-1 \choose K} \frac{1}{(2\pi i m)^K}
    \int_N^\infty e^{2\pi i m t} t^{-s-K} dt,
\end{equation}
so the neglected terms contribute altogether less than
\begin{equation}
    \notag
    \frac{1}{N^\sigma}
    \frac{|s+K-1|}{\sigma+K-1}
    \left( \prod_{j=0}^{K-2} \frac{|s+j|}{2\pi N} \right)
    \left( \sum_{M+1}^\infty \frac{2}{m^K} \right) .
\end{equation}
Here we have combined the $\pm m$ terms together. Comparing to an integral,
the sum above is $<2/((K-1)M^{K-1})$ and so the neglected terms contribute
less than
\begin{equation}
    \notag
    \frac{2}{(K-1)N^\sigma}
    \frac{|s+K-1|}{\sigma+K-1}
    \prod_{j=0}^{K-2} \frac{|s+j|}{2\pi MN}.
\end{equation}
We start to win when $2\pi MN$ exceeds $|s|,\ldots,|s+K-2|$.
For $\sigma \geq 1/2$,
choose $K>\text{Digits} + \log_{10}(|s+K-1|) +1$ and $M=N$ with
\begin{equation}
    \notag
    2\pi MN \geq 10 |s+K-2|.
\end{equation}
Asymptotically, we can improve the above choices so as to achieve 
$M=N\sim |s|^{1/2}/(2\pi)$, the same as in the Riemann-Siegel formula.
The only drawback is that extra precision as described above is needed.
The individual terms summed in (\ref{eq:zeta tail sum})
are somewhat large in comparison to the final result, this coming form
the binomial coefficients which have numerator $(s+k-2)\ldots(s+1)s$,
and this leads to cancellation.

Finally to compute the contribution to the Fourier expansion from the 
terms with $|m|\leq M$, we assume that $4|K$ so that the terms $\pm m$ together
involve in (\ref{eq:contr Fourier series term}) the integral
$$
    \int_N^\infty \cos(2\pi m t) t^{-s-K} dt
    =
    (2\pi m)^{s+K-1}
    \int_{2\pi mN}^\infty \cos(u) u^{-s-K} du.
$$
This can be expressed in terms of the incomplete $\Gamma$ function
$$
    \int_w^\infty \cos(u) u^{z-1} du
    =
    \frac{1}{2}
    \left(
        e^{-\pi i z/2} \Gamma(z,iw) +
        e^{\pi i z/2} \Gamma(z,-iw) 
    \right).
$$
See Section 3 which describes properties of the incomplete $\Gamma$ function and
methods for its evaluation.

The Euler-Maclaurin formula can also be used to
evaluate Dirichlet
$L$-functions. It works in that case due to the periodic nature of 
the corresponding Dirichlet coefficients.
For general $L$-functions, there are smoothed Riemann-Siegel
type formulae. These are described later.

\subsection{Mobius inversion with an application to sums and products over primes}

Computations in analytic number theory often involve evaluating
sums or products over primes. For example, let $\pi_2(x)$ denote
the number of twin primes $(p,p+2)$, with $p$ and $p+2$ both prime and less than or equal to $x$.
The famous conjecture of Hardy and Littlewood predicts that
\begin{equation}
    \notag
    \pi_2(x) \sim 
    2 \prod_{p>2} \frac{p(p-2)}{(p-1)^2}
    \frac{x}{(\log{x})^2}.
\end{equation}
Generally, it is easier to deal with a sum rather than a product,
so we turn this product over primes into a sum by expressing it as
\begin{equation}
    \notag
    \exp \left(
        \sum_{p>2} \log(1-2/p) - 2 \log(1-1/p)
    \right).
\end{equation}
Letting $f(p)=\log(1-2/p) - 2 \log(1-1/p)$, we have
\begin{equation}
    \notag
    f(p) = - \sum_{m=1}^\infty \frac{2^m-2}{mp^m}
\end{equation}
hence
\begin{equation}
    \label{eq:twin 1}
    \sum_{p>2} f(p) = 
    - \sum_{m=1}^\infty \frac{2^m-2}{m} (h(m)-1/2^m)
\end{equation}
with
\begin{equation}
    \notag
    h(s) = \sum_p p^{-s}, \quad \Re{s}>1.
\end{equation}
We therefore need an efficient method for computing $h(m)$. This will be
dealt with below.
Notice that $h(m)-1/2^m \sim 1/3^m$ so the sum on the r.h.s.
of~(\ref{eq:twin 1}) converges exponentially fast. We can achieve
faster convergence by writing
\begin{equation}
    \notag
    \sum_{p>2} f(p) =
    \sum_{2<p\leq P} f(p) + \sum_{p>P} f(p),
\end{equation}
summing the terms in the first sum, and expressing the second sum as
\begin{equation}
    \notag
    - \sum_{m=1}^\infty \frac{2^m-2}{m} (h(m)-1/2^m-\ldots -1/P^m).
\end{equation}

A second example involves the computation of constants that arise in conjectures
for moments of $\zeta(s)$. The Keating-Snaith conjecture~\cite{KeS} asserts that
\begin{equation}
     \label{eq:moment 1}
     M_k(T) := \frac{1}{T} \int_o^T |\zeta(1/2+it)|^{2k} dt 
     \sim \frac{a_k g_k}{k^2!} (\log{T})^{k^2}
\end{equation}
\begin{eqnarray}
    \label{eq:zeta_ak}
    a_k& =  &\prod_p \left(1-\frac{1}{p}\right)^{k^2}
    \sum_{m=0}^{\infty} {m+k-1 \choose m }^2
    p^{-m} \\ 
    & =  &\prod_p\left(1-\frac{1}{p}\right)^{(k-1)^2}
    \ \sum_{j=0}^{k-1}\binom{k-1}{j}^2p^{-j}  
    \notag
\end{eqnarray}
and
\begin{equation}
    \notag
    g_k=  k^2 ! \prod_{j=0}^{k-1}\frac{j!}{(k+j)!}.
\end{equation}
The placement of $k^2!$ is to ensure that $g_k$ is an integer
\cite{CF}. Keating and Snaith also provide a conjecture for complex
values, $\Re{k} > -1/2$, of which the above is a special case.
Keating and Snaith used random matrix theory to
identify the factor $g_k$. The form of~(\ref{eq:moment 1}), without identifying
$g_k$, was conjectured by Conrey and Ghosh~\cite{CG}.

The above conjecture gives the leading term for the asymptotics for
the moments of $|\zeta(1/2+it)|$.
In~\cite{CFKRS} a conjecture is given for the full asymptotics
of $M_k(T)$:
\begin{equation}
    \notag
    M_k(T) \sim \sum_{r=0}^{k^2} c_r(k) (\log{T})^{k^2-r}
\end{equation}
where $c_0(k)=  a_k g_k /k^2!$ coincides with the Keating-Snaith leading term
and where the degree $k^2$ polynomial is given implicitly as an elaborate
multiple residue.
Explicit expressions for $c_r(k)$ are worked out in~\cite{CFKRS3} and are
given as $c_0(k)$ times complicated rational functions in $k$,
generalized Euler constants, and sums over primes involving $\log(p)$,
$_2F_1(k,k,1;p^{-1})$ and its derivatives.
One method for computing the $c_r(k)$'s 
involves as part of a single step the computation of sums of the form
\begin{equation}
    \label{eq:der h(s)}
    \sum_p \frac{(\log{p})^r}{p^m}, \quad m=2,3,4,\ldots \quad r=0,1,2,\ldots.
\end{equation}

We now describe how to efficiently compute $h(s)= \sum_p p^{-s}$ and the sums 
in (\ref{eq:der h(s)}). Take
the logarithm of
\begin{equation}
    \notag
    \zeta(s) = \prod_p (1-p^{-s})^{-1}, \quad \Re{s}>1
\end{equation}
and apply the Taylor series for $\log(1-x)$ to get
\begin{equation}
    \label{eq:h(s) 1}
    \log{\zeta(s)} = \sum_{m=1}^\infty \frac{1}{m} h(ms), \quad \Re{s}>1.
\end{equation}
Let $\mu(n)$, the Mobius $\mu$ function, denote the Dirichlet coefficients
of $1/\zeta(s)$:
\begin{equation}
    \notag
    1/\zeta(s) = \prod_p (1-p^{-s}) = \sum_1^\infty \mu(n) n^{-s}.
\end{equation}
We have
\begin{equation}
    \notag
    \mu(n)=
    \begin{cases}
        0 \quad \text{if $n$ is divisible by the square of an integer $>1$ } \\
        (-1)^{\text{number of prime factors of $n$}} \quad \text{if $n$ is squarefree}
    \end{cases}
\end{equation}
and 
\begin{equation}
    \notag
    \sum_{n|r} \mu(n) = 
    \begin{cases}
        1 \quad \text{if $r=1$} \\
        0 \quad \text{otherwise}.
    \end{cases}
\end{equation}
The last property can be proven by writing the sum of the left as
$\prod_{p|r}(1-1)$, 
and it allows us to invert equation~(\ref{eq:h(s) 1})
\begin{eqnarray}
    \sum_{m=1}^\infty \frac{\mu(m)}{m} \log\zeta(ms) 
    &=&
    \sum_{m=1}^\infty \frac{\mu(m)}{m} \sum_{n=1}^\infty \frac{h(mns)}{n} \notag \\
    &=& \sum_{r=1}^\infty \frac{h(rs)}{r} \sum_{m|r} \mu(m) = h(s),
    \notag
\end{eqnarray}
i.e.
\begin{eqnarray}
    \label{eq:mobius 1}
    \sum_p p^{-s} = \sum_{m=1}^\infty \frac{\mu(m)}{m} \log\zeta(ms).
\end{eqnarray}
This is an example of Mobius inversion, and expresses $h(s)$ as a sum involving
$\zeta$. Mobius inversion can be interpreted as a form of the sieve of Eratosthenes. 

Notice that $\zeta(ms)=1+2^{-ms}+3^{-ms}+\ldots$ tends to 1, and hence
$\log\zeta(ms)$ tends to 0, exponentially fast
as $m \to \infty$. Therefore, the number of terms needed on the r.h.s. 
of~(\ref{eq:mobius 1}) is proportional to the desired precision. 

To compute the series appearing in (\ref{eq:der h(s)}) we can 
differentiate $h(s)$ $r$~times, obtaining
\begin{equation} 
    \label{eq:mobius 2}
    \sum_p \frac{(\log{p})^r}{p^s} =
    (-1)^r \sum_{m=1}^\infty 
    \frac{\mu(m)}{m} (\log\zeta(ms))^{(r)}.
\end{equation}

In both (\ref{eq:mobius 1}) and (\ref{eq:mobius 2}), we can use Euler-Maclaurin summation to 
compute $\zeta$ and its derivatives. 
The paper of Henri Cohen~\cite{C} is a good reference for computations
involving sums or products of primes.

\subsection{Poisson summation as a tool for numerical integration}

Let $f \in L^1({\mathbb{R}})$ and let
\begin{equation}
    \notag
    \hat{f}(y) = \int_{-\infty}^\infty f(t) e^{-2\pi i y t} dt.
\end{equation}
denote its Fourier transform. 
The Poisson summation formula asserts, for $f,\hat{f} \in L^1({\mathbb{R}})$
and of bounded variation, that
\begin{equation}
    \notag
    \sum_{n=-\infty}^{\infty} f(n) = \sum_{n=-\infty}^{\infty} \hat{f}(n).
\end{equation}
We often encounter the Poisson summation formula as a potent theoretical
tool in analytic number theory. For example, the functional equations of the 
Riemann $\zeta$ function and of the Dedekind $\eta$ function can be derived
by exploiting Poisson summation.
However, Poisson summation is often overlooked
in the setting of numerical integration where it provides justification
for carrying out certain numerical integrals in a very naive way.

Let $\Delta >0$. By a change of variable
\begin{equation}
    \notag
    \Delta \sum_{n=-\infty}^{\infty} f(n\Delta) = \sum_{n=-\infty}^{\infty} \hat{f}(n/\Delta)
    = \hat{f}(0) + \sum_{n \neq 0} \hat{f}(n/\Delta)
\end{equation}
so that
\begin{equation}
    \notag
    \int_{-\infty}^\infty f(t)dt
    -\Delta \sum_{n=-\infty}^{\infty} f(n\Delta)
    =-\sum_{n \neq 0} \hat{f}(n/\Delta)
\end{equation}
tells us how closely the Riemann sum $\Delta \sum_{n=-\infty}^{\infty} f(n\Delta)$
approximates the integral $\int_{-\infty}^\infty f(t)dt$.

The main point is that if $\hat{f}$ is rapidly decreasing then we get enormous 
accuracy from the Riemann sum, even with $\Delta$ not too small.
For example, with $\Delta=1/10$, the 
first contribution comes from $\hat{f}(\pm 10)$ which can be extremely small
if $\hat{f}$ decreases sufficiently fast.

As a simple application, 
let $f(t)=\exp(-t^2/2)$. Then $\hat{f}(y) = \sqrt{2\pi} \exp(-2\pi^2y^2)$,
and so 
\begin{equation}
    \notag
    \sum_{n \neq 0} \hat{f}(n/\Delta) =  O(\exp(-2\pi^2/\Delta^2)).
\end{equation}
Therefore
\begin{equation} 
    \notag
    \int_{-\infty}^\infty \exp(-t^2/2)dt
    -\Delta \sum_{n=-\infty}^{\infty} \exp(-(n\Delta)^2/2)
    = O(\exp(-2\pi^2/\Delta^2)).
\end{equation}
As everyone knows, the integral on the l.h.s. equals $\sqrt{2\pi}$.
Taking $\Delta=1/10$, we therefore get
\begin{equation}
    \notag
    \Delta \sum_{n=-\infty}^{\infty} \exp(-(n\Delta)^2/2)
    = \sqrt{2\pi} + \epsilon
\end{equation}
with $\epsilon \approx 10^{-857}$. We can truncate the sum over $n$ roughly when
\begin{equation}
    \notag
    \frac{(n\Delta)^2}{2}
    > 
    \frac{2\pi^2}{\Delta^2},
\end{equation}
i.e. when $n > 2\pi/\Delta^2$. So only 628 terms (combine $\pm n$) are needed to
evaluate $\sqrt{2\pi}$ to about 857 decimal place accuracy!

This method can be applied to the problem of computing certain
probability distributions that arise in random matrix theory.
Let $U$ be an $N\times N$ unitary matrix, with eigenvalues
$\exp(i\theta_1),\ldots\exp(i\theta_N)$, and characteristic equation
\begin{equation}
    \notag
    Z(U,\theta) = \prod_1^N (\exp(i\theta) -\exp(i\theta_n))
\end{equation}
evaluated on the unit circle at the point $\exp(i\theta)$. In making their conjecture for the
moments of $|\zeta(1/2+it)|$, Keating and Snaith~\cite{KeS}
studied the analogous random matrix theory problem of evaluating the 
moments of $|Z(U,\theta)|$, averaged according to Haar measure on $\text{U}(N)$.
The characteristic function of a matrix is a class function that only
depends on the eigenvalues of the matrix. For class functions,
the Weyl integration formula gives Haar measure in terms of the
eigenangles, the invariant probability measure on $\text{U}(N)$ being
\begin{equation}  
    \notag
     \frac{1}{(2\pi)^N N!} 
     \prod_{1\leq j < m \leq N}
     |e^{i\theta_j}-e^{i\theta_m}|^2
     d\theta_1 \ldots d\theta_N.
\end{equation}
Therefore, $M_N(r)$, the $r$th moment of $|Z(U,\theta)|$, is given by
\begin{equation}
     \notag
     M_N(r) =
     \frac{1}{(2\pi)^N N!} 
     \int_0^{2\pi}\ldots
     \int_0^{2\pi}
     \prod_{1\leq j < m \leq N}
     |e^{i\theta_j}-e^{i\theta_m}|^2
     |Z(U,\theta)|^r
     d\theta_1 \ldots d\theta_N,
\end{equation} 
for $\Re{r} > -1$.
This integral happens to be a special case of Selberg's integral, and
Keating and Snaith consequently determined that
\begin{equation}
    \notag
    M_N(r) = \prod_{j=1}^N \frac{\Gamma(j) \Gamma(j+r)}{\Gamma(j+r/2)^2}.
\end{equation}
Notice that this does not depend on $\theta$.

Say we are interested in computing the probability distribution
of $|Z(U,\theta)|$. One can recover the probability density function
from the moments as follows.
We can express the moments of $|Z(U,\theta)|$ in terms of
its probability density function. Let
\begin{equation}
    \notag
    \text{prob}( 0 \leq a \leq |Z(U,\theta)| \leq b) = \int_a^b p_N(t) dt.
\end{equation}
Then 
\begin{equation}
    \label{eq:mellin transform M_N(r)}
    M_N(r) = \int_0^\infty p_N(t) t^r dt
\end{equation}
is a Mellin transform, and taking the inverse Mellin transform we get
\begin{equation} 
    \label{eq:inverse mellin p(t)}
    p_N(t) = \frac{1}{2\pi i t} 
    \int_{\nu-i\infty}^{\nu+i\infty} \prod_{j=1}^N \frac{\Gamma(j) \Gamma(j+r)}{\Gamma(j+r/2)^2}
    t^{-r} dr
\end{equation}
with $\nu$ to the right of the poles of $M_N(r)$, $\nu>-1$. There is an extra $1/t$ in
front of the integral since the Mellin transform (\ref{eq:mellin transform M_N(r)}) 
is evaluated at $r$ rather than at $r-1$. 

To compute $p_N(t)$ we could shift the line integral to the left
picking up residues at the poles of $M_N(r)$, but as $N$ grows this becomes 
burdensome. Instead, we can compute the inverse Mellin 
transform~(\ref{eq:inverse mellin p(t)}) as a simple Riemann sum.

Changing variables we have
\begin{equation} 
    \notag
    p_N(t) = 
    \frac{1}{2\pi t} 
    \int_{-\infty}^{\infty} \prod_{j=1}^N \frac{\Gamma(j) \Gamma(j+\nu+iy)}{\Gamma(j+(\nu+iy)/2)^2}
    t^{-\nu-iy} dy.
\end{equation}
Let 
\begin{equation} 
    \notag
    f_t(y) = 
    \frac{1}{2\pi} 
    \prod_{j=1}^N \frac{\Gamma(j) \Gamma(j+\nu+iy)}{\Gamma(j+(\nu+iy)/2)^2}
    t^{-\nu-1-iy}.
\end{equation}
This function also depends on $\nu$ and $N$, but we do not include them 
explicitly on the l.h.s.  so as to simplify our notation.
The above integral equals
\begin{equation} 
    \label{eq:pN(t) ft integral}
    p_N(t) = 
    \int_{-\infty}^{\infty} f_t(y) dy.
\end{equation}
To estimate the error in computing this integral as a Riemann sum using
increments of size $\Delta$,
we need bounds on the Fourier transform
\begin{equation} 
    \label{eq:hat ft(y)}
    \hat{f}_t(u) = 
    \int_{-\infty}^{\infty}
    f_t(y) e^{-2\pi iuy} dy.
\end{equation}
However,
$$
    f_t(y) e^{-2 \pi i uy} = f_{te^{2\pi u}}(y) e^{2\pi u(\nu+1)}
$$
and so
$$
    \hat{f}_t(u) = e^{2\pi u(\nu+1)} p_N(t e^{2\pi u}).
$$
Now, $p_N(t)$ is supported in $[0,2^N]$, because
$0 \leq |Z(U,\theta)| \leq 2^N$.
Hence if
$u > (N\log{2}-\log{t})/(2\pi)$ then
$\hat{f}_t(u)=0$. Thus, for $0 < t < 2^N$, if we evaluate
(\ref{eq:pN(t) ft integral}) as a Riemann sum with step size 
$\Delta < 2\pi/(N\log{2}-\log{t})$ the error is 
$$
    \sum_{n \neq 0} \hat{f}_t(n/\Delta) =
    \sum_{n < 0} \hat{f}_t(n/\Delta) 
$$
since the terms with $n>0$ are all zero. On the other hand, with $n <0$
we get
$$
    \hat{f}_t(-|n|/\Delta) = e^{-2\pi(\nu+1)|n|/\Delta} p_N(t e^{-2\pi |n|/\Delta})
    \leq e^{-2\pi(\nu+1)|n|/\Delta} p_{\text{max}}
$$
where $p_{\text{max}}$ denotes the maximum of $p_N(t)$ (an upper bound for
$p_N(t)$ can be obtained from~(\ref{eq:inverse mellin p(t)})).

Therefore, choosing
$$
    \Delta = \frac{2\pi}{\text{Digits}\log{10} + N \log{2} - \log{t}}
$$
and setting $\nu=0$ we have
$$
    \hat{f}_t(-|n|/\Delta) < (10^{-\text{Digits}}2^{-N}t)^{|n|} p_{\text{max}}
$$
Summing over $n=-1,-2,-3,\ldots$ we get an overall bound of
$$
    10^{-\text{Digits}} p_{\text{max}}/(1-10^{-\text{Digits}}) \approx 10^{-\text{Digits}} p_{\text{max}}.
$$

We could  choose $\nu$ to be larger, i.e. shift our line integral 
(\ref{eq:inverse mellin p(t)}) to the right, and thus achieve 
more rapid decay of $\hat{f}_t(u)$ as $u \to -\infty$. However, this leads to
precision issues. As $\nu$ increases, the integrand in~(\ref{eq:inverse mellin p(t)})
increases in size, yet $p_N(t)$ remains constant for given $N$ and $t$. Therefore
cancellation must occur when we evaluate the Riemann sum and higher precision is
needed to capture this cancellation. We leave it as an excercise to determine
the amount of precision needed for a given value of $\nu$.

Another application appears in~\cite{RS} where 
Poisson summation is used to compute, on a logarithmic scale,
the probability that $\pi(x)$, the number of primes up to $x$,
exceeds $\text{Li}(x) = \int_2^\infty dt/\log(t)$.
The answer turns out to be $.00000026\ldots$

Later in this paper, we apply this method to computing certain 
complicated integrals that arise in the 
theory of general $L$-functions.

\section{Analytic aspects of $L$-function computations}
\label{sec:L algorithms}

\subsection{Riemann-Siegel formula}

The Riemann Siegel formula expresses the Riemann $\zeta$ function
as a main sum involving a truncated Dirichlet series and correction terms.
The formula is often presented with $\Re{s}=1/2$, but can be given
for $s$ off the critical line. See~\cite{OS} for a nice
presentation of the formula for $1/2 \leq \Re{s} \leq 2$ and references.
Here we stick to $\Re{s}=1/2$.

Let
\begin{eqnarray}
    \label{eq:Z(t)}
    Z(t) &=& e^{i\theta(t)} \zeta(1/2+it) \notag \\
    e^{i\theta(t)} &=&
    \left(
        \frac{\Gamma(1/4+it/2)}
        {\Gamma(1/4-it/2)}
    \right)^{1/2} \pi^{-it/2}.
\end{eqnarray}
The rotation factor $e^{i\theta(t)}$ is chosen so that $Z(t)$ is real.

For $t>2\pi$, let $a=(t/(2\pi))^{1/2}$, $N=\lfloor a \rfloor$, 
$\rho = \{a\}=a-\lfloor a \rfloor$ the fractional
part of $a$. Then
\begin{equation}
    \notag
    Z(t) = 
    2 \sum_{n=1}^N n^{-1/2}\cos(t\log(n)-\theta(t)) +
    R(t)
\end{equation}
where
\begin{equation}
    \notag
    R(t) =    
    \frac{(-1)^{N+1}}{a^{1/2}}
    \sum_{r=0}^m \frac{C_r(\rho)}{a^r}
    +R_m(t)
\end{equation}
with
\begin{eqnarray}
     C_0(\rho) &=& \psi (\rho) := \cos(2\pi(\rho^2-\rho-1/16))/\cos(2\pi\rho) \notag \\
     C_1(\rho) &=& -\frac{1}{96 \pi^2}\psi^{(3)}(\rho) \notag \\
     C_2(\rho) &=& \frac{1}{18432 \pi^4} \psi^{(6)}(\rho) + \frac{1}{64 \pi^2}\psi^{(2)}(\rho).
    \notag
\end{eqnarray}
In general~\cite{E}, 
$C_j(\rho)$ can be expressed as a linear combination of the derivatives of
$\psi$. We also have
\begin{equation}
    \notag
    R_m(t) = O(t^{-(2m+3)/4}).
\end{equation}
Gabcke~\cite{G} showed that
\begin{equation}
    \notag
    |R_1(t)| \leq .053 t^{-5/4}, \quad t\geq 200.
\end{equation}

The bulk of computational time in evaluating $\zeta(s)$ using the Riemann-Siegel formula
is spent on the main sum $\sum_{n=1}^N n^{-1/2}\cos(t\log(n)-\theta(t))$.
Odlyzko and Sch\"{o}nhage~\cite{OS}~\cite{O} developed an algorithm to compute the main
sum for $T \leq t \leq T+T^{1/2}$ in $O(t^\epsilon)$ operations providing that
a precomputation involving $O(T^{1/2+\epsilon})$ operations
and bits of storage are carried out beforehand. This algorithm lies behind Odlyzko's monumental
$\zeta$ computations~\cite{O}~\cite{O2}. An earlier implementation
proceeded by using the Fast Fourier Transform to compute the main sum and its derivatives at equally 
spaced grid points to then compute the main sum in between using Taylor series.
This was then improved~\cite[4.4]{O} to using 
just the values of the main sum at equally spaced points and an 
interpolation formula from the theory of band-limited functions.

Riemann used the saddle point method to obtain $C_j$, for $j\leq 5$. The reason that a nice
formula works using a sharp cutoff, truncating the sum over $n$ at $N$, is that all the
Dirichlet coefficients are equal to one. Riemann starts with an expression for
$\zeta(s)$ which involves the geometric series identity $1/(1-x) = \sum x^n$, the Taylor
coefficients on the right being the Dirichlet coefficients of $\zeta(s)$.
For general $L$-functions smoothing works better.

\subsection{Smoothed approximate functional equations}

Let
$$
   L(s) = \sum_{n=1}^{\infty} \frac{b(n)}{n^s}
$$
be a Dirichlet series that converges absolutely in a half plane, $\Re(s) > \sigma_1$,
and hence uniformly convergent in any half plane $\Re(s) \geq \sigma_2 > \sigma_1$ by
comparison with the series for $L(\sigma_2)$.

Let
\begin{equation}
    \label{eq:lambda}
    \Lambda(s) = Q^s
                 \left( \prod_{j=1}^a \Gamma(\kappa_j s + \lambda_j) \right)
                 L(s),
\end{equation}
with $Q,\kappa_j \in {\mathbb{R}}^+$, $\Re\lambda_j \geq 0$,
and assume that:
\begin{enumerate}
    \item  $\Lambda(s)$ has a meromorphic continuation to all of ${\mathbb{C}}$ with
           simple poles at $s_1,\ldots, s_\ell$ and corresponding
           residues $r_1,\ldots, r_\ell$.
    \item (functional equation)
          $\Lambda(s) = \omega \cj{\Lambda(1-\cj{s})}$ for some
          $\omega \in {\mathbb{C}}$, $\omega \neq 0$.
    \item For any $\alpha \leq \beta$, $L(\sigma +i t) = O(\exp{t^A})$ for some $A>0$,
          as $\abs{t} \to \infty$, $\alpha \leq \sigma \leq \beta$, with $A$ and the constant in
          the `Oh' notation depending on $\alpha$ and $\beta$. \label{page:condition 3}
\end{enumerate}

\begin{rems}
    a) The 3rd condition, $L(\sigma +i t) = O(\exp{t^A})$, is very mild.
    Using the fact that $L(s)$ is bounded in $\Re{s} \geq \sigma_2 > \sigma_1$,
    the functional equation and the estimate~(\ref{eq:gamma_asymptotics}), and
    the Phragm\'{e}n-Lindel\"{o}f Theorem~\cite{MR88k:00002}
    we can show that in any vertical strip $\alpha \leq \sigma \leq \beta$,
    \begin{equation}
        \notag
        L(s) = O(t^b), \quad \text{for some $b>0$}
    \end{equation}
    where both $b$ and the constant in the `Oh' notation depend
    on $\alpha$ and $\beta$.\\
    b) If $b(n), \lambda_j \in \R$, then the second assumption reads
    $\Lambda(s) = \omega \Lambda(1-s)$.\\
    c) In all known examples the $\kappa_j$'s can
    be taken to equal $1/2$. It is useful to know the Legendre duplication formula
    \begin{equation}
        \label{eq:legendre duplication}
        \Gamma(s)= (2\pi)^{-1/2}2^{s-1/2}\Gamma(s/2)\Gamma((s+1)/2).
    \end{equation}
    However, it is sometimes more convenient to work
    with~(\ref{eq:lambda}), and we avoid specializing prematurely
    to $\kappa_j=1/2$. \\
    d) The assumption that $L(s)$ have at most simple poles is not
       crucial and is only made to simplify the presentation. \\
    e) From the point of view of computing $\Lambda(s)$ given the
       Dirichlet coefficients and functional equation, we do not need 
       to assume an Euler product for $L(s)$. Without an
       Euler product, however, it is unlikely that $L(s)$ will satisfy a Riemann
       Hypothesis.
    %b) In~(\ref{eq:lambda}), $Q$ does not stand for the conductor of $L(s)$
    %   (though it is easily related to it) and $a$ does not stand for the
    %   degree of $L(s)$ (since $\kappa_j \in \cbr{1/2,1}$).
\end{rems}

To obtain a smoothed approximate functional equation with desirable
properties we introduce an auxiliary function.
Let $g: \C \to \C$ be an entire function that, for fixed $s$, satisfies
$$
    \abs{\Lambda(z+s) g(z+s) z^{-1}} \to 0
$$
as $\abs{\Im{z}} \to \infty$, in vertical strips,
$-\alpha \leq \Re{z} \leq \alpha$. The smoothed approximate functional
equation has the following form. 
\begin{theorem}
    \label{thm:formula}
    For $s \notin \cbr{s_1,\ldots, s_\ell}$, and $L(s)$, $g(s)$ as above,
    \begin{align}
         \label{eq:formula}
         \Lambda(s) g(s) =
         \sum_{k=1}^{\ell} \frac{r_k g(s_k)}{s-s_k}
         + Q^s &\sum_{n=1}^{\infty} \frac{b(n)}{n^s} f_1(s,n) \notag \\
         + \omega Q^{1-s} &\sum_{n=1}^{\infty} \frac{\cj{b}(n)}{n^{1-s}} f_2(1-s,n)
    \end{align}
    where
    \begin{align}
        \label{eq:mellin}
        f_1(s,n) &= \frac{1}{2\pi i}
                    \int_{\nu - i \infty}^{\nu + i \infty}
                    \prod_{j=1}^a \Gamma(\kappa_j (z+s) + \lambda_j)
                    z^{-1}
                    g(s+z)
                    (Q/n)^z
                    dz \notag \\
        f_2(1-s,n) &= \frac{1}{2\pi i}
                    \int_{\nu - i \infty}^{\nu + i \infty}
                    \prod_{j=1}^a \Gamma(\kappa_j (z+1-s) + \cj{\lambda_j})
                    z^{-1}
                    g(s-z)
                    (Q/n)^z
                    dz
    \end{align}
    with $\nu > \max \cbr{0,-\Re(\lambda_1/\kappa_1+s),\ldots,-\Re(\lambda_a/\kappa_a+s)}$.
\end{theorem}

%\begin{rem}
%    The third condition on $L(s)$ is not required for this Theorem. It is only
%    required if we wish to allow certain $g(s)$'s. See Section~\ref{subsection:choosing g(s)}.
%\end{rem}

\begin{proof}
Let $C$ be the rectangle with verticies $(-\alpha,-iT)$, $(\alpha,-iT)$,
$(\alpha,iT)$, $(-\alpha,iT)$,
let $s \in \C - \cbr{s_1,\ldots, s_\ell}$, and
consider
\begin{equation}
    \label{eq:lamda g}
    \frac{1}{2\pi i}
    \int_C \Lambda(z+s) g(z+s) z^{-1} dz.
\end{equation}
(integrated counter-clockwise). $\alpha$ and $T$ are chosen big
enough so that all the poles of the integrand are contained within
the rectangle. We will also require, soon, that $\alpha > \sigma_1 - \Re{s}$.
On the one hand~(\ref{eq:lamda g}) equals
\begin{equation}
    \label{eq:sum residues}
    \Lambda(s) g(s)
    + \sum_{k=1}^{\ell} \frac{r_k g(s_k)}{s_k-s}
\end{equation}
since the poles of the integrand are included in the set $\cbr{0, s_1-s,\ldots, s_\ell-s}$,
and are all simple.
Typically, the set of poles will coincide with this set. However,
if $\Lambda(s) g(s)=0$, then $z=0$ is no longer a pole of
the integrand. But then $\Lambda(s) g(s)$ contributes nothing to~(\ref{eq:sum residues})
and the equality remains valid. And if $g(s_k)=0$, then there is no pole at $z=s_k-s$
but also no contribution from $r_k g(s_k)/(s_k-s)$.

On the other hand, we may break the integral over $C$ into four integrals:
\begin{align}
    \int_C &= \int_{\alpha-iT}^{\alpha+iT} +
             \int_{\alpha+iT}^{-\alpha+iT} +
             \int_{-\alpha+iT}^{-\alpha-iT} +
             \int_{-\alpha-iT}^{\alpha-iT} \notag \\
           &= \int_{C_1} + \int_{C_2} + \int_{C_3} + \int_{C_4}. \notag
\end{align}
The integral over $C_1$, assuming that $\alpha$ is big enough to write
$L(s+z)$ in terms of its Dirichlet series i.e. $\alpha > \sigma_1 - \Re{s}$,
is
$$
    Q^s \sum_{n=1}^{\infty} \frac{b(n)}{n^s}
    \frac{1}{2\pi i} \int_{\alpha-iT}^{\alpha+iT}
    \prod_{j=1}^a \Gamma(\kappa_j (z+s) + \lambda_j)
    z^{-1} g(s+z)
    (Q/n)^z dz.
$$
We are justified in rearranging summation and integration since the series for
$L(z+s)$ converges uniformly on $C_1$.
Further, by the functional equation, the integral over $C_3$ equals
\begin{align}
    &\frac{\omega}{2\pi i}
    \int_{-\alpha+iT}^{-\alpha-iT}
    \cj{\Lambda(1-\cj{z+s})} g(z+s) z^{-1} dz \notag \\
    &=
    \omega Q^{1-s} \sum_{n=1}^{\infty} \frac{\cj{b}(n)}{n^{1-s}}
    \frac{1}{2\pi i} \int_{-\alpha+iT}^{-\alpha-iT}
    \prod_{j=1}^a \Gamma(\kappa_j (1-s-z) + \cj{\lambda_j})
    z^{-1}
    g(s+z)
    (Q/n)^{-z}
    dz \notag \\
    &=
    \omega Q^{1-s} \sum_{n=1}^{\infty} \frac{\cj{b}(n)}{n^{1-s}}
    \frac{1}{2\pi i} \int_{\alpha-iT}^{\alpha+iT}
    \prod_{j=1}^a \Gamma(\kappa_j (1-s+z) + \cj{\lambda_j})
    z^{-1}
    g(s-z)
    (Q/n)^{z}
    dz. \notag
\end{align}
Letting $T \to \infty$, the integrals over $C_2$ and $C_4$ tend to zero
by our assumption on the rate of growth of $g(s)$, and we obtain~(\ref{eq:formula}).
The integrals in~(\ref{eq:mellin}) are, by Cauchy's Theorem, independent of
the choice of $\nu$, so long as 
$\nu > \max \cbr{0,-\Re(\lambda_1/\kappa_1+s),\ldots,-\Re(\lambda_a/\kappa_a+s)}$.

%\begin{flushright}$\square$\end{flushright}
\end{proof}

\subsection{Choice of $g(z)$}
\label{sec:g(z)}

Formulae of the form~(\ref{eq:formula}) are well known~\cite{L}~\cite{Fr}. 
Usually, one finds
it in the literature with $g(s)=1$. For example, for the Riemann zeta function
this leads to Riemann's formula~\cite[pg 179]{MR91j:01070b}~\cite[pg 22]{MR88c:11049}
\begin{align}
    \pi^{-s/2} \Gamma(s/2) \zeta(s) =
    -\frac{1}{s} -\frac{1}{1-s}
    + \pi^{-s/2} &\sum_{n=1}^{\infty} \frac{1}{n^s}
    \Gamma(s/2, \pi n^2) \notag \\
    + \pi^{(s-1)/2} &\sum_{n=1}^{\infty} \frac{1}{n^{1-s}}
    \Gamma((1-s)/2, \pi n^2) \notag
\end{align}
where $\Gamma(s,w)$ is the incomplete gamma function (see
Section~\ref{subsection:fi(s,n),a=1}).

However, the choice $g(s)=1$ is not well suited for computing $\Lambda(s)$
as $\abs{\Im(s)}$ grows. By Stirling's
formula~\cite[pg 294]{MR55:8655}
\begin{equation}
    \label{eq:gamma_asymptotics}
    \left| \Gamma(s) \right|
       \sim (2\pi)^{1/2}|s|^{\sigma-1/2}e^{-|t|\pi/2}, \quad s=\sigma +i t
\end{equation}
as $|t| \rightarrow \infty$, and
so decreases very quickly as $|t|$ increases. Hence, with $g(s)=1$, the
\lhs of~(\ref{eq:formula}) is extremely small for large $\abs{t}$ and fixed
$\sigma$.
On the other hand, we can show that the terms on the r.h.s., though
decreasing as $n \to \infty$, start off relatively large compared
to the l.h.s.. Hence a tremendous amount of cancellation must occur on the
\rhs and and we would need an unreasonable amount of precision.
This problem is analogous to what happens if we try to sum 
$\exp(-x) = \sum (-x)^n/n!$ in a naive way.
If $x$ is positive and large, the \lhs is exponentially small, yet the terms on the
\rhs are large before they become small and high precision is needed to capture
the ensuing cancellation.

One way to control this cancellation is to choose $g(s)$ equal to
$\delta^{-s}$ with $|\delta|=1$ and chosen to cancel out most of the
exponentially small size of the $\Gamma$ factors. This idea appears
in the work of Lavrik~\cite{L}, and was also suggested
by Lagarias and Odlyzko~\cite{MR80g:12010} who did not implement it
since it led to complications regarding the computation of~(\ref{eq:mellin}). 
This method was successfully applied in the author's PhD thesis~\cite{Ru}
to compute Dirichlet $L$-functions and $L$-functions associated to cusp forms
and is used extensively in the author's $L$-function package~\cite{Ru3}
More recently, this approach was used in the computation of 
Akiyama and Tanigawa~\cite{AT} to compute several elliptic curve $L$-functions.

In fact when there are multiple
$\Gamma$ factors it is better to choose
a different $\delta$ for each $\Gamma$ and multiply these together.
For a given $s$ let
\begin{eqnarray}
    \label{eq:delta}
    &&t_j = \Im(\kappa_j s + \lambda_j) \notag \\
    &&\theta_j = 
    \begin{cases}
        \pi/2, \quad &\text{if $\abs{t_j} \leq 2c/(a\pi)$} \\
        c/(a |t_j|), \quad &\text{if $\abs{t_j} > 2c/(a\pi)$}
    \end{cases} \notag \\
    &&\delta_j = \exp(i\ \sgn(t_j)(\pi/2-\theta_j)).
\end{eqnarray}
Here $c>0$ is a free parameter. Larger $c$ means faster convergence 
of the sums in~(\ref{eq:formula}), but also more cancellation and loss of precision.

Next, we set
\begin{equation}
    \label{eq:g(s) 2}
    g(z) := \prod_{j=1}^a \delta_j^{-\kappa_j z - \Im{\lambda_j}} =  \beta \delta^{-z}.
\end{equation}
Because $\delta_j$ depends on $s$, the constants $\delta$ and $\beta$ depends on $s$.
We can either use a fresh $\delta$ for each new $s$ value, or else
modify the above choice of $t_j$ so as to use the same $t_j$ 
for other nearby $s$'s. The latter is prefered if we wish to
carry out precomputations that can be recycled as we vary $s$.
For simplicity, here we assume that a fresh $\delta$ is chosen as 
above for each new $s$.

The choice of $g$ controls the exponentially small size of
the $\Gamma$ factors. Notice that the constant factor 
$\beta = \prod_{j=1}^a \delta_j^{-\Im{\lambda_j}}$ in 
(\ref{eq:g(s) 2}) appears in every term in~(\ref{eq:formula}),
and hence can be dropped from $g(z)$ without any effect on cancellation
or the needed precision. However, to analyze the size of the 
the l.h.s. of~(\ref{eq:formula}) and the terms on the r.h.s. this
factor is helpful and we leave it in for now, but with the understanding
that it can be omitted. 

To see the effect of the function $g(z)$ on the l.h.s of~(\ref{eq:formula}) 
we have, by~(\ref{eq:gamma_asymptotics})
and~(\ref{eq:delta}) 
\begin{eqnarray}
    \abs{\Lambda(s_0)g(s_0)} &\sim&
        \ast \cdot 
        \abs{L(s_0)}
        \prod_{|t_j| \leq 2c/\pi} \exp\br{-\abs{t_j}\pi/2} 
        \prod_{|t_j| > 2c/\pi} \exp\br{-c/a} \notag \\
        &\geq& \ast \cdot
        \abs{L(s_0)}
        \exp(-c)
        \notag
\end{eqnarray}
where
\begin{equation}
    \notag
    \ast = Q^{\sigma_0} (2\pi)^{a/2}
           \prod_{j=1}^a \abs{\kappa_j s_0 + \lambda_j}^{\kappa_j \sigma_0 + \Re{\lambda_j}-1/2}.
\end{equation}
We have thus managed to control the exponentially small size of $\Lambda(s)$
up to a factor of $\exp(-c)$ which we can regulate via the choice of $c$.
We can also show that this choice of $g(z)$ leads to well balanced terms on
the r.h.s. of~(\ref{eq:formula}).

\subsection{Approximate functional equation in the case of one $\Gamma$-factor}
\label{subsection:fi(s,n),a=1}

We first treat the case $a=1$ separately because it is the simplest,
the greatest number of tools have been developed to handle this case,
and many popular $L$-functions have $a=1$.

Here we are assuming that
\begin{equation}
    \notag
    \Lambda(s) = Q^s \Gamma(\gamma s + \lambda) L(s).
\end{equation}

According to~(\ref{eq:g(s) 2}) we should set
$$
     g(s)=\delta^{-s}
$$
(we omit the factor $\beta$ as described following~(\ref{eq:g(s) 2}))
with 
$$
    \delta = {\delta_1}^\gamma
$$
and
\begin{eqnarray}
    &&t_1 = \Im(\gamma s + \lambda) \notag \\
    &&\theta_1 =
    \begin{cases}
        \pi/2, \quad &\text{if $\abs{t_1} \leq 2c/\pi$} \\
        c/|t_1|, \quad &\text{if $\abs{t_1} > 2c/\pi$}
    \end{cases} \notag \\
    &&\delta_1 = \exp(i\ \sgn(t_1)(\pi/2-\theta_1)).
    \notag
\end{eqnarray}

In that case, the function $f_1(s,n)$ that appears in
Theorem~\ref{thm:formula} equals
\begin{align}
    f_1(s,n) &= \frac{\delta^{-s}}{2\pi i}
                \int_{\nu - i \infty}^{\nu + i \infty}
                \Gamma(\gamma (z+s) + \lambda) z^{-1}
                \br{Q/(n \delta)}^z
                dz \notag \\
             &= \frac{\delta^{-s}}{2\pi i}
                \int_{\gamma \nu - i \infty}^{\gamma \nu + i \infty}
                \Gamma(u+\gamma s + \lambda) u^{-1}
                \br{Q/(n \delta)}^{u/\gamma}
                du.\notag
\end{align}
Now 
\begin{equation}
    \label{eq:invmellin g(v+u)/u}
    \Gamma(v+u) u^{-1} = \int_0^\infty \Gamma(v,t) t^{u-1} dt,
    \quad \Re{u}>0, \quad \Re(v+u)>0
\end{equation}
where
\begin{eqnarray}
    \notag
    \Gamma(z,w) =\int_w^\infty e^{-x} x^{z-1} dx \quad |\arg{w}|<\pi  \notag \\
    =w^{z} \int_1^\infty e^{-wx}x^{z-1} dx, \quad \Re(w)>0. \notag
\end{eqnarray}
$\Gamma(z,w)$ is known as the incomplete gamma function.
By Mellin inversion
\begin{equation}
    \notag
    f_1(s,n) = \delta^{-s} \Gamma\br{\gamma s +\lambda,\br{n \delta/Q}^{1/\gamma}}.
\end{equation}
Similarly
\begin{equation}
    \notag
    f_2(1-s,n) = \delta^{-s} \Gamma\br{\gamma (1-s) +\cj{\lambda},\br{n/(\delta Q)}^{1/\gamma}}.
\end{equation}
%\enlargethispage*{18ex}
We may thus express, when $a=1$ and $g(s)=\delta^{-s}$,~(\ref{eq:formula}) as
\begin{eqnarray}
%\boxed{
%\begin{align}
     Q^s \Gamma(\gamma s +\lambda) L(s)
     \delta^{-s} =
     &&\sum_{k=1}^{\ell} \frac{r_k \delta^{-s_k}}{s-s_k} \notag \\
     +&&\br{\delta/Q}^{\lambda/\gamma}
      \sum_{n=1}^{\infty}
      b(n) n^{\lambda/\gamma} 
      G\br{\gamma s +\lambda, \br{n \delta/Q}^{1/\gamma}}
      \notag \\
     +&&\frac{\omega}{\delta} (Q\delta)^{-\cj{\lambda}/\gamma}
      \sum_{n=1}^{\infty}
      \cj{b}(n)n^{\cj{\lambda}/\gamma}
      G\br{\gamma (1-s) +\cj{\lambda}, \br{n/(\delta Q)}^{1/\gamma}}
      \notag
%\end{align}
%}
\end{eqnarray}
{\vspace{-4ex} \begin{equation} \label{eq:formula a=1} \end{equation}}
where
\begin{equation}
    \label{eq:G(z,w)}
    G(z,w) = w^{-z} \Gamma(z,w) = \int_1^\infty e^{-wx} x^{z-1} dx,
    \quad \Re(w)>0.
\end{equation}
Note, from~(\ref{eq:g(s) 2}) with $a=1$, we have $\Re \delta^{1/\gamma}>0$,
so both $\br{n \delta/Q}^{1/\gamma}$ and $\br{n/(\delta Q)}^{1/\gamma}$
have positive $\Re$ part.

%-------------------------------------------------------
%\begin{flushleft}
%\hrulefill
%\end{flushleft}

\subsubsection{Examples}

\begin{flushleft}
1) Riemann zeta function, $\zeta(s)$:
the necessary background can be found in~\cite{MR88c:11049}.
Formula~(\ref{eq:formula a=1}), for $\zeta(s)$, is
\begin{align}
    \label{eq:zeta}
    \pi^{-s/2} \Gamma(s/2) \zeta(s)
    \delta^{-s} =
    -\frac{1}{s} - \frac{\delta^{-1}}{1-s}
    +&\sum_{n=1}^{\infty}
    G\br{s/2, \pi n^2 \delta^2 }
    \notag \\
    +\delta^{-1}
    &\sum_{n=1}^{\infty}
    G\br{(1-s)/2, \pi n^2/\delta^2 }
\end{align}
\end{flushleft}

%-------------------------------------------------------
\begin{flushleft}
2) Dirichlet $L$-functions, $L(s,\chi)$:
(see~\cite[chapter 9]{D}). When $\chi$ is primitive and even, $\chi(-1)=1$, we get
\begin{align}
    \br{\frac{q}{\pi}}^{s/2} \Gamma(s/2) L(s,\chi)
    \delta^{-s} =
    &\sum_{n=1}^{\infty}
    \chi(n) G\br{s/2, \pi n^2 \delta^2 /q}
    \notag \\
    +\frac{\tau(\chi)}{\delta q^{1/2}}
    &\sum_{n=1}^{\infty}
    \cj{\chi}(n) G\br{(1-s)/2, \pi n^2 /(\delta^2 q) }
    \notag
\end{align}
and when  $\chi$ is primitive and odd, $\chi(-1)=-1$, we get
\begin{align}
    \br{\frac{q}{\pi}}^{s/2} \Gamma(s/2+1/2) L(s,\chi)
    \delta^{-s} =
    \delta \br{\frac{\pi}{q}}^{1/2}
    &\sum_{n=1}^{\infty}
    \chi(n) n G\br{s/2+1/2, \pi n^2 \delta^2 /q}
    \notag \\
    +\frac{\tau(\chi) \pi^{1/2}}{i q \delta^2}
    &\sum_{n=1}^{\infty}
    \cj{\chi}(n) n G\br{(1-s)/2 +1/2, \pi n^2 /(\delta^2 q) }
    \notag
\end{align}
Here, $\tau(\chi)$ is the Gauss sum
$$
   \tau(\chi) = \sum_{m=1}^q \chi(m) e^{2\pi i m/q}.
$$
\end{flushleft}

%-------------------------------------------------------
\begin{flushleft}
3) Cusp form $L$-functions: (see~\cite{MR41:1648}).
Let $f(z)$ be a cusp form of weight $k$ for SL$_2(\Z)$, $k$ a positive
even integer:
\end{flushleft}
\begin{enumerate}
   \item $f(z)$ is entire on $\UH$, the upper half plane.
   \item $f(\sigma z) = (cz+d)^k f(z)$, $\sigma=\begin{pmatrix} a & b
         \\ c & d \end{pmatrix} \in $ SL$_2(\Z)$, $z \in \UH$.
   \item $\lim_{t \to \infty} f(it)=0$.
\end{enumerate}
Assume further that $f$ is a Hecke eigenform, i.e. an eigenfunction
of the Hecke operators.
We may expand $f$ in a Fourier series
$$
   f(z) = \sum_{n=1}^\infty a_n e^{2\pi i n z},
   \quad \Im(z)>0
$$
and associate to $f(z)$ the Dirichlet series
$$
    L_f(s) := \sum_1^\infty \frac{a_n}{n^{(k-1)/2}}n^{-s}.
$$
We normalize $f$ so that $a_1=1$.
This series converges absolutely when $\Re(s) > 1$ because, as
proven by Deligne~\cite{MR49:5013},
\begin{equation}
    \notag
    |a_n| \leq \sigma_0(n) n^{(k-1)/2},
\end{equation}
where $\sigma_0(n) := \sum_{d|n}1 = O(n^\epsilon)$ for any
$\epsilon>0$.

$L_f(s)$ admits an analytic continuation to all of $\C$
and satisfies the functional equation
$$
    \Lambda_f(s) := (2\pi)^{-s} \Gamma(s+(k-1)/2) L_f(s) =
    (-1)^{k/2} \Lambda_f(1-s).
$$
With our normalization, $a_1=1$, the $a_n$'s are real since they are
eigenvalues of self adjoint operators, the Hecke operators with respect
to the Petersson inner product
(see~\cite[III-12]{MR41:1648}).
Furthermore, the required rate of growth on $L_f(s)$, condition 3 on
page~\pageref{page:condition 3}, follows from the modularity of $f$.

Hence, in this example, formula~(\ref{eq:formula a=1}) is
\begin{align}
    (2 \pi)^{-s} \Gamma(s+(k-1)/2) L_f(s)
    \delta^{-s} =
    \br{\delta 2 \pi}^{(k-1)/2}
    &\sum_{n=1}^{\infty}
    a_n G\br{s+(k-1)/2, 2 \pi n \delta}
    \notag \\
    +\frac{(-1)^{k/2}}{\delta}
    \br{\frac{2\pi}{\delta}}^{(k-1)/2}
    &\sum_{n=1}^{\infty}
    a_n G\br{1-s+(k-1)/2, 2 \pi n/\delta}
    \notag
\end{align}

%-------------------------------------------------------

\begin{flushleft}
4) Twists of cusp forms:
$L_f(s,\chi)$, $\chi$ primitive, $f(z)$ as in the previous example. 
$L_f(s,\chi)$ is given by the Dirichlet series
\end{flushleft}
$$
    L_f(s,\chi) = \sum_1^\infty \frac{a_n \chi(n)}{n^{(k-1)/2}} n^{-s}.
$$
$L_f(s,\chi)$ extends to an entire function and satisfies the
functional equation
\begin{eqnarray}
    \Lambda_f(s,\chi) &:=& 
    \br{\frac{q}{2\pi}}^{s} \Gamma(s+(k-1)/2) L_f(s,\chi) \notag \\
    &=&
    (-1)^{k/2} \chi(-1) \frac{\tau(\chi)}{\tau(\cj{\chi})} \Lambda_f(1-s,\cj{\chi}).
    \notag
\end{eqnarray}
In this example, formula~(\ref{eq:formula a=1}) is
\begin{eqnarray}
    \notag
    &&\br{\frac{q}{2\pi}}^{s} \Gamma(s+(k-1)/2) L_f(s,\chi)
    \delta^{-s} = \notag \\
    &&\br{\frac{2 \pi \delta}{q}}^{(k-1)/2}
    \sum_{n=1}^{\infty}
    a_n \chi(n) G\br{s+(k-1)/2, 2 \pi n \delta /q}
    \notag \\
    &&+\frac{(-1)^{k/2}}{\delta}
    \chi(-1) 
    \frac{\tau(\chi)}{\tau(\cj{\chi})}
    \br{\frac{2\pi}{q \delta}}^{(k-1)/2}
    \sum_{n=1}^{\infty}
    a_n \cj{\chi}(n) G\br{1-s+(k-1)/2, 2 \pi n /(\delta q)}. \notag
\end{eqnarray}

%-------------------------------------------------------
\begin{flushleft}
5) Elliptic curve $L$-functions: (see~\cite[especially chapters X,XII]{MR93j:11032}).
Let $E$ be an elliptic curve over $\Q$, which we write in global
minimal Weierstrass form
\end{flushleft}
$$
   y^2 + c_1xy + c_3y = x^3 +c_2x^2 +c_4x +c_6
$$
where the $c_j$'s are integers and the disciminant $\Delta$ is minimal.

To the elliptic curve $E$ we may associate an Euler product
\begin{equation}
   \label{eq:euler_prdct}
   L_E(s) := \prod_{p|\Delta} (1 - a_p p^{-1/2-s})^{-1}
             \prod_{p \nmid \Delta}(1- a_p p^{-1/2-s}+p^{-2s})^{-1}
\end{equation}
where, for $p \nmid \Delta$,
$a_p = p+1-\# E_p(\Z_p)$, with $\# E_p(\Z_p)$ being the number
of points $(x,y)$ in $\Z_p \times \Z_p$ on the curve $E$ considered
modulo $p$, together with the point at infinity. 
When $p|\Delta$, $a_p$ is either $1$, $-1$, or $0$.
If $p \nmid \Delta$, a theorem of Hasse states that $\abs{a_p} < 2 p^{1/2}$.
Hence,~(\ref{eq:euler_prdct}) converges when $\Re(s)>1$, and for
these values of $s$ we may expand $L_E(s)$ in an absolutely convergent
Dirichlet series
\begin{equation}
    \label{eq:dirichlet series E}
    L_E(s)=\sum_1^\infty \frac{a_n}{n^{1/2}} n^{-s}.
\end{equation}

The Hasse-Weil conjecture asserts that $L_E(s)$ extends to an entire
function and has the functional equation
$$
    \Lambda_E(s) :=
    \br{\frac{N^{1/2}}{2\pi}}^{s} \Gamma(s+1/2) L_E(s) =
    -\varepsilon \Lambda_E(1-s).
$$
where $N$ is the conductor of $E$, and $\varepsilon$, which depends
on $E$, is either $\pm 1$. The Hasse-Weil conjecture and also the required
rate of growth on $L_E(s)$  follows from
the Shimura-Taniyama-Weil conjecture, which has been
proven by Wiles and Taylor~\cite{MR96d:11072}~\cite{MR96d:11071} for elliptic curves 
with square free conductor and has been extended, 
by  Breuil, Conrad, Diamond and Taylor to all elliptic curves over $\Q$
\cite{BCDT}.

Hence  we have
\begin{align}
    \br{\frac{N^{1/2}}{2\pi}}^{s} \Gamma(s+1/2) L_E(s)
    \delta^{-s} =
    \br{\frac{2 \pi \delta}{N^{1/2}}}^{1/2}
    &\sum_{n=1}^{\infty}
    a_n G\br{s+1/2, 2 \pi n \delta / N^{1/2}}
    \notag \\
    -\frac{\varepsilon}{\delta}
    \br{\frac{2\pi}{N^{1/2} \delta}}^{1/2}
    &\sum_{n=1}^{\infty}
    a_n G\br{1-s+1/2, 2 \pi n /(\delta N^{1/2})}.
    \notag
\end{align}

%-----------------------------------------------------
\begin{flushleft}
6) Twists of elliptic curve $L$-functions:
$L_E(s,\chi)$, $\chi$ a primitive character of conductor $q$, $(q,N)=1$.
Here $L_E(s,\chi)$ is given by the Dirichlet series
\end{flushleft}
$$
    L_E(s,\chi) = \sum_1^\infty \frac{a_n}{n^{1/2}} \chi(n) n^{-s}.
$$
The Weil conjecture asserts, here, that $L_E(s)$ extends to
an entire function and satisfies
$$
    \Lambda_E(s,\chi) :=
    \br{\frac{q N^{1/2}}{2\pi}}^{s} \Gamma(s+1/2) L_E(s,\chi) =
    -\varepsilon \chi(-N) \frac{\tau(\chi)}{\tau(\cj{\chi})} \Lambda_E(1-s,\cj{\chi}).
$$
Here $N$ and $\varepsilon$ are the same as for $E$.
In this example the conjectured formula is
\begin{align}
    \br{\frac{qN^{1/2}}{2\pi}}^{s} \Gamma(s+1/2) L_E(s)
    \delta^{-s} =
    \br{\frac{2 \pi \delta}{q N^{1/2}}}^{1/2}
    &\sum_{n=1}^{\infty}
    a_n \chi(n) G\br{s+1/2, 2 \pi n \delta / (q N^{1/2})}
    \notag \\
    -\frac{\varepsilon}{\delta}
    \chi(-N) \frac{\tau(\chi)}{\tau(\cj{\chi})}
    \br{\frac{2\pi}{q N^{1/2} \delta}}^{1/2}
    &\sum_{n=1}^{\infty}
    a_n \cj{\chi}(n) G\br{1-s+1/2, 2 \pi n /(\delta q N^{1/2})}.
    \notag
\end{align}
%------------------------------------------------
\begin{flushleft}
\hrulefill
\end{flushleft}

We have reduced in the case $a=1$ the computation of $\Lambda(s)$ to one of evaluating
two sums of incomplete gamma functions. The $\Gamma(\gamma s +\lambda) \delta^{-s}$
factor on the left of~(\ref{eq:formula a=1}) and elsewhere is easily evaluated
using several terms of Stirling's asymptotic formula 
and also the
recurrence $\Gamma(z+1) = z \Gamma(z)$ applied a few times.
The second step is needed for small $z$. Some care needs to be taken to
absorb the $e^{-\pi\abs{\Im(\gamma s +\lambda)}/2}$ factor of $\Gamma(\gamma s +\lambda)$ into 
the $e^{\pi \abs{\Im(\gamma s +\lambda)}/2}$
factor of $\delta^{-s}$. Otherwise our effort to control the size of
$\Gamma(\gamma s +\lambda)$ will have been in vain, and lack of precision 
will wreak havoc.

To see how many terms in (\ref{eq:formula a=1}) are needed 
we can use the rough bound
\begin{equation}
  \notag
  |G(z,w)| < e^{-\Re(w)} \int_0^\infty e^{-(\Re(w)-\Re(z)+1)t} dt
       = \frac{e^{-\Re(w)}}{\Re(w)-\Re(z)+1},
\end{equation}
valid for $\Re(w) > \Re(z)-1>0$. 
We have put $t=x-1$ in~(\ref{eq:G(z,w)}) and have used $t+1 \leq e^t$.
Also, for $\Re(w)>0$ and $\Re(z) \leq 1$, 
\begin{equation}
  \notag
  |G(z,w)| < \frac{e^{-\Re(w)}}{\Re(w)}.
\end{equation}
These inequalities tells us that the terms in~(\ref{eq:formula a=1}) decrease exponentially fast
once $n$ is sufficiently large.

For example, in equation (\ref{eq:zeta}) for $\zeta(s)$
we get exponential drop off roughly when
$$
    \Re{\pi n^2\delta^2} >> 1.
$$
But 
$$
    \Re{\pi n^2 \delta^2} = \pi n^2 \Re{\delta^2} \sim 2 \pi n^2 c/t
$$
so the number of terms needed is roughly
$$
     >> (t/c)^{1/2}.
$$

\subsubsection{Computing $\Gamma(z,w)$}

Recall the definitions
\begin{eqnarray}
    \Gamma(z,w) &=& \int_w^\infty e^{-t} t^{z-1} dt, \quad |\arg{w}|<\pi \notag \\
    G(z,w) &=& w^{-z} \Gamma(z,w). \notag 
\end{eqnarray}
Let
$$
   \gamma(z,w) := \Gamma(z)-\Gamma(z,w) = \int_0^w e^{-x} x^{z-1} dx,
   \quad \Re{z}>0, \quad \abs{\arg{w}} < \pi
$$
be the complimentary incomplete gamma function, and set
\begin{equation}
    \label{eq:g defn}
    g(z,w) = w^{-z} \gamma(z,w) = \int_0^1 e^{-wt}t^{z-1} dt
\end{equation}
so that $G(z,w)+g(z,w)=w^{-z}\Gamma(z)$. 
The function $g(z,w)/\Gamma(z)$ is entire in $z$ and $w$.

The incomplete $\Gamma$ function undergoes a transition when
$|w|$ is close to $|z|$. This will be described using Temme's uniform asymptotics
for $\Gamma(z,w)$. The transition explains the difficulty in computing
$\Gamma(z,w)$ without resorting to several different expressions
or using uniform asymptotics.

A combination of series, asymptotics, and continued
fractions are useful when $|z|$ is somewhat bigger than or smaller than $|w|$.
When the two parameters are close in size to one another, we can employ Temme's more
involved uniform asymptotics. We can also apply the Poisson summation method described
in Section 2, or an expansion due to Nielsen. Below we look at a few useful approaches.

Integrating by parts we get
\begin{equation}
   \notag
   g(z,w) =
    e^{-w}\sum_{j=0}^{\infty} \frac{w^j}{(z)_{j+1}}
\end{equation}
where
\begin{equation}
   \notag
   (z)_j = \begin{cases}
                  z(z+1)\ldots(z+j-1)  & \text{if $j>0$;} \\
                  1                    & \text{if $j=0$.}
           \end{cases}
\end{equation}
(The case $j=0$ occurs below in an expression for $G(z,w)$).
While this series converges for $z\neq 0,-1,-2,\ldots$ and all $w$, 
it is well suited, say if $\Re{z} >0$ and $|w|< \alpha |z|$ 
with $0<\alpha <1$. Otherwise, not only does the series
take too long to converge, but precision issues arise.

The following continued fraction converges for
$\Re{z} >0$ 
\begin{equation}
    \notag
    g(z,w) = \cfrac{e^{-w}}
             {z-\cfrac{zw}
                     {z+1+\cfrac{w}{z+2-\cfrac{(z+1)w}{z+3+\cfrac{2w}{z+4-\cfrac{(z+2)w}{z+5+\dotsb}}}}}
             }
\end{equation}
The paper of Akiyama and Tanigawa~\cite{AT} contains an analysis of the truncation error for
this continued fraction, as well as the continued fraction in (\ref{eq:cfrac2}) below,
and show that the above is most useful when $|w| < |z|$, with poorer performance
as $|w|$ approaches $|z|$.

Another series, useful when $|w|<<1$, is 
\begin{equation}
   \notag
   g(z,w) =
   \sum_{j=0}^{\infty} \frac{(-1)^j}{j!}\frac{w^j}{z+j}.
\end{equation}
This is obtained from (\ref{eq:g defn}) by expanding $e^{-wt}$ in a Taylor
series and integrating termwise. As $|w|$ grows, cancellation and precision become
an issue in the same way it does for the sum $e^{-w} = \sum(-w)^j/j!$.

Next, integrate $G(z,w)$ by parts to obtain the asymptotic series
\begin{equation}
   \notag
   G(z,w) = \frac{e^{-w}}{w}
          \sum_{j=0}^{M-1} \frac{(1-z)_j}{(-w)^j} + \epsilon_M(z,w)
\end{equation}
with
\begin{equation}
   \notag
   \epsilon_M(z,w) = \frac{(1-z)_M}{(-w)^M} G(z-M,w). 
\end{equation}
This asymptotic expansion works well if $|w| > \beta |z|$ with $\beta>1$
and $|z|$ large. In that region the following continued fraction
also works well
\begin{equation}
    \label{eq:cfrac2}
    G(z,w) = \cfrac{e^{-w}}
             {w+\cfrac{1-z}
                     {1+\cfrac{1}{w+\cfrac{2-z}{1+\cfrac{2}{w+\cfrac{3-z}{1+\cfrac{3}{w+\dotsb}}}}}}
             }
\end{equation}

Temme's uniform asymptotics for $\Gamma(z,w)$ provide a powerful tool for
computing the function in its transition zone and elsewhere.
Following the
notation in~\cite{T}, let
\begin{eqnarray}
    Q(z,w) &=& \Gamma(z,w)/\Gamma(z) \notag \\
    \lambda &=& w/z \notag \\
    \eta^2/2 &=& \lambda -1 -\log{\lambda}
    \notag
\end{eqnarray}
where the sign of $\eta$ is chosen to be positive for $\lambda >1$.
Then
\begin{equation}
   \notag
    Q(z,w) = \frac{1}{2} \text{erfc}(\eta (z/2)^{1/2}) + R_z(\eta)
\end{equation}
where
\begin{equation}
   \notag
    \text{erfc} = \frac{2}{\pi^{1/2}} \int_z^\infty e^{-t^2} dt,
\end{equation}
and $R_z$ is given by the asymptotic series, as $z \to \infty$,
\begin{equation}
    \label{eq:R_z}
    R_z(\eta)
    = \frac{e^{-z\eta^2 /2}}{(2\pi z)^{1/2}} 
    \sum_{n=0}^\infty \frac{c_n(\eta)}{z^n}.
\end{equation}
Here
\begin{eqnarray}
    c_0(\eta) &=& \frac{1}{\lambda-1} -\frac{1}{\eta} \notag \\
    c_1(\eta) &=& \frac{1}{\eta^3} 
                  - \frac{1}{(\lambda-1)^3} -\frac{1}{(\lambda-1)^2} -\frac{1}{12(\lambda-1)} \notag \\
    \eta c_n(\eta) &=& \frac{d}{d\eta} c_{n-1}(\eta) +\frac{\eta}{\lambda-1}  \gamma_n, \quad n \geq 1
   \notag
\end{eqnarray}
with
\begin{equation}
   \notag
    \Gamma^*(z) = \sum_0^\infty \frac{(-1)^n\gamma_n}{z^n}
\end{equation}
being the asymptotic expansion of 
$$
    \Gamma^*(z) = (z/(2\pi))^{1/2} (e/z)^z \Gamma(z).
$$
The first few terms are $\gamma_0=1$, $\gamma_1=-1/12$, $\gamma_2=1/288$, $\gamma_3=139/51840$.
The singularities at $\eta=0$, i.e. $\lambda=1, z=w$, are removable.
Unfortunately, explicit
estimates for the remainder in truncating (\ref{eq:R_z})
when the parameters are complex
have not been worked out, but in practice the expansion seems to work very well.

To handle the intermediate region $|z| \approx |w|$ we could also 
use the following expansion
of Nielsen to step through the troublesome region
\begin{equation}
    \label{eq:E2 a=1}
    \gamma(z,w+d) = \gamma(z,w) +
    w^{z-1} e^{-w} \sum_{j=0}^\infty \frac{(1-z)_j}{(-w)^j} (1-e^{-d}e_j(d)),
    \quad  \abs{d} < \abs{w} \\
\end{equation}
where
$$
   e_j(d) = \sum_{m=0}^j \frac{d^m}{m!}.
$$
A proof can be found in~\cite{MR84h:33001b}. 
This expansion is very well suited, for example, for $L$-functions associated
to modular forms, since in that case we increment $w$ in equal steps from term to
term in (\ref{eq:formula a=1}) and precomputations can be arranged to
recycle data. Numerically, this expansion is unstable if $|d|$ is big. 
This can be overcome by taking many smaller steps, but this then makes Nielsen's expansion 
an inefficient choice for $\zeta(s)$ or Dirichlet $L$-functions.

In computing~(\ref{eq:E2 a=1}) some care needs to be taken to avoid numerical
pitfalls. One pitfall is that, as $j$ grows, $e^{-d}e_j(d) \to 1$. So once 
$\abs{1-e^{-d}e_j(d)}<10^{-\text{Digits}}$, the error in computation of $1-e^{-d}e_j(d)$
is bigger than its value, and this gets magnified when we multiply by 
$(1-z)_j/(-w)^j$. So in computing
$\br{(1-z)_j/(-w)^j}(1-e^{-d}e_j(d))$ one must avoid the temptation to view this
as a product of $(1-z)_j/(-w)^j$ and $1-e^{-d}e_j(d)$. Instead, we let
$$
    a_j(z,w,d)=  \frac{(1-z)_j}{(-w)^j} (1-e^{-d}e_j(d)).
$$
Now, $1-e^{-d}e_j(d)= e^{-d}(e^d-e_j(d))$, and we get
\begin{align}
    a_{j+1}(z,w,d) &= a_j(z,w,d)
                      \frac{z-(j+1)}{w}
                      \br{\sum_{j+2}^\infty d^m/m!} \left/
                      \br{\sum_{j+1}^\infty d^m/m!} \right. \notag \\
                   &= a_j(z,w,d)
                      \frac{z-(j+1)}{w}
                      \br{1-1/\beta_j(d)}, \quad j=1,2,3,\ldots \notag
\end{align} 
where
\begin{equation}
   \notag
   \beta_j(d) = \sum_{m=0}^\infty d^m/(j+2)_m. 
\end{equation}
Furthermore
$$
    \beta_j(d) - 1 \sim d/(j+2), \quad \text{as $\abs{d}/j \to 0$}.
$$
Hence, for $\abs{w} \approx \abs{z}$, we {\it approximately} have (as $\abs{d}/j \to 0$)
$$
    \abs{
        \frac{z-(j+1)}{w}
        \br{1-1/\beta_j(d)}
    } \leq
    \br{1+\frac{j+1}{\abs{w}}}\frac{\abs{d}}{j+2} \leq 
    \frac{\abs{d}}{j+2} +\frac{\abs{d}}{\abs{w}}.
$$
Thus, because $\abs{d/w}<1$, we have, for $j$ big enough, that the above is $<1$,
and so the sum in~(\ref{eq:E2 a=1}) converges geometrically fast, and hence only
a handful of terms are required. 

One might be tempted to compute the $\beta_j(d)$'s using the recursion
$$
   \beta_{j+1}(d) = (\beta_j(d)-1)(j+2)/d
$$
but this leads to numerical instability. The $\beta_j(d)$'s are all equal to
$1+O_d(1/(j+2))$ and are thus all roughly of comparable size. Hence, a small
error due to roundoff in $\beta_j(d)$ is turned into a much larger error in
$\beta_{j+1}(d)$, $(j+2)/\abs{d}$ times larger, and this quickly 
destroys the numerics.

There seems to be some potential in an
asymptotic expression due to Ramanujan~\cite[pg 193, entry 6]{MR90b:01039}
$$
    G(z,w) \sim
    w^{-z} \Gamma(z)/2
    +e^{-w} \sum_{k=0}^M p_k(w-z+1)/w^{k+1},
    \quad \text{as $\abs{z} \to \infty$},
$$
for $\abs{w-z}$ relatively small,
where $p_k(v)$ is a polynomial in $v$ of degree $2k+1$, though this
potential has not been investigated substantially. 

We list the first few $p_k(v)$'s here:
\begin{align}
    p_0(v)   =
    &-v+2/3 \notag \\
    p_1(v)   =
    &-{\frac {{v}^{3}}{3}}+{\frac {{v}^{2}}{3}}-{\frac {4}{135}}
    \notag \\
    p_2(v)   =
    &-{\frac {{v}^{5}}{15}}+{\frac {{v}^{3}}{9}}-{\frac {2\,{v}^{2}}{135}}-
    {\frac {4\,v}{135}}+{\frac {8}{2835}}
    \notag \\
    p_3(v)   =
    &-{\frac {{v}^{7}}{105}}-{\frac {{v}^{6}}{45}}+{\frac {{v}^{5}}{45}}+{
    \frac {7\,{v}^{4}}{135}}-{\frac {8\,{v}^{3}}{405}}-{\frac {16\,{v}^{2}
    }{567}}+{\frac {16\,v}{2835}}+{\frac {16}{8505}}
    \notag \\
    p_4(v)   =
    &-{\frac {{v}^{9}}{945}}-{\frac {2\,{v}^{8}}{315}}-{\frac {2\,{v}^{7}}{
    315}}+{\frac {8\,{v}^{6}}{405}}+{\frac {11\,{v}^{5}}{405}}-{\frac {62
    \,{v}^{4}}{2835}}-{\frac {32\,{v}^{3}}{1215}}+{\frac {16\,{v}^{2}}{
    1701}}+{\frac {16\,v}{2835}}-{\frac {8992}{12629925}}
    \notag \\
    p_5(v)   =
    &-{\frac {{v}^{11}}{10395}}-{\frac {{v}^{10}}{945}}-{\frac {2\,{v}^{9}}
    {567}}-{\frac {2\,{v}^{8}}{2835}}+{\frac {43\,{v}^{7}}{2835}}+{\frac {
    41\,{v}^{6}}{2835}}-{\frac {968\,{v}^{5}}{42525}}-{\frac {68\,{v}^{4}}
    {2835}}+{\frac {368\,{v}^{3}}{25515}} \notag \\
    &+{\frac {138064\,{v}^{2}}{
    12629925}}-{\frac {35968\,v}{12629925}}-{\frac {334144}{492567075}}
    \notag
\end{align}

It is worth noting that when many evaluations of $\Lambda(s)$ are required,
we can reduce through precomputations the bulk of the work to that of 
computing a main sum.
This comes from the identity
$$
     G(z,w) = w^{-z}\Gamma(z) - g(z,w).
$$
The above discussion indicates that, in (\ref{eq:formula a=1}), we should use $g(z,w)$
and this identity to compute $G(z,w)$ roughly when $|w|$ is smaller than $|z|$. For example,
with $\zeta(1/2+it)$, the region $|w|<|z|$ corresponds in (\ref{eq:zeta}) to
$
    |\pi n^2 \delta^2| < |1/4+it/2|
$
and 
$
    |\pi n^2/ \delta^2| < |1/4-it/2|.
$
Because $|\delta|=1$ this leads to a main sum consisting of approximately 
$|t/(2\pi)|^{1/2}$ terms, the same as in the Riemann-Siegel formula.

\subsection{The approximate functional equation when there is more than one $\Gamma$-factor, and $\kappa_j=1/2$}

In this case, the function $f_1(s,n)$ that appears in Theorem~\ref{thm:formula} is
\begin{equation}
    \label{eq:f case2}
    f_1(s,n) = \frac{\delta^{-s}}{2\pi i}
                \int_{\nu - i \infty}^{\nu + i \infty}
                \prod_{j=1}^a \Gamma((z+s)/2 + \lambda_j)
                z^{-1}
                \br{Q/(\delta n)}^z
                dz.
\end{equation}
This is a special case of the Meijer $G$ function and we develop
some of its properties.

Let $\mellin{\phi(t)}{z}$ denote the Mellin transform of $\phi$
$$
    \mellin{\phi(t)}{z} = \int_0^\infty \phi(t) t^{z-1}.
$$

We will express 
$\prod_{j=1}^a \Gamma((z+s)/2 + \lambda_j) z^{-1}$
as a Mellin transform analogous to~(\ref{eq:invmellin g(v+u)/u}).

Letting $\phi_1 \ast \phi_2$ denote the convolution of two
functions 
$$
    (\phi_1 \ast \phi_2)(v) = \int_0^\infty \phi_1(v/t) \phi_2(t) \frac{dt}{t}
$$
we have (under certain conditions on $\phi_1,\phi_2$)
$$
    \mellin{\phi_1 \ast \phi_2}{z}
    = \mellin{\phi_1}{z} \cdot \mellin{\phi_2}{z}.
$$
Thus
\begin{equation}
    \label{eq:mellin conv}
    \prod_{j=1}^a \mellin{\phi_j}{z}
    = \int_0^\infty (\phi_1 \ast \dots \ast \phi_a)(t) t^{z-1} dt,
\end{equation}
with
$$
    (\phi_1 \ast \dots \ast \phi_a)(v) = 
    \int_0^\infty \dots \int_0^\infty 
    \phi_1(v/t_1) \phi_2(t_1/t_2) \dots \phi_{a-1}(t_{a-2}/t_{a-1}) \phi_a(t_{a-1}) 
    \frac{dt_1}{t_1} \dots \frac{dt_{a-1}}{t_{a-1}}.
$$
Now
$$
    \prod_{j=1}^a \Gamma((z+s)/2 + \lambda_j)
    z^{-1}
    = 
    \br{\prod_{j=1}^{a-1} \Gamma((z+s)/2 + \lambda_j)}
    \br{\Gamma((z+s)/2 + \lambda_a) z^{-1}}.
$$
But
$$
    \Gamma((z+s)/2 + \lambda) = \mellin{2 e^{-t^2} t^{2\lambda +s}}{z},
$$
and~(\ref{eq:invmellin g(v+u)/u}) gives
$$
    \Gamma((z+s)/2 + \lambda) z^{-1} = 
    \mellin{\Gamma(s/2 +\lambda, t^2)}{z}.
$$
So letting
$$
   \phi_j(t) = \begin{cases}
                  2 e^{-t^2} t^{2\lambda_j +s}     & \text{$j=1,\ldots a-1$;} \\
                  \Gamma(s/2 +\lambda_a, t^2)        & \text{$j=a$,}
               \end{cases}
$$
and applying Mellin inversion, we find that~(\ref{eq:f case2}) equals
\begin{equation}
    \label{eq:f case2b}
    f_1(s,n) = \delta^{-s} (\phi_1 \ast \dots \ast \phi_a)(n\delta/Q), 
\end{equation}
where
\begin{align}
    (\phi_1 \ast \dots \ast \phi_a)(v) = 
    v^{2\lambda_1+s}
    \int_0^\infty \dots \int_0^\infty
    &2^{a-1} 
    \prod_{j=1}^{a-1} t_j^{2(\lambda_{j+1}-\lambda_j)} 
    e^{-\br{
              \frac{v^2}{t_1^2} + \frac{t_1^2}{t_2^2} + \dots + \frac{t_{a-2}^2}{t_{a-1}^2}
        }
    } \notag \\
    &\br{
           \int_1^\infty e^{-t_{a-1}^2 x} x^{s/2 +\lambda_a - 1} dx
    }
    \frac{dt_1}{t_1} \dots \frac{dt_{a-1}}{t_{a-1}}. 
    \notag
\end{align}
Substituting $u_j = \frac{(v^2 x)^{j/a}}{v^2} t_j^2$ and rearranging order of integration
this becomes
$$
    v^{2 \mu +s}
    \int_1^\infty
    \Elambda{xv^2}
    x^{s/2+ \mu -1}
    dx,
$$
where
\begin{equation}
    \label{eq:mu}
    \mu = \frac{1}{a}\sum_{l=1}^a \lambda_j,
\end{equation}
\begin{equation}
    \label{eq:Elambda}
    \Elambda{w} =
    \int_0^\infty \dots \int_0^\infty
    \prod_{j=1}^{a-1} u_j^{\lambda_{j+1}-\lambda_j}
    e^{ -w^{1/a}
       \br{
           \frac{1}{u_1} + \frac{u_1}{u_2} + \dots +
           \frac{u_{a-2}}{u_{a-1}} + u_{a-1}
       }
    }
    \frac{du_1}{u_1} \dots \frac{du_{a-1}}{u_{a-1}}.
\end{equation}
So, returning to~(\ref{eq:f case2b}), we find that
$$
    f_1(s,n) = \br{n\delta/Q}^{2 \mu}
               \br{n/Q}^s
               \int_1^\infty \Elambda{x\br{n \delta/Q}^2}
               x^{s/2+ \mu -1}
               dx.
$$
Note that because~(\ref{eq:f case2}) is symmetric in the $\lambda_j$'s, so 
is $E_{\bl}$.

Similarly
$$
    f_2(1-s,n) = \delta^{-1}
               \br{n/(\delta Q)}^{2\cj{\mu}}
               \br{n/Q}^{1-s}
               \int_1^\infty \Ecjlambda{x\br{n/(\delta Q)}^2}
               x^{(1-s)/2+ \cj{\mu} -1}
               dx.
$$
Hence,
%\enlargethispage*{18ex}
\begin{eqnarray}
%\boxed{
%\begin{align}
    Q^s 
    \prod_{j=1}^a \Gamma(s/2 + \lambda_j)
    L(s)
    \delta^{-s} &=&
    \sum_{k=1}^{\ell} \frac{r_k \delta^{-s_k}}{s-s_k} \notag \\
    &+&\br{\delta/Q}^{2\mu}
     \sum_{n=1}^{\infty}
     b(n) n^{2\mu}
     \Glambda{s/2+\mu}{\br{n \delta/Q}^2}
     \notag \\
    &+&\frac{\omega}{\delta} (\delta Q)^{-2\cj{\mu}}
     \sum_{n=1}^{\infty}
     \cj{b}(n)n^{2\cj{\mu}}
     \Gcjlambda{(1-s)/2+\cj{\mu}}{\br{n/(\delta Q)}^2} \notag
%\end{align}
%}
\end{eqnarray}
{\vspace{-4ex} \begin{equation} \label{eq:formula a>1} \end{equation}}
with
$$
    \Glambda{z}{w} = \int_1^\infty \Elambda{xw} x^{z-1} dx
$$
($\mu$ and $E_{\bl}$ are given by~(\ref{eq:mu}),~(\ref{eq:Elambda})).

%-------------------------------------------------------
%\begin{flushleft}
%\hrulefill
%\end{flushleft}

\subsubsection{Examples}

When $a=2$ 
\begin{align}
    &\Elambda{xw} = 
    \int_0^\infty 
    t^{\lambda_2-\lambda_1} 
    e^{-(wx)^{1/2}(1/t+t)}\frac{dt}{t} \notag \\
    \label{eq:K}
    &=2K_{\lambda_2-\lambda_1}\br{2(wx)^{1/2}}
    =2K_{\lambda_1-\lambda_2}\br{2(wx)^{1/2}},
\end{align}
$K$ being the $K$-Bessel function, so that $G_{\bl}$ is an incomplete
integral of the $K$-Bessel function.

Note further that if $\lambda_1=\lambda/2$, $\lambda_2=(\lambda+1)/2$ 
then~(\ref{eq:K}) is
$$
    2K_{1/2}\br{2(wx)^{1/2}} = \br{\pi^{1/2}/(wx)^{1/4}} e^{-2(wx)^{1/2}}
$$
(see~\cite{MR84h:33001b}), so
$G_{(\lambda/2,(\lambda+1)/2)}(z,w) = 2 (2\pi)^{1/2} (4w)^{-z} \Gamma(2z-1/2,2w^{1/2})$,
i.e. the incomplete gamma function.
This is what we expect since, using~(\ref{eq:legendre duplication}), we can write
the gamma factor $\Gamma((s+\lambda)/2) \Gamma((s+\lambda+1)/2)$ in
terms of $\Gamma(s+\lambda)$, for which the $a=1$ expansion,~(\ref{eq:formula a=1}), applies.

\begin{flushleft}

Maass cusp form $L$-functions:
(background material can be found in~\cite{MR97k:11080}).
Let $f$ be a Maass cusp form with eigenvalue $\lambda=1/4-v^2$, 
i.e. $\Delta f=\lambda f$, where 
$\Delta= -y^2(\partial/\partial x^2+\partial/\partial y^2)$,
and Fourier expansion 
\end{flushleft}
$$
    f(z) = \sum_{n\neq 0} a_n y^{1/2} K_v(2\pi \abs{n} y) e^{2\pi i n x},
$$
with $a_{-n}=a_n$ for all $n$, or $a_{-n}=-a_n$ for all $n$.
Let
$$
    L_f(s) = \sum_{n=1}^{\infty} \frac{a_n}{n^s}, \quad \Re{s} > 1
$$
(absolute convergence in this half plane can be proven via the
Rankin-Selberg method), and let
$\e=0$ or $1$ according to whether $a_{-n}=a_n$ or $a_{-n}=-a_n$.
We have that
$$
    \Lambda_f(s) := \pi^{-s} \Gamma((s+\e+v)/2) \Gamma((s+\e-v)/2) L_f(s)
$$
extends to an entire function and satisfies
$$
    \Lambda_f(s) = (-1)^\e \Lambda_f(1-s).
$$
Hence, formula~(\ref{eq:formula a>1}), for $L_f(s)$, is
\begin{align}
    \notag
    &\pi^{-s}
    \Gamma((s+\e+v)/2) \Gamma((s+\e-v)/2) L_f(s)
    \delta^{-s} = \\ \notag
    &\br{\delta \pi}^{\e}
     \sum_{n=1}^{\infty}
     a_n n^{\e}
     G_{((\e+v)/2,(\e-v)/2)}\br{s/2+\e/2, \br{n \delta \pi}^2}
     \notag \\
    &+\frac{(-1)^\e}{\delta} (\pi/\delta)^{\e}
     \sum_{n=1}^{\infty}
     a_n n^{\e}
     G_{((\e+\cj{v})/2,(\e-\cj{v})/2)}\br{(1-s)/2+\e/2,\br{n\pi/\delta}^2}
     \notag
\end{align}
where, by~(\ref{eq:K}),
\begin{align}
    G_{((\e+v)/2,(\e-v)/2)}\br{s/2+\e/2, \br{n \delta \pi}^2} 
     &=  4 \int_1^\infty K_v(2n\delta \pi t) t^{s+\e-1} dt
     \notag \\
    G_{((\e+\cj{v})/2,(\e-\cj{v})/2)}\br{(1-s)/2+\e/2, \br{n \pi/\delta}^2} 
     &=  4 \int_1^\infty K_{\cj{v}}(2n \pi t/\delta) t^{-s+\e} dt.
     \notag
\end{align}

%\begin{flushleft}
%\hrulefill
%\end{flushleft}

%-------------------------------------------------------

Next, let
\begin{align}
    \Gammalambda{z,w} 
    &= w^z \Glambda{z}{w} = \int_w^\infty \Elambda{t}t^{z-1} dt, \notag \\ 
    \label{eq:Gammalambda b}
    \Gammalambda{z}
    &=\int_0^\infty \Elambda{t}t^{z-1} dt,  \\
    \gammalambda{z,w} 
    &= \int_0^w \Elambda{t}t^{z-1} dt, \notag
\end{align}
with $E_{\bl}$ given by~(\ref{eq:Elambda}).

\begin{lemma}
$$
    \Gammalambda{z} = \prod_{j=1}^a \Gamma(z-\mu+\lambda_j)
$$
where $\mu =\frac{1}{a} \sum_{j=1}^{a} \lambda_j$.
\end{lemma}

\begin{proof}
Let $\psi_j(t) = e^{-t} t^{\lambda_j}$, $j=1,\ldots,a$, and consider
\begin{align}
    (\psi_1 \ast \dots \ast \psi_a)(v) 
    &= v^{\lambda_1}
    \int_0^\infty \dots \int_0^\infty
    \prod_{j=1}^{a-1} t_j^{\lambda_{j+1}-\lambda_j}
    e^{ -
       \br{
           \frac{v}{t_1} + \frac{t_1}{t_2} + \dots +
           \frac{t_{a-2}}{t_{a-1}} + t_{a-1}
       }
    }
    \frac{dt_1}{t_1} \dots \frac{dt_{a-1}}{t_{a-1}} \notag \\
    &= v^\mu 
    \int_0^\infty \dots \int_0^\infty
    \prod_{j=1}^{a-1} x_j^{\lambda_{j+1}-\lambda_j}
    e^{ -v^{1/a}
       \br{
           \frac{1}{x_1} + \frac{x_1}{x_2} + \dots +
           \frac{x_{a-2}}{x_{a-1}} + x_{a-1}
       }
    }
    \frac{dx_1}{x_1} \dots \frac{dx_{a-1}}{x_{a-1}} \notag.
\end{align}
(we have put $t_j=v^{1-j/a} x_j$).
Thus, from~(\ref{eq:Elambda})
$$
    \Elambda{v} =v^{-\mu} (\psi_1 \ast \dots \ast \psi_a)(v),
$$
and hence~(\ref{eq:Gammalambda b}) equals
$$
    \int_0^\infty (\psi_1 \ast \dots \ast \psi_a)(t) t^{z-\mu-1} dt
$$
which, by~(\ref{eq:mellin conv}) is
$\prod_{j=1}^a \Gamma(z-\mu+\lambda_j)$.

\end{proof}

Inverting, we get
\begin{equation}
    \notag
    \Elambda{t} 
    =\frac{1}{2\pi i} \int_{\nu-i\infty}^{\nu+i\infty}
    \Gammalambda{z}
    t^{-z} 
    dz
\end{equation}
with $\nu$ to the right of the poles of $\Gammalambda{z}$. Shifting
the line integral to the left, we can express $\Elambda{t}$ as a sum of residues,
and hence obtain through termwise integration a series expansion for 
$\gammalambda{z,w}$. An algorithm for doing so is detailed in~\cite{Do},
though with different notation.
Such an expansion is useful for $|w| << 1$.
That paper also describes how to obtain an asymptotic expansion for
$\Elambda{t}$ and hence, by termwise integration, for $\Gammalambda{z,w}$,
useful for $|w|$ large in comparison to $|z|$. 
The paper has, implictly, $g(z)=1$ and 
does not control for cancellation.  Consequently, it does not provide a
means to compute $L$-functions away from the real axis other than 
increasing precision. 

If one wishes to use the methods of this paper
to control for cancellation, then one will have $w$ varying over a wide range of
values for which the series expansion in~\cite{Do} is not adequate. We thus need
an alternative method to compute $\Glambda{z}{w}$ especially in the transition
zone $|z| \approx |w|$. It would be useful to have Temme's uniform asymptotics
generalized to handle $\Glambda{z}{w}$. Alternatively, we can apply the
naive but powerful Riemann sum technique described in section 2.

\subsection{The functions $f_1(s,n), f_2(1-s,n)$ as Riemann sums}

Substituting $z=v+iu$ into (\ref{eq:mellin}) we have
$$
    f_1(s,n) =  
    \frac{1}{2\pi}
    \int_{-\infty}^{\infty}
    \prod_{j=1}^a \Gamma(\kappa_j (s+v+iu) + \lambda_j)
    \frac{g(s+v+iu)}{v+iu}
    (Q/n)^{v+iu}
    du.
$$
Let
$$ 
    h(u) =
    \frac{1}{2\pi}
    \prod_{j=1}^a \Gamma(\kappa_j (s+v+iu) + \lambda_j)
    \frac{g(s+v+iu)}{v+iu}
    (Q/n)^{v+iu}.
$$
With the choice of $g(z)$ as in (\ref{eq:g(s) 2}), an analysis similar to that 
following~(\ref{eq:hat ft(y)}) shows that $\hat{h}(y)$ decays exponentially fast 
as $y \to -\infty$, and doubly exponentially fast as $y \to \infty$. Hence, we
can successfully evaluate $f_1(s,n)$, and similarly $f_2(1-s,n)$ as simple 
Riemann sums, with step size inversely proportional to the number of digits of 
precision required.

The Riemann sum approach gives us tremendous flexibility. We are no longer
bound in our choice of $g(z)$ to functions for which~(\ref{eq:mellin})
has nice series or asymptotic expansions. For example, 
we can, with $A>0$, set
\begin{equation}
    \notag
    g(z) = \exp(A(z-s)^2) \prod_{j=1}^a \delta_j^{-\kappa_j z}.
\end{equation}
The extra factor $\exp(A(z-s)^2)$ is chosen so as to cut down on the domain
of integration. Recall that in $f_1(s,n)$ and $f_2(1-s,n)$, $g$ appears as 
$g(s \pm (v+iu))$, hence $\exp(A(z-s)^2)$ decays in the integral like
$\exp(-Au^2)$. Ideally, we would like to have $A$ large. However, this
would cause the Fourier transform $\hat{h}(y)$ to decay too slowly. The Fourier
transform of a product is a convolution of Fourier transforms, and
the Fourier transform of $\exp(A(v+iu)^2)$ equals
$$
    (\pi/A)^{1/2} \exp(\pi y(2 A v - \pi y)/A).
$$ 
A large value of $A$ leads to
to a small $1/A$ and this results in poor performance of $\hat{h}(y)$.
We also need to specify $v$, for the line of integration. 
Larger $v$ means more rapid decay of $\hat{h}(y)$ but more
cancellation in the Riemann sum and hence loss of precision.

Another advantage to the Riemann sum approach is that we can rearrange
sums, putting the Riemann sum on the outside and the sum over $n$ on the inside.
Both sums are finite since we truncate them once the tails are within
the desired precision.
This then expresses, to within an error that we can control by our choice of
stepsize and truncation,
$\Lambda(s)$ as a sum of finite Dirichlet series evaluated at equally
spaced points and hence gives a sort of interpolation formula for $\Lambda(s)$.
Details related to this approach will appear in a future paper.

\subsection{Looking for zeros}

To look for zeros of an $L$-function, we can rotate it so that it is
real on the critical line, for example working with $Z(t)$, see~(\ref{eq:Z(t)}), 
rather than $\zeta(1/2+it)$.

We can then advance in small steps, say one quarter the average gap size between consecutive zeros,
looking for sign changes of this real valued function, zooming in each time a sign change
occurs. Along the way, we need to 
determine if any zeros have been missed, and, if so, go back and look for
them, using more refined step sizes. We can also use more sophisticated 
interpolation techniques to make the search for zeros more efficient
\cite{O}. If this search fails to turn up the missing zeros,
then presumably a bug has crept into one's code, or else
one should look for zeros of the $L$-function nearby but off the critical line in violation
of the Riemann hypothesis. 

To check for missing zeros, we could use the argument principle and numerically integrate
the logarithmic derivative of the $L$-function along a rectangle, rounding to the closest
integer. However, this is inefficient and difficult to make numerically rigorous.

It is better to use a test devised by Alan Turing~\cite{Tu} for $\zeta(s)$ but which
seems to work well in general. 
Let $N(T)$ denote the number of zeros of $\zeta(s)$
in the critical strip above the real axis and up to height $T$:
\begin{equation}
    \notag
    N(T) = 
    \left|
    \{ 
        \rho = \beta +i \gamma | \zeta(\rho)=0, 0 \leq \beta \leq 1, 0 < \gamma \leq T
    \}
    \right|.
\end{equation}
A theorem of von Mangoldt states that 
\begin{equation}
    \label{eq:N(T) asympt}
    N(T) = \frac{T}{2\pi} \log(T/(2\pi)) -\frac{T}{2\pi} + \frac{7}{8} +S(T) +O(T^{-1})
\end{equation}
with
$$
    S(T) = O(\log{T}).
$$
However, a stronger inequality due to Littlewood and with explicit constants due
to Turing~\cite{Tu}~\cite{Le} is given by
\begin{equation}
    \label{eq:littlewood-turing}
    \left|
        \int_{t_1}^{t_2} S(t) dt
    \right|
    \leq 2.3 + .128 \log(t_2/\pi)
\end{equation}
for all $t_2 > t_1 > 168 \pi$, i.e. $S(T)$ is 0 on average. Therefore,
if we miss one sign change (at least two zeros), we'll quickly detect
the fact.
To illustrate this, Table~\ref{tab:missing zeros} contains a list of the imaginary parts of the
zeros of $\zeta(s)$ found naively by searching for sign changes of $Z(t)$ taking step 
sizes equal to two.  We notice that near the ninth zero on our list a missing pair is
detected, and similary near the twenty fifth zero. A more refined
search reveals the pairs of zeros with imaginary parts equal to 
$48.0051508812$, $49.7738324777$, and $94.6513440405$, 
$95.8706342282$ respectively.

It would be useful to have a general form of the explicit inequality~(\ref{eq:littlewood-turing})
worked out for any $L$-function. The papers of Rumely~\cite{Rum}
and Tollis~\cite{To} generalize this inequality to Dirichlet $L$-functions
and Dedekind zeta functions respectively. 

The main term, analogous 
to~(\ref{eq:N(T) asympt}), for a general $L$-function is easy to derive.
Let $L(s)$ be an $L$-function with functional equation as described 
in~(\ref{eq:lambda}).
Let $N_L(T)$ denote the number of zeros of $L(s)$
lying within the rectangle $|\Im{s}| \leq T$, $0 < \Re{s} < 1$.
Notice here we are considering zeros lying both above and below the real axis
since the zeros of $L(s)$ will not be located symmetrically about the real axis
if its Dirichlet coefficients $b(n)$ are non-real.

Assume for simplicity that $L(s)$ is entire. 
The arguement principle and the functional equation for $L(s)$
suggests a main term for $N_L(T)$ equal to 
\begin{equation}
     \notag
     N_L(T) \sim 
     \frac{2 T}{\pi} \log(Q) +   
     \frac{1}{\pi} \sum_{j=1}^a 
     \Im \left(
         \log\left(
              \frac{\Gamma((1/2+iT)\kappa_j+\lambda_j)}
                   {\Gamma((1/2-iT)\kappa_j+\lambda_j)}
             \right)
     \right).
\end{equation}
If we assume further that the $\lambda_j$'s are all real, then the above
is, by Stirling's formula, asymptotically equal to
\begin{equation}
     \notag
     N_L(T) \sim 
     \frac{2 T}{\pi} \log(Q) +   
     \sum_{j=1}^a 
     \left(
         \frac{2T\kappa_j}{\pi} \log(T\kappa_j/e) 
         +(\kappa_j/2+\lambda_j-1/2)
     \right).
\end{equation}
A slight modification of the above
is needed if $L(s)$ has poles, as in the case of $\zeta(s)$.
See Davenport~\cite[chapters 15,16]{D} where
rigorous proofs are presented for $\zeta(s)$ and Dirichlet $L$-functions (the original
proof is due to von Mangoldt).

\begin{table}[h!]
\centerline{
\begin{tabular}{|c|c|c|}
\hline
$j$ & $t_j$ & $\tilde{N}((t_j+t_{j-1})/2)-j+1$ \cr \hline
1 & 14.1347251417 & -0.11752 \cr 
2 & 21.0220396388 & -0.04445 \cr 
3 & 25.0108575801 & -0.03216 \cr 
4 & 30.4248761259 & 0.01102 \cr 
5 & 32.9350615877 & -0.01000 \cr 
6 & 37.5861781588 & -0.05699 \cr 
7 & 40.9187190121 & 0.07354 \cr 
8 & 43.3270732809 & -0.07314 \cr 
9 & 52.9703214777 & 0.81717 \cr 
10 & 56.4462476971 & 2.01126 \cr 
11 & 59.3470440026 & 2.12394 \cr 
12 & 60.8317785246 & 1.90550 \cr 
13 & 65.1125440481 & 1.95229 \cr 
14 & 67.0798105295 & 2.11039 \cr 
15 & 69.5464017112 & 1.94654 \cr 
16 & 72.0671576745 & 1.90075 \cr 
17 & 75.7046906991 & 2.09822 \cr 
18 & 77.1448400689 & 2.10097 \cr 
19 & 79.3373750202 & 1.82662 \cr 
20 & 82.9103808541 & 1.99205 \cr 
21 & 84.7354929805 & 2.09800 \cr 
22 & 87.4252746131 & 2.03363 \cr 
23 & 88.8091112076 & 1.88592 \cr 
24 & 92.4918992706 & 1.95640 \cr 
25 & 98.8311942182 & 3.10677 \cr 
26 & 101.3178510057 & 4.03517 \cr 
27 & 103.7255380405 & 4.11799 \cr 
\hline
\end{tabular}
}
\caption{Checking for missing zeros. The second column lists the imaginary parts of
the zeros of $\zeta(s)$ found by looking for sign changes of $Z(t)$, 
advancing in step sizes equal to two. The third column compares the number of zeros 
found to the main term of $N(T)$, namely to
$\tilde{N}(T):=(T/(2\pi)) \log(T/(2\pi e)) + 7/8$,
evaluated at the midpoint between consecutive zeros, with $t_0$ taken to be 0.
This detects a pair of missing zeros near the ninth and twenty fifth
zeros on our list.
}\label{tab:missing zeros}
\end{table}

\section{Experiments involving $L$-functions}

Here we describe some of the experiments that reflect the random matrix theory philosophy, 
namely that the zeros and values of $L$-functions behave like the zeros and values of
characteristic functions from the classical compact groups~\cite{KS2}. Consequently, 
we are interested in questions concerning the distribution of zeros, horizontal and 
vertical, and the value distribution of $L$-functions.

\subsection{Horizontal distribution of the zeros}

Riemann himself computed the first few zeros of $\zeta(s)$, and detailed
numerical studies were initiated almost as soon as computers were invented. See
Edwards~\cite{E} for a historical survey of these computations.
To date, the most impressive computations for $\zeta(s)$ have been those of Odlyzko
~\cite{O}~\cite{O2} and Wedeniwski~\cite{W}. The latter 
adapted code of van de Lune, te Riele, and Winter~\cite{LRW} for
grid computing over the internet. Several thousand computers have been used to verify
that the first $8.5 \cdot 10^{11}$ nontrivial zeros
of $\zeta(s)$ fall on the critical line. Odlyzko's computations
have been more concerned with examining the distribution of the spacings
between neighbouring zeros, although the Riemann Hypothesis has also been
checked for the intervals examined. In~\cite{O}, Odlyzko computed 175 million 
consecutive zeros of $\zeta(s)$ lying near the $10^{20}$th zero, and more recently, 
billions of zeros in a higher region~\cite{O2}. The Riemann-Siegel
formula has been at the heart of these computations. Odlyzko also uses
FFT and interpolation algorithms to allow for many evaluations of $\zeta(s)$
at almost the same cost of a single evaluation.

Dirichlet $L$-functions were not computed on machines until 1961 when
Davies and Haselgrove~\cite{DH} looked at several
$L(s,\chi)$ with conductor $\leq 163$.
Rumely~\cite{Rum},
using summation by parts, computed the first several thousand
zeros for many Dirichlet $L$-functions with small moduli.
He both verified RH and looked at statistics of neighbouring zeros.

Yoshida~\cite{Y}~\cite{Y2} has also used summation by parts,
though in a different manner, to compute the first few zeros of certain higher degree,
with two or more $\Gamma$-factors in the functional equation,
$L$-functions.

Lagarias and Odlyzko~\cite{MR80g:12010} have computed the low lying zeros of
several Artin $L$-functions
using expansions involving the incomplete gamma function.
They noted that one could compute higher up in the critical strip
by introducing the parameter $\delta$, as explained in section~\ref{sec:g(z)},
but did not implement it since it led to difficulties
concerning the computation of $G(z,w)$ with both $z$ and $w$ complex.

Other computations of $L$-functions include those of
Berry and Keating~\cite{BK} and Paris~\cite{P} ($\zeta(s)$),
Tollis~\cite{To} (Dedekind zeta functions), 
Keiper~\cite{Ke} and Spira~\cite{Sp} (Ramanujan $\tau$ $L$-function), 
Fermigier~\cite{F} and Akiyama-Tanigawa~\cite{AT} (elliptic curve $L$-functions),
Strombergsson~\cite{St} and Farmer-Kranec-Lemurell~\cite{FKL}
(Maass waveform $L$-functions), and
Dokchister~\cite{Do} (general $L$-functions near the critical line).

The author has verified the Riemann hypothesis for various $L$-functions.
These computations use the methods described in section~\ref{sec:L algorithms}
and are not rigorous in the sense that no attempt is made to obtain explicit bounds
for truncation errors on some of the asymptotic expansions and continued fractions used,
and no interval arithmetic to bound round off errors is carried out. Tables of the zeros
mentioned may be obtained from the author's homepage~ \cite{Ru4}.
These include the first tens of millions zeros of 
all $L(s,\chi)$ with the conductor of $\chi$ less than 20, the first $300000$ zeros
of $L_\tau(s)$, the Ramanujan $\tau$ $L$-function, the first $100000$ zeros of
the $L$-functions associated to elliptic curves of conductors $11,14,15,17,19$, 
the first $1000$ zeros for elliptic curves of conductors less than 1000, the first
100 zeros of elliptic curves with conductor less than 8000, 
and hundreds/millions of zeros of many other $L$-functions.

In all these computations,
no violations of the Riemann hypothesis have been found. 

\subsection{Vertical distribution: correlations and spacings distributions}

The random matrix philosophy predicts that various 
statistics of the zeros of $L$-functions
will mimic the same statistics for the eigenvalues of matrices
in the classical compact groups.

Montgomery~\cite{Mo} achieved the first result connecting zeros of
$\zeta(s)$ with eigenvalues of unitary matrices. Write a typical non-trivial
zero of $\zeta$ as
$$
     1/2 + i\gamma.
$$
Assume the Riemann Hypothesis, so that the $\gamma$'s are real. Because the zeros 
of $\zeta(s)$ come in conjugate pairs, we can restrict our attention
to those lying above the real axis and order them
$$
    0 < \gamma_1 \leq \gamma_2 \leq \gamma_3 \ldots
$$
We can then ask how the spacings between consecutive zeros,
$\gamma_{i+1}-\gamma_i$, are distributed, but first, we need to 'unfold' the zeros
to compensate for the fact that the zeros on average become closer as one goes
higher in the critical strip. We set
\begin{equation}
    \label{eq:normalized spacings}
    \tilde{\gamma}_i = \gamma_i \frac{\log(\gamma_i/(2\pi e))}{2\pi}
\end{equation}
and investigate questions involving the $\tilde{\gamma}$'s. 
This normalization
is chosen so that the mean spacing between consecutive 
$\tilde{\gamma}$'s equals one. Summing the consecutive differecnces,
we get a telespcoping sum
$$
    \sum_{\gamma_i  \leq T} (\tilde{\gamma}_{i+1} -\tilde{\gamma}_i)
    = \tilde{\gamma}(T)+O(1)=
    \gamma(T) \frac{\log(\gamma(T)/(2\pi e))}{2\pi}+O(1)
$$
where $\gamma(T)$ is the largest $\gamma$ less than or equal to $T$.
By~(\ref{eq:N(T) asympt}), the r.h.s above equals
$$
     N(\gamma(T))+O(\log(\gamma(T))) = N(T)+O(\log(T)),
$$
hence $\tilde{\gamma}_{i+1} -\tilde{\gamma}_i$ has mean spacing equal to one.

From a theoretical point
of view, studying the consecutive spacings distribution is difficult
since this assumes the ability to sort the zeros.
The tool that is used for studying spacings questions about the zeros, 
namely the explicit formula, involves a sum over all zeros of
$\zeta(s)$, and it is easier to consider the pair correlation, a statistic
incorporating differences between all pairs of zeros. 
Montgomery conjectured that
for $0 \leq \alpha < \beta$ and $M \to \infty$, 
\begin{eqnarray}
    \label{eq:pair correlation zeta}
    M^{-1} 
    | \{ 1 \leq i<j \leq M : \tilde{\gamma}_j -\tilde{\gamma}_i 
     \in [\alpha,\beta) \} | \notag \\
    \sim 
    \int_\alpha^\beta 
    \left( 1 - \left(\frac{\sin \pi t}{\pi t} \right)^2 \right) dt.
\end{eqnarray}
Notice that $M^{-1}$, and not, say, ${M \choose 2}^{-1}$, is
the correct normalization. For any $j$ there, are just a handful of $i$'s
with $\tilde{\gamma}_j-\tilde{\gamma}_i \in [\alpha,\beta)$.

Montgomery was able to prove that
\begin{equation}
     \label{eq:pair correlation montgomery}
     M^{-1} \sum_{1 \leq i < j \leq M} f(\tilde{\gamma}_j-\tilde{\gamma}_i)
     \to
     \int_0^\infty f(t) \left( 1 - \left(\frac{\sin \pi t}{\pi t} \right)^2 \right) dt.
\end{equation}
as $M \to \infty$,
for test functions $f$ satisfying the stringent restriction that $\hat{f}$
be supported in $(-1,1)$. 

An equivalent way to state the conjecture as $M \to \infty$, 
and one which Odlyzko uses in his numerical experiments,
is to let
\begin{equation}
    \label{eq:normalized spacings2}
    \delta_i = ( \gamma_{i+1} - \gamma_i ) 
    \frac{\log(\gamma_i/(2\pi))}{2\pi}.
\end{equation}
and replace the condition $\tilde{\gamma}_j-\tilde{\gamma}_i \in [\alpha,\beta)$
with the condition
$\delta_i + \delta_{i+1} + \cdots + \delta_{i+k} \in [\alpha,\beta)$ 
for $1 \leq i \leq M , k \geq 0$. The main difference is the
absence of the $1/e$ in the logarithm. 
This is done so as to maintain
a mean spacing tightly asymptotic to one. Set
$$
    C(T) = \sum_{\gamma_i \leq T} ( \gamma_{i+1} - \gamma_i),
$$
and sum by parts
\begin{equation}
     \notag
    \sum_{\gamma_i \leq T} 
    \delta_i
    = C(T) \frac{\log(T/(2\pi))}{2\pi}
    -\frac{1}{2\pi} \int_{\gamma_1}^T C(t) \frac{dt}{t}.
\end{equation}
Now, $C(t)$ telescopes, and von Mangoldt's formula
~(\ref{eq:N(T) asympt}) implies that $C(t)=t+O(1)$, so that
the r.h.s above equals $N(T)+O(\log(T))$, and
$\delta_i$ is on average equal to one. 
In carrying out
numerical experiments with zeros one can either use the normalization
given in~(\ref{eq:normalized spacings}) or~(\ref{eq:normalized spacings2}).
For the theoretical purpose of examining leadings asymptotics of, say,
the pair correlation, the factors appearing in these normalizations
in the logarithm, $1/(2\pi e)$ or $1/(2\pi)$,
are not important as they only affect lower order terms.
However, for the purpose of comparing numerical data to theoretical predictions
it is crucial to include them.

\begin{figure}[htp]
    \centerline{
        \psfig{figure=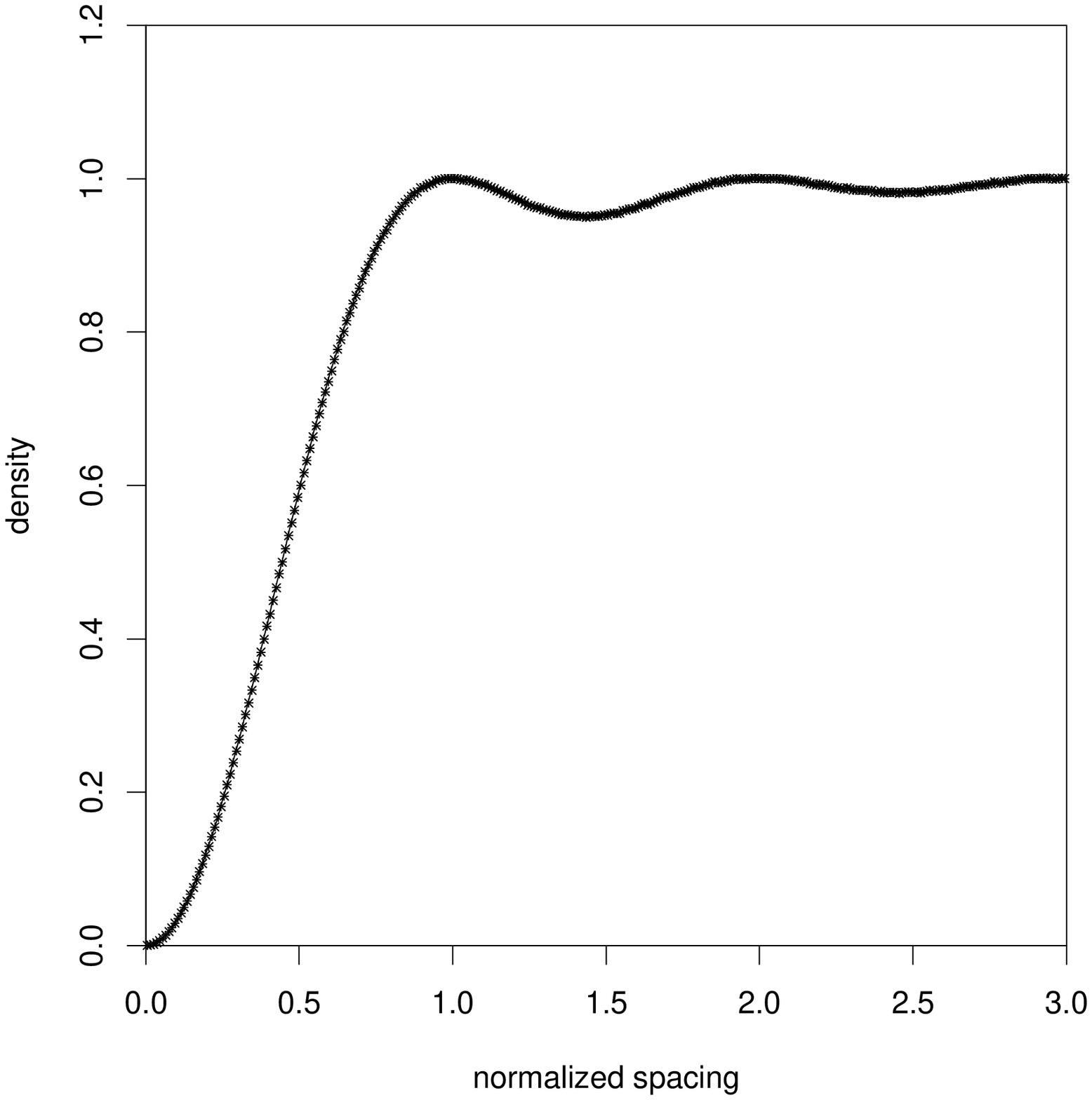,width=3in,angle=-90} 
        \psfig{figure=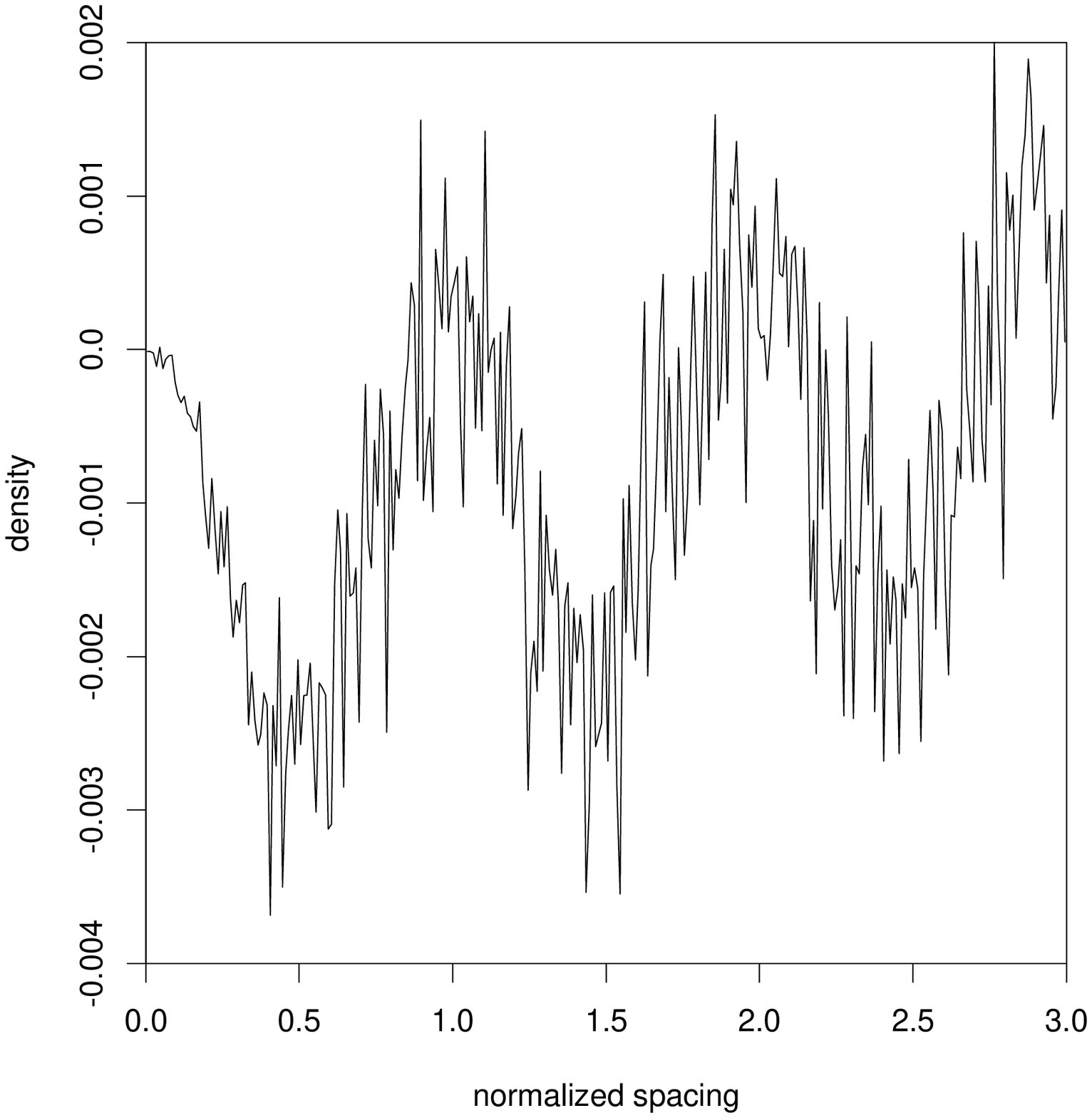,width=3in,angle=-90} 
    }
    \caption
    {The first graph depicts Odlyzko's pair correlation picture for $2\times 10^8$ zeros of $\zeta(s)$ near the 
     $10^{23}$rd zero. The second graph shows the difference between the histogram in the first graph and
     $1 - \left((\sin \pi t)/(\pi t) \right)^2$. In the interval displayed, the two
     agree to within about $.002$.
    }
    \label{fig:odlyzko1}
\end{figure}

On a visit by Montgomery to the the Institute for Advanced Study,
Freeman Dyson out that large unitary matrices have the same 
pair correlation. Let 
$$
    e^{i \theta_1}, e^{i \theta_2}, \ldots, e^{i \theta_N}
$$
be the eigenvalues of a matrix in $\text{U}(N)$, sorted so that
$$
    0 \leq \theta_1 \leq \theta_2 \ldots \leq \theta_N < 2 \pi.
$$
Normalize the eigenangles
\begin{equation}
     \label{eq:normalized eigenangles}
     \tilde{\theta}_i = \theta_i N / (2\pi)
\end{equation}
so that $\tilde{\theta}_{i+1}-\tilde{\theta}_i$ equals one on average.
Then, a classic result in random matrix theory~\cite{M} asserts that
\begin{equation}
     \notag
    N^{-1}
    | \{ 1 \leq i < j \leq N ,
    \tilde{\theta}_j-\tilde{\theta}_i \in [\alpha,\beta) \} |
\end{equation}
equals, when averaged according to Haar measure over $\text{U}(N)$ and letting
$N \to \infty$,
$$
    \int_\alpha^\beta
    \left( 1 - \left(\frac{\sin \pi t}{\pi t} \right)^2 \right) dt.
$$
Odlyzko~\cite{O}~\cite{O2} has carried out numerics to verify Montgomery's
conjecture~(\ref{eq:pair correlation zeta}). His most extensive data to date 
involves billions of zeros near the $10^{23}$rd zero of $\zeta(s)$.
With kind 
permission we reproduce~\cite{O4} Odlyzko's pair correlation picture in figure~\ref{fig:odlyzko1}.

This picture compares the l.h.s. of ~(\ref{eq:pair correlation zeta}) for many 
bins $[a,b)$ of size $b-a=.01$ to the curve
$$
     1 - \left(\frac{\sin \pi t}{\pi t} \right)^2.
$$
Odlyzko's histogram fits the theoretical prediction beautifully. 
Bogomolny and Keating~\cite{K}~\cite{BoK},
using conjectures of Hardy and Littlewood,
have explained the role played by secondary 
terms in the pair correlation of the zeros of $\zeta(s)$ and these terms are 
related to $\zeta(s)$ on the one line. A nice description of these results are
contained in~\cite{BK2}.
Recently, Conrey and Snaith~\cite{CS} obtained the main and lower terms 
of the pair correlation using 
a conjecture for the full asymptotics of the average value of a ratio of 
four zeta functions rather than the Hardy-Littlewood conjectures.

Montgomery's pair correlation theorem~(\ref{eq:pair correlation montgomery}) 
has been generalized by Rudnick and Sarnak~\cite{RudS} 
to any primitive $L$-function, i.e. one which does not factor as a product
of other $L$-functions, as well as to higher correlations which are defined in
a way similar to the pair correlation. Again, there are severe restrictions on the
fourier transform of the allowable test functions, and further, for $L$-functions
of degree greater than three, Rudnick and Sarnak assume 
a weak form of the the Ramanujan conjectures. 
Bogomolny and Keating provide a heuristic derivation of the higher correlations
of the zeros of $\zeta(s)$ using the Hardy-Littlewood conjectures~\cite{BoK2}.

\begin{figure}[htp]
    \centerline{
            \psfig{figure=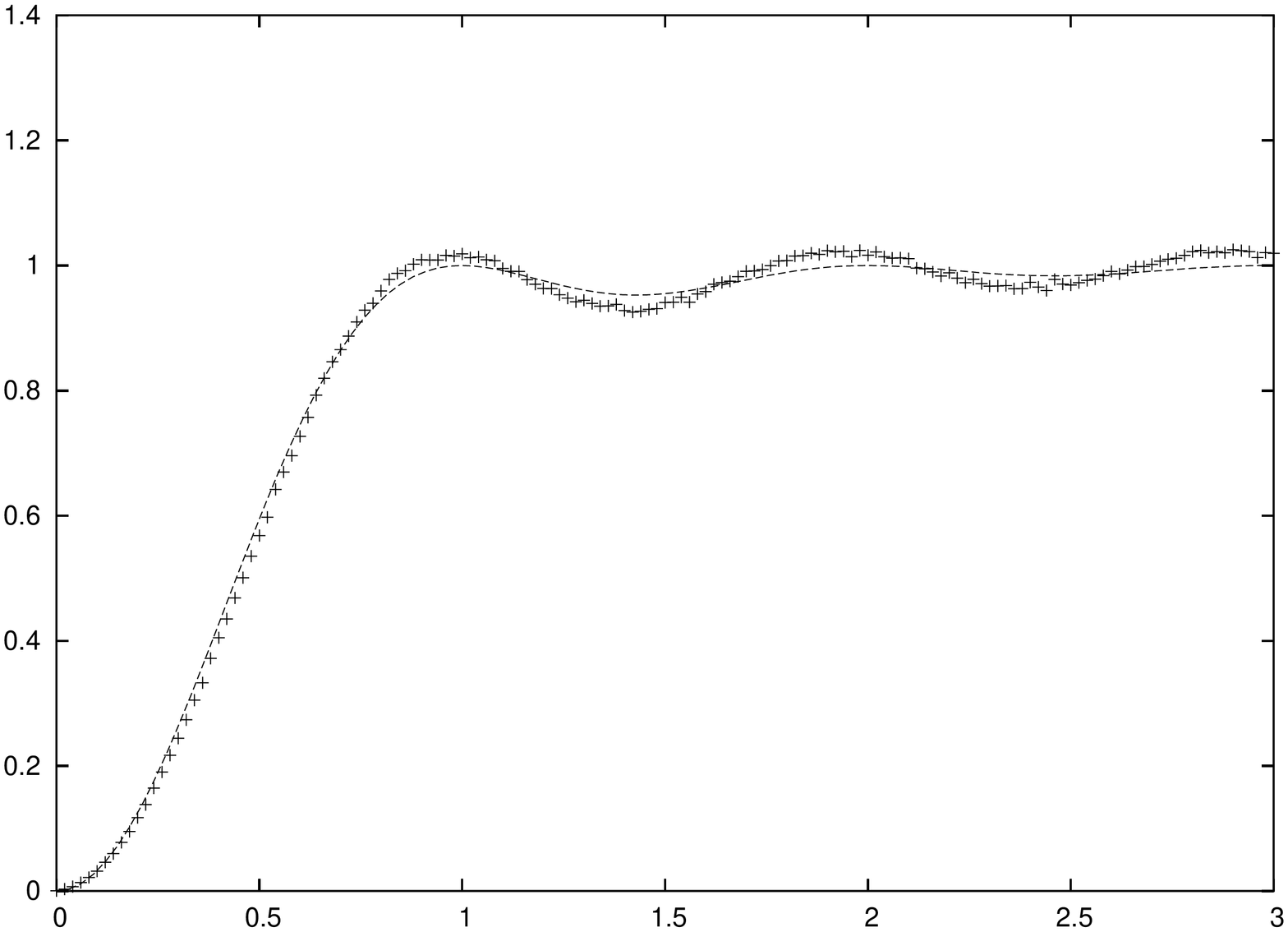,width=2in,angle=-90}
            \psfig{figure=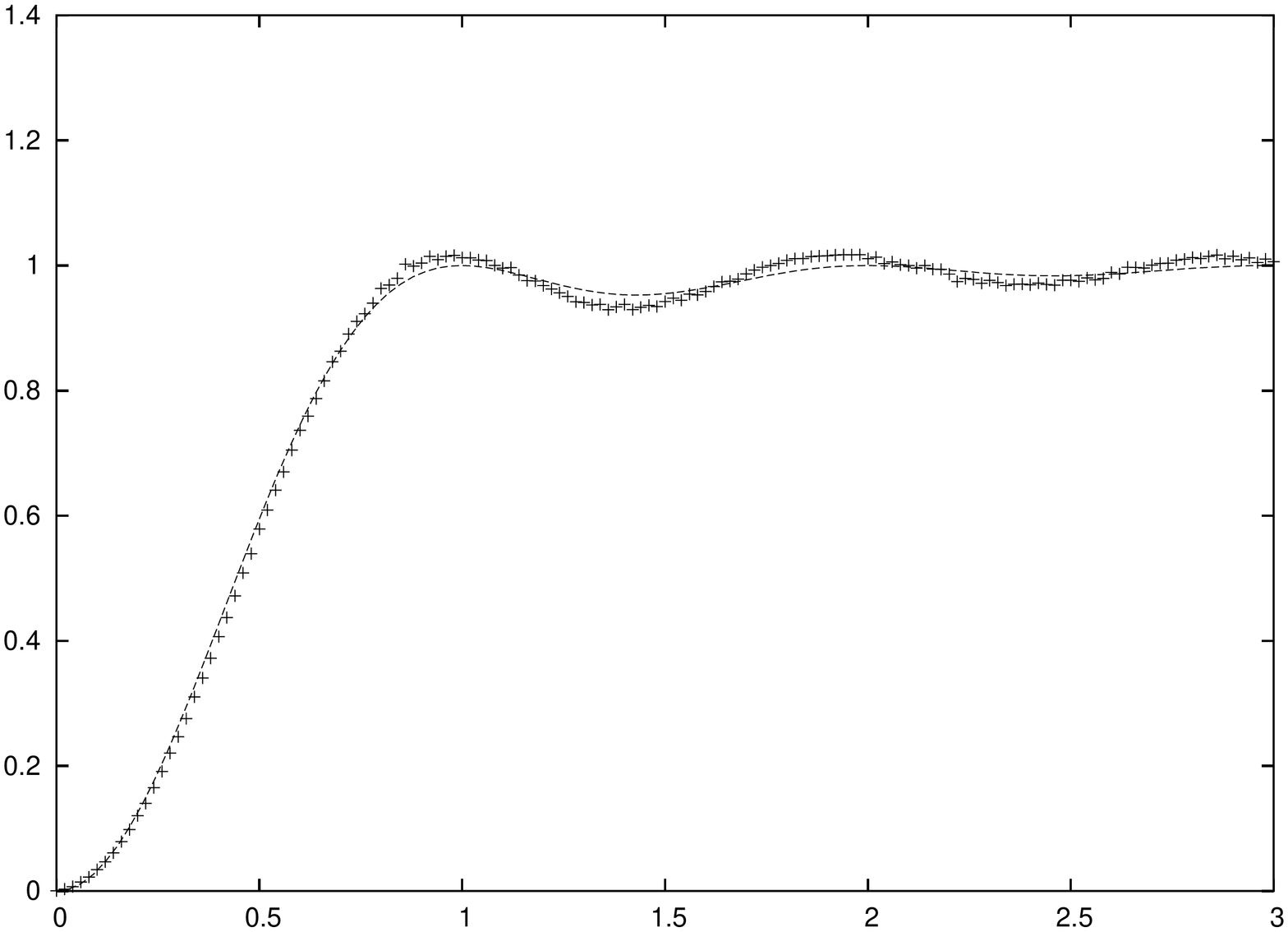,width=2in,angle=-90}
            \psfig{figure=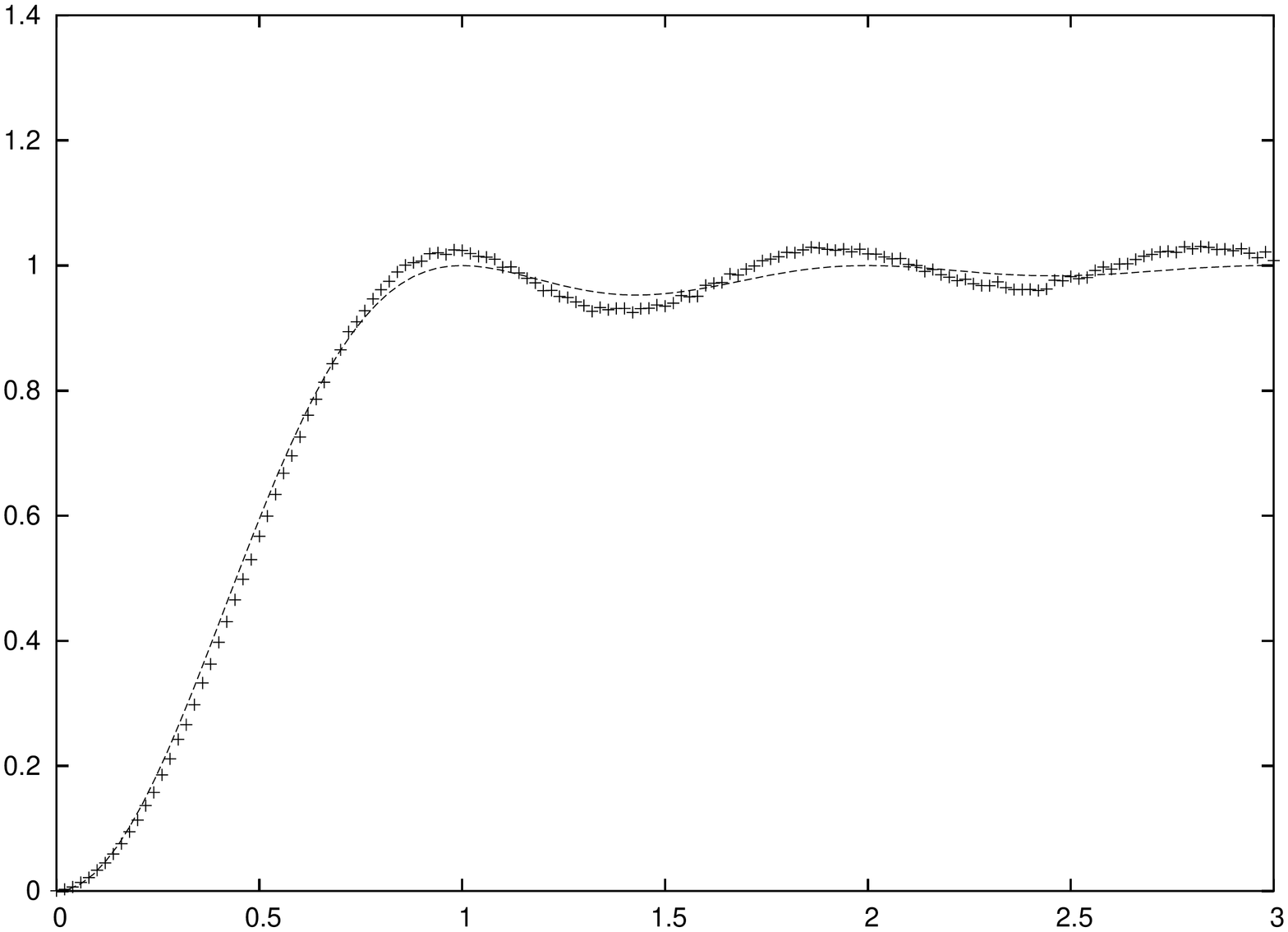,width=2in,angle=-90}
    }
    \centerline{
            \psfig{figure=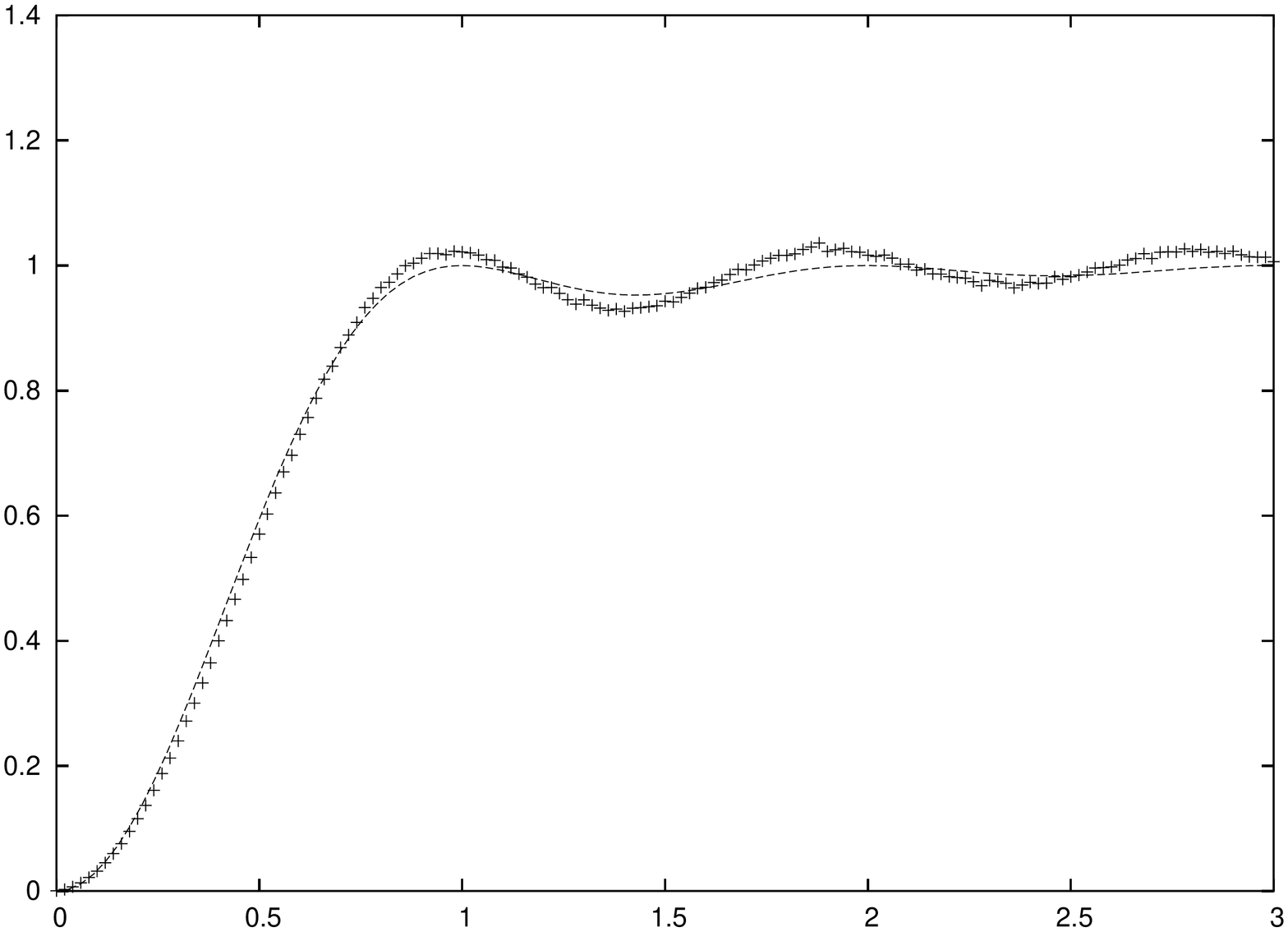,width=2in,angle=-90}
            \psfig{figure=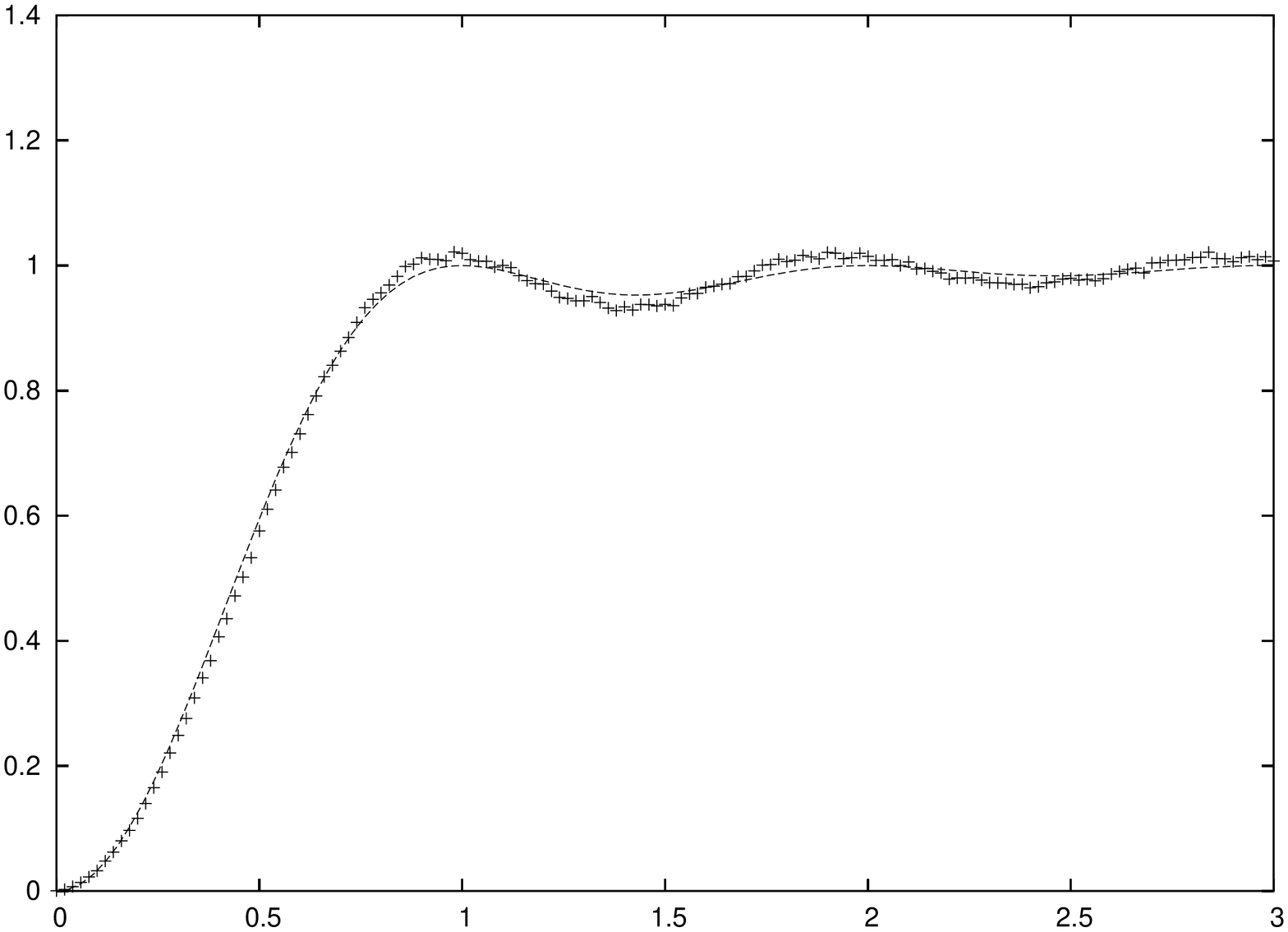,width=2in,angle=-90}
            \psfig{figure=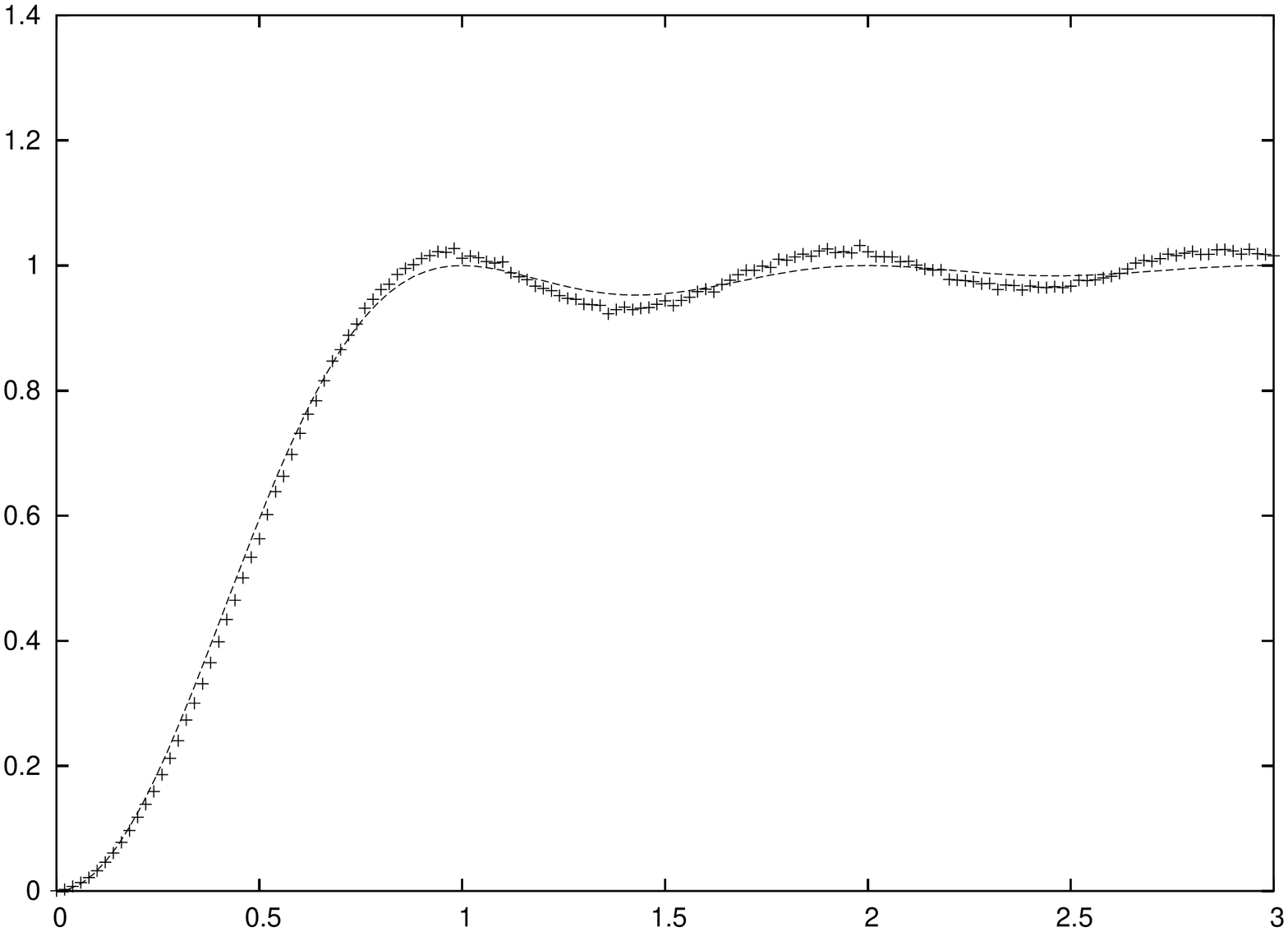,width=2in,angle=-90}
    }
    \centerline{
            \psfig{figure=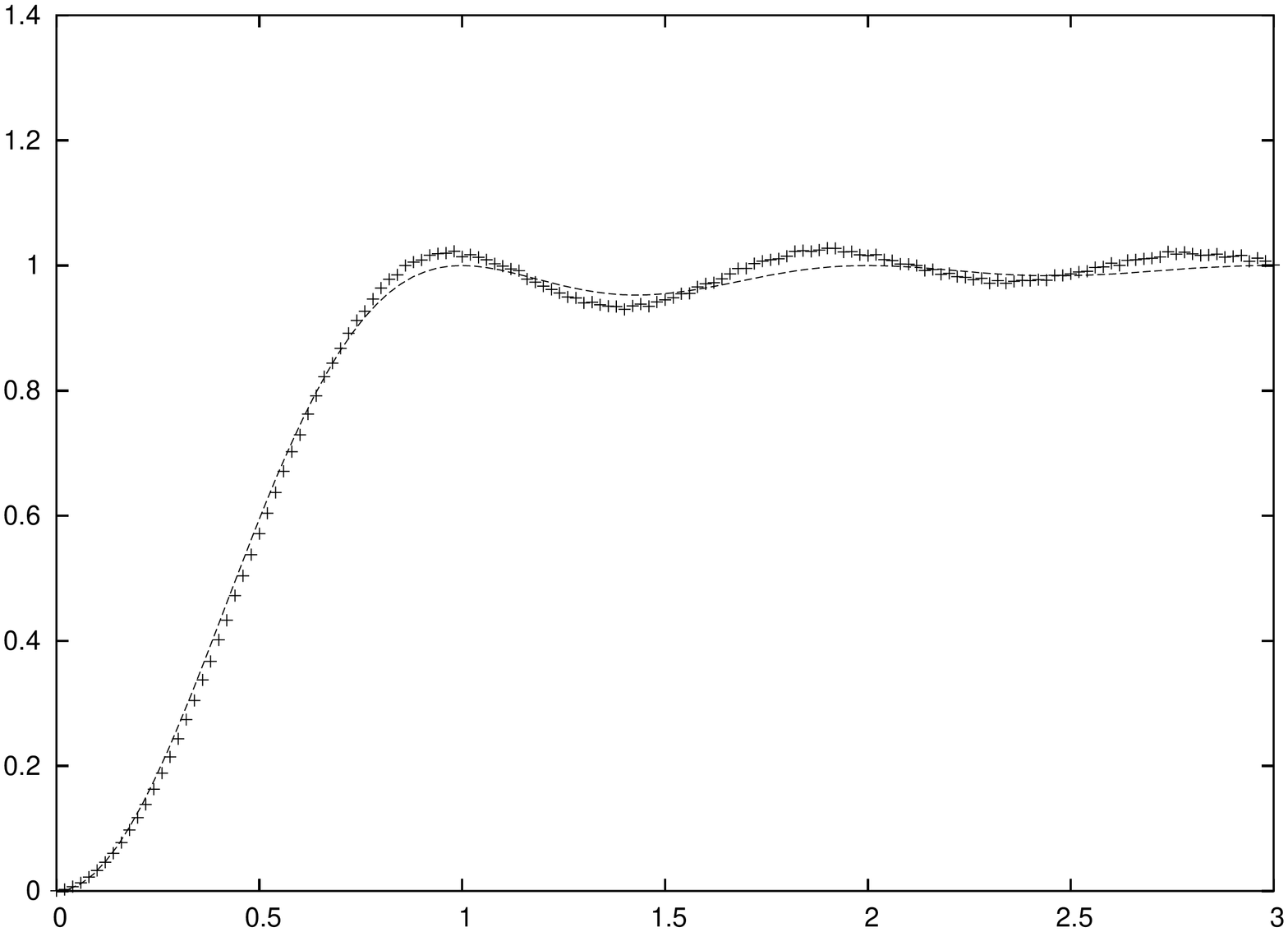,width=2in,angle=-90}
            \psfig{figure=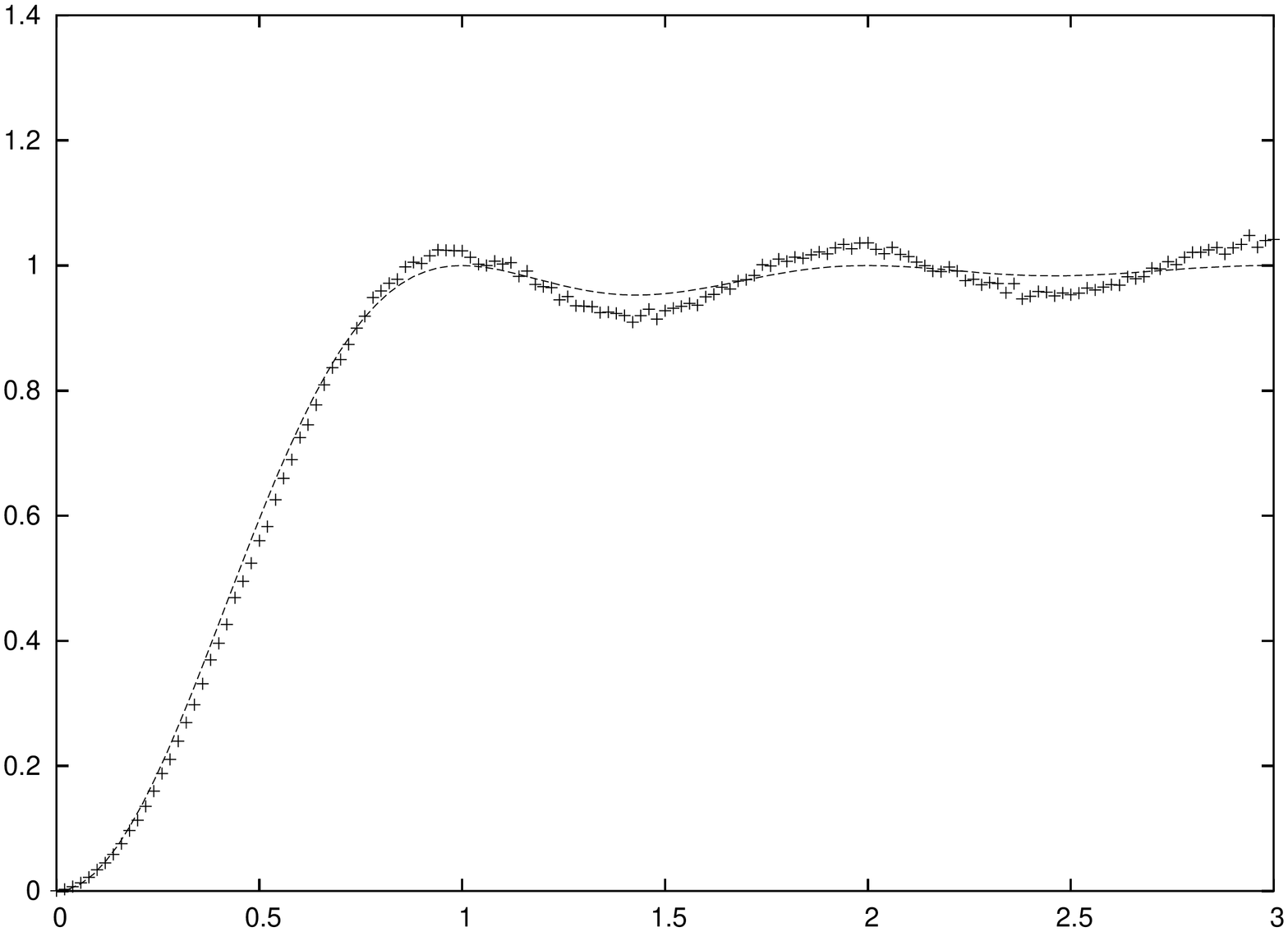,width=2in,angle=-90}
            \psfig{figure=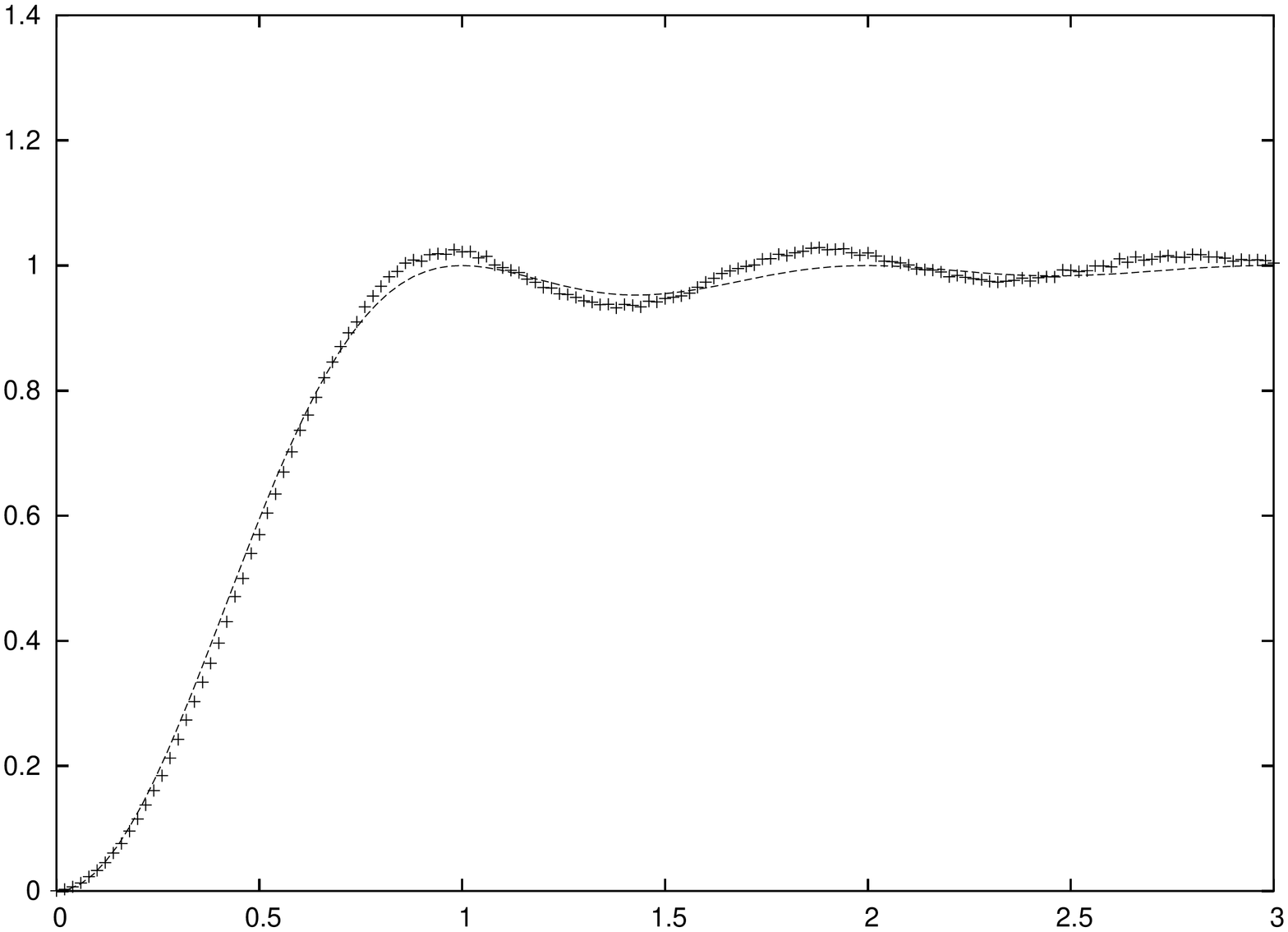,width=2in,angle=-90}
    }
    \centerline{
            \psfig{figure=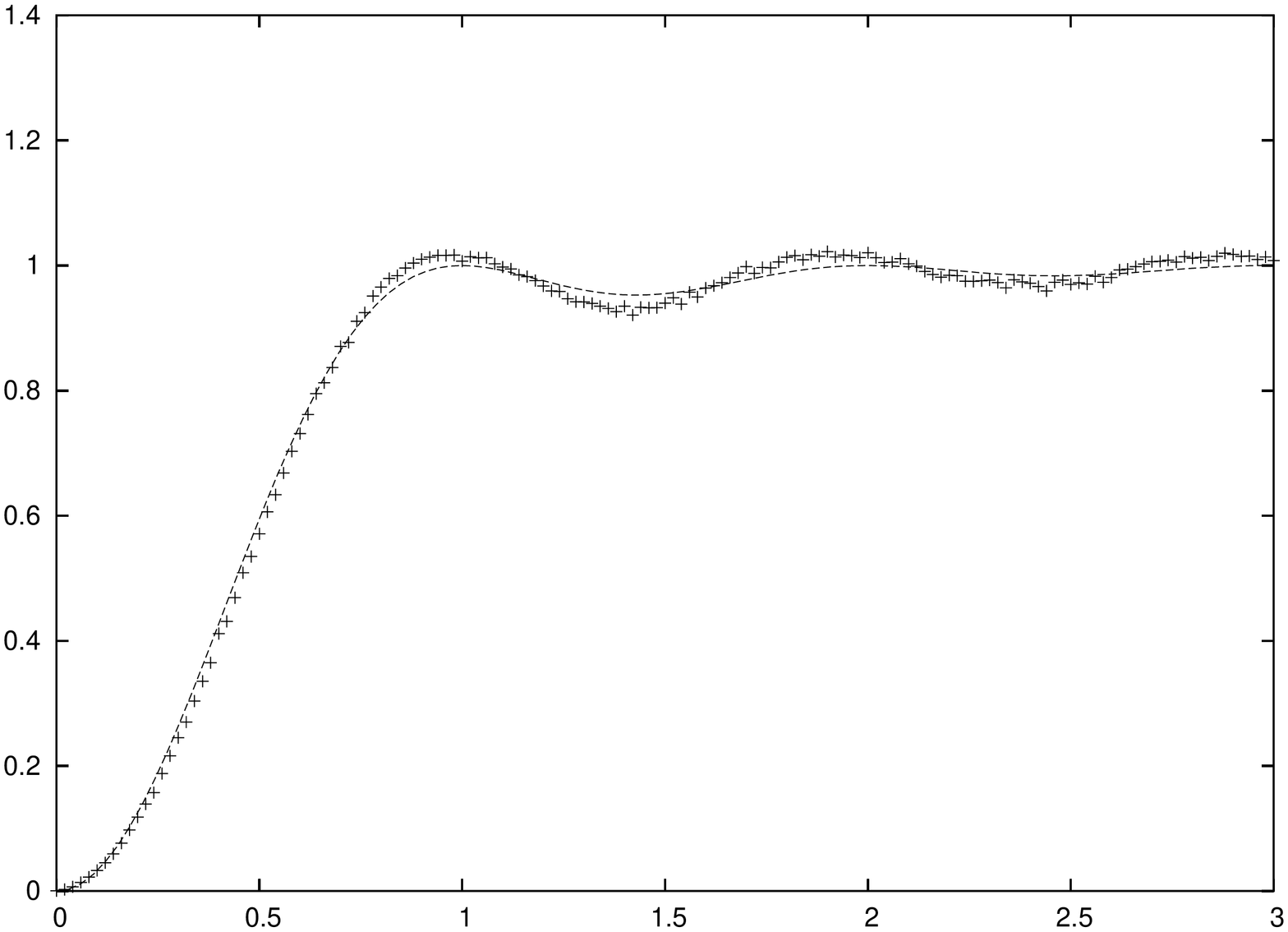,width=2in,angle=-90}
            \psfig{figure=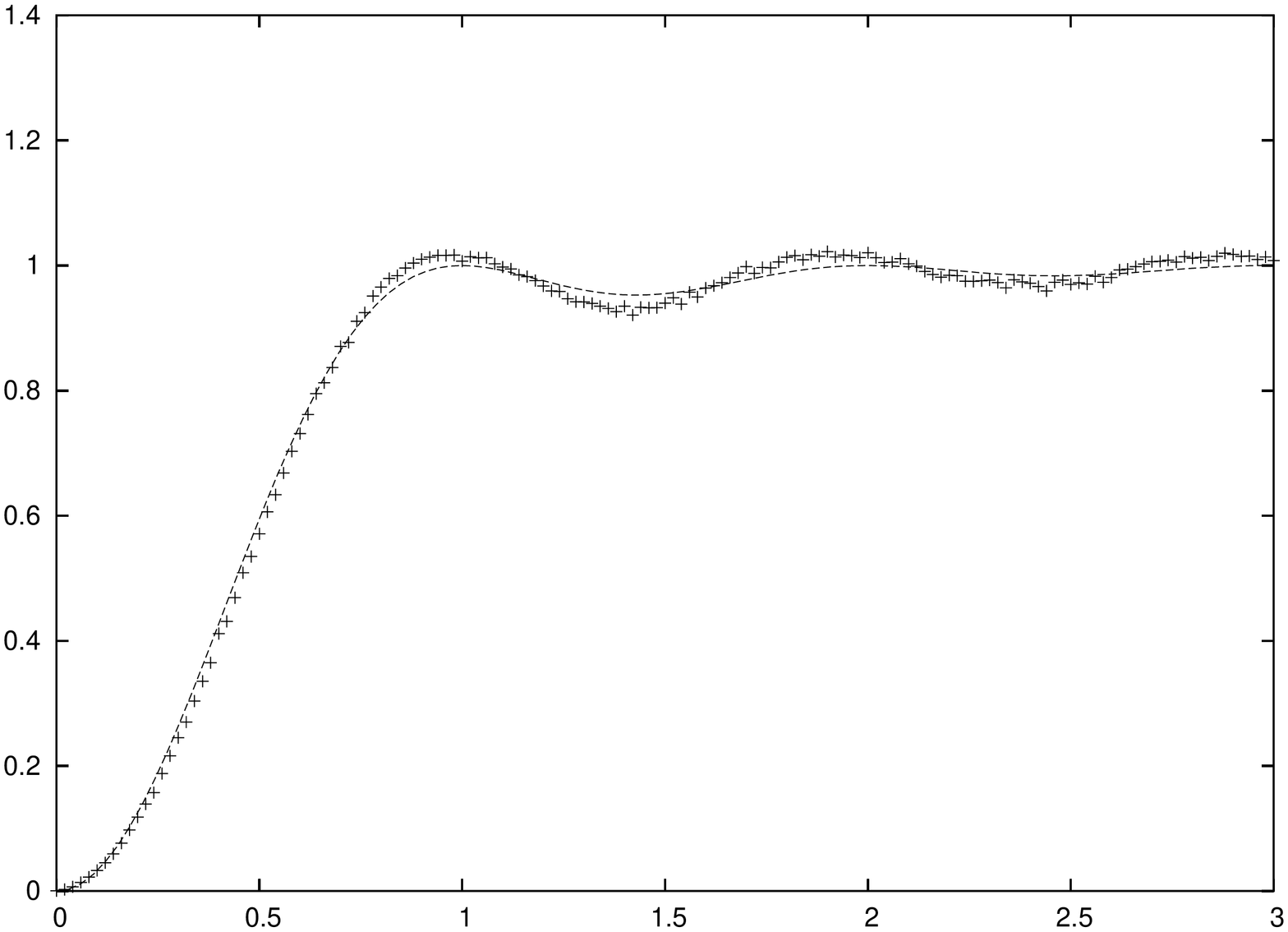,width=2in,angle=-90}
            \psfig{figure=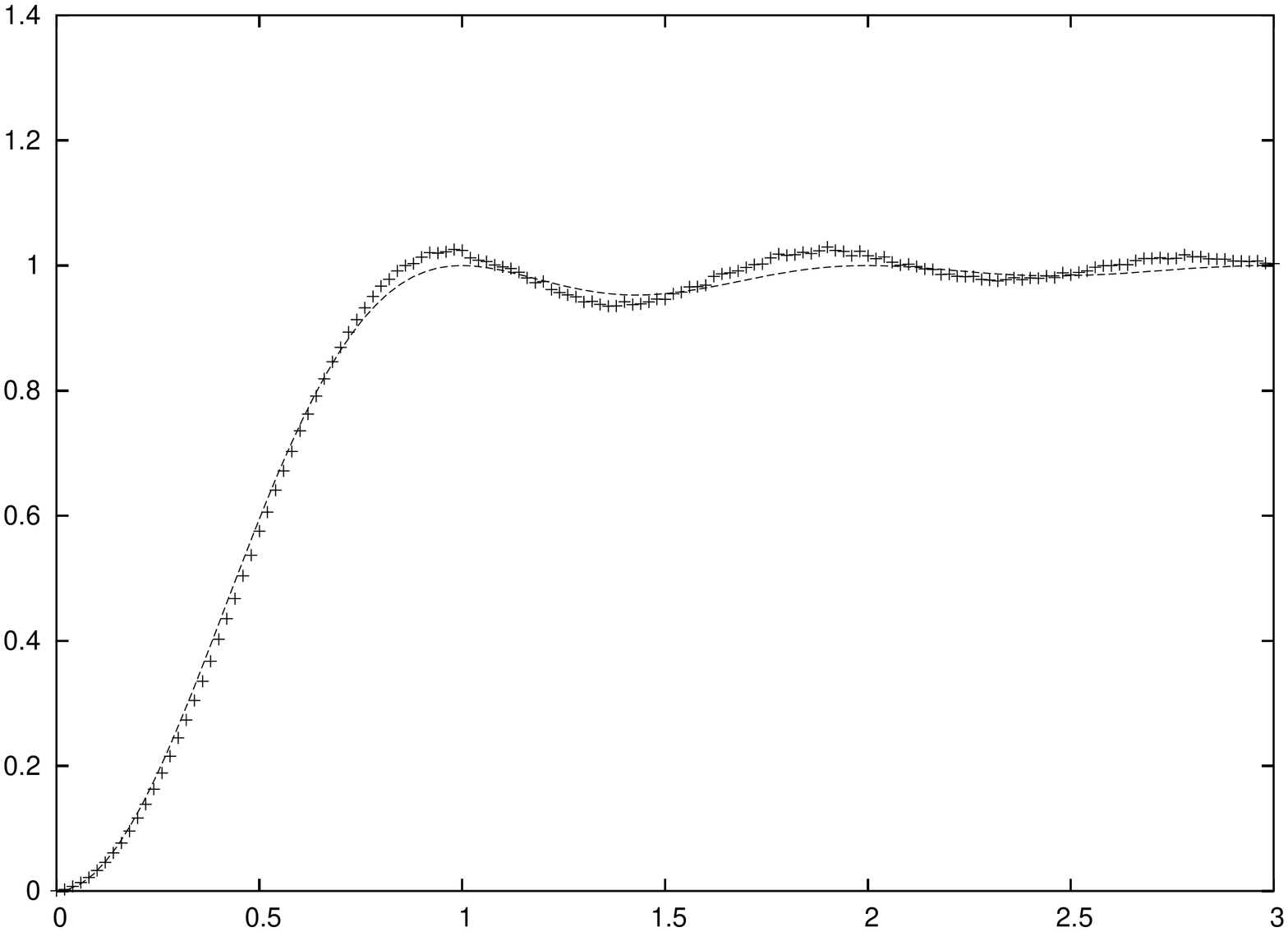,width=2in,angle=-90}
    }
    \centerline{
            \psfig{figure=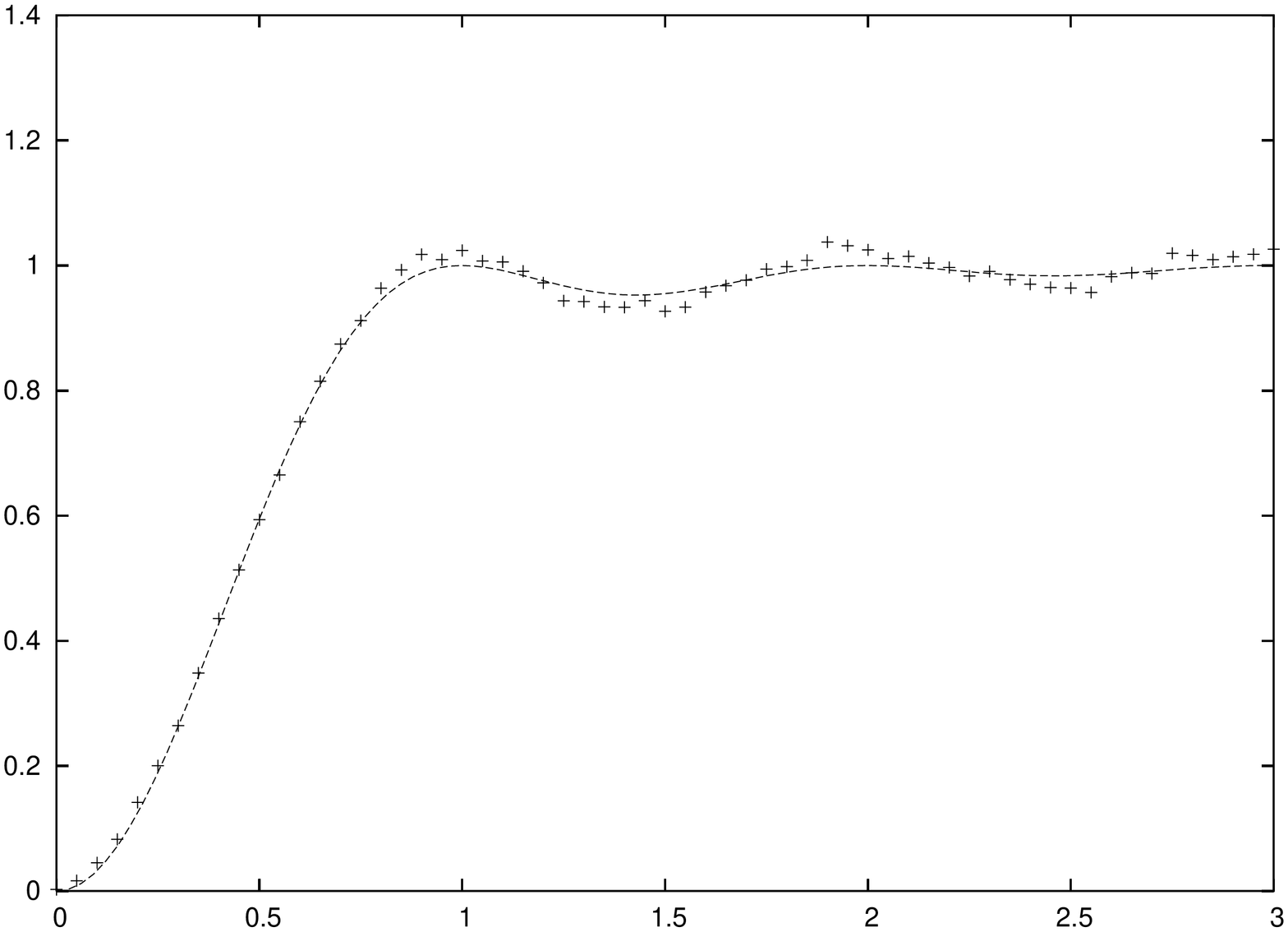,width=2in,angle=-90}
            \psfig{figure=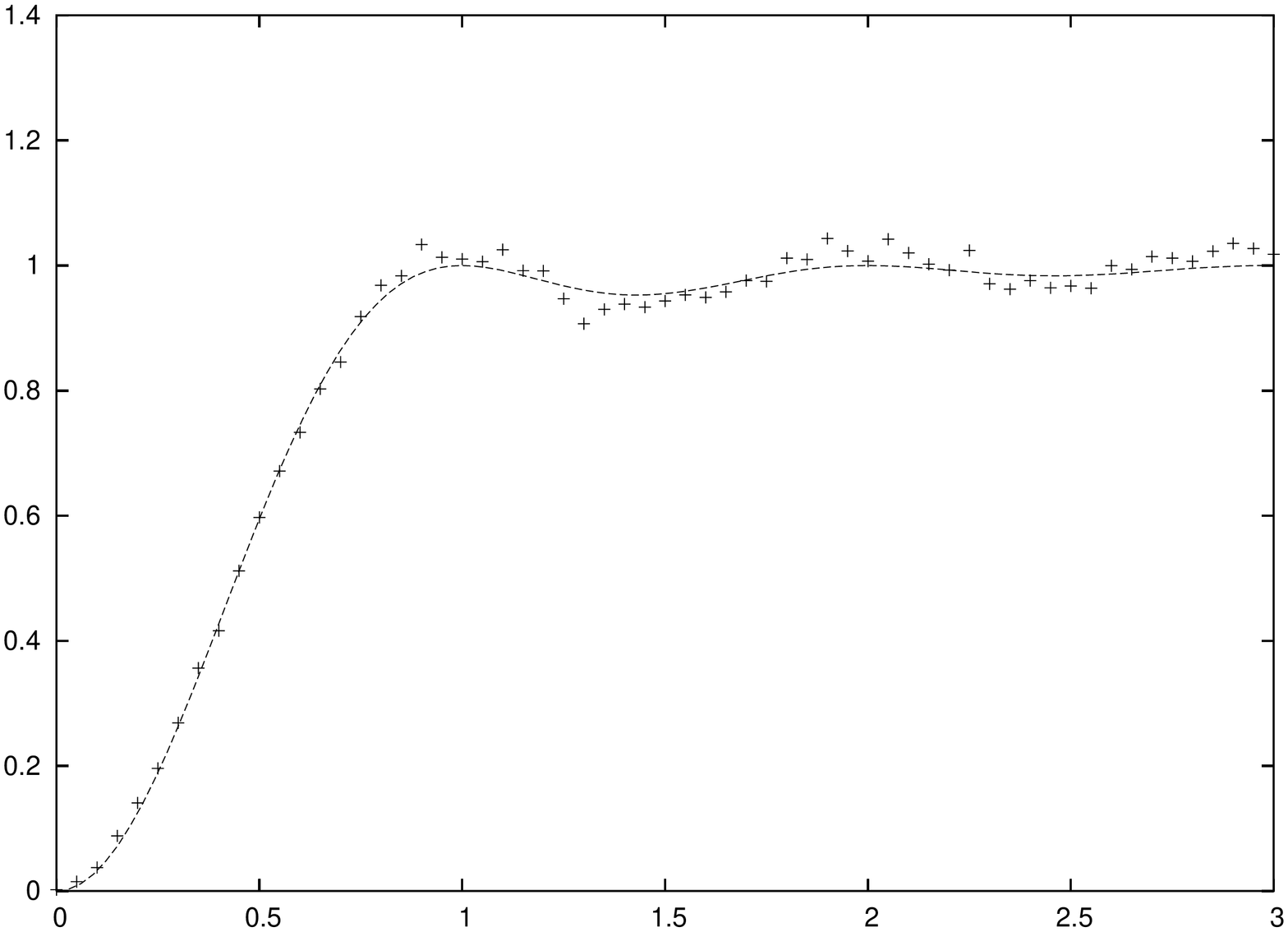,width=2in,angle=-90}
            \psfig{figure=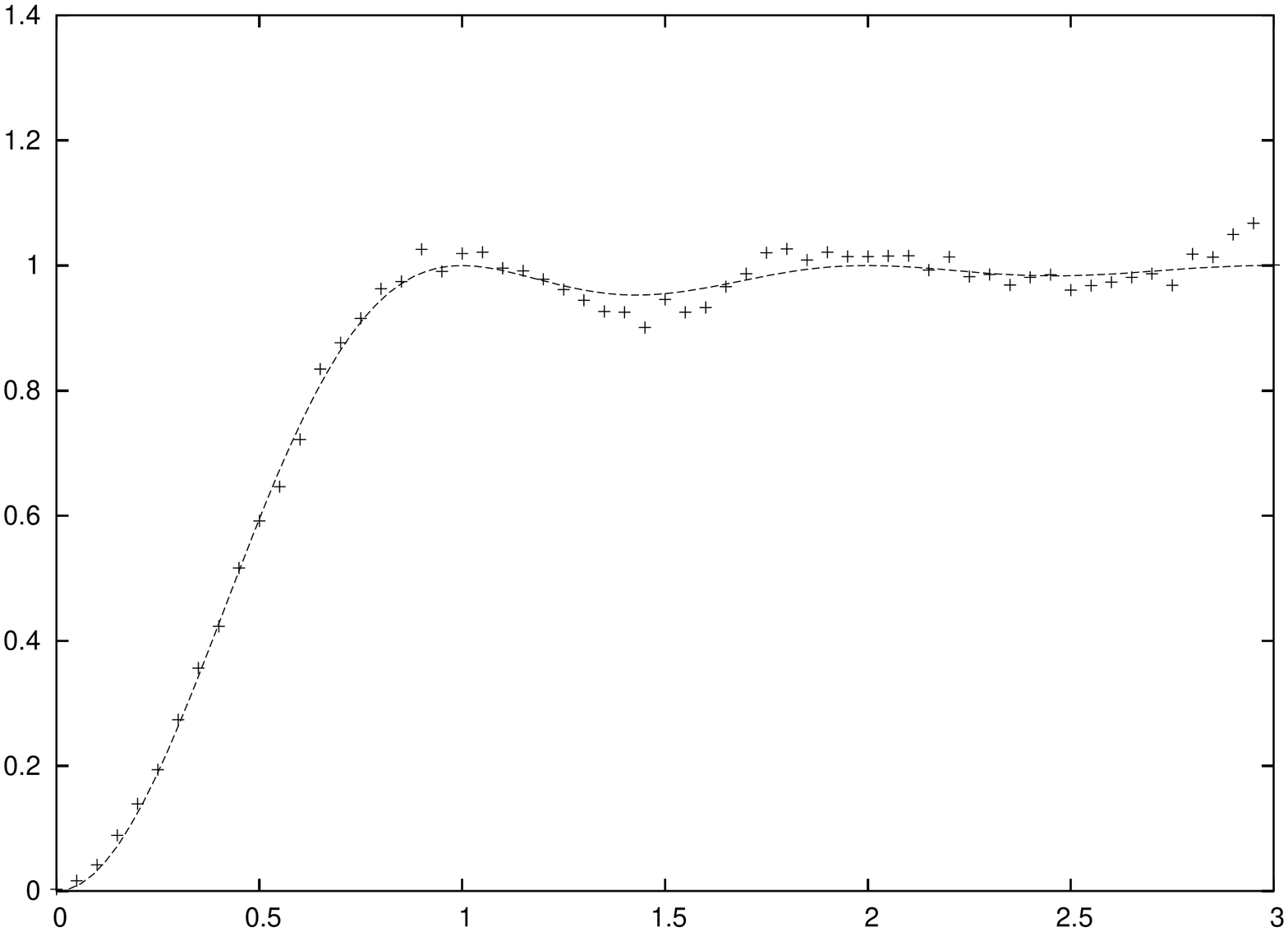,width=2in,angle=-90}
    }
    \centerline{
            \psfig{figure=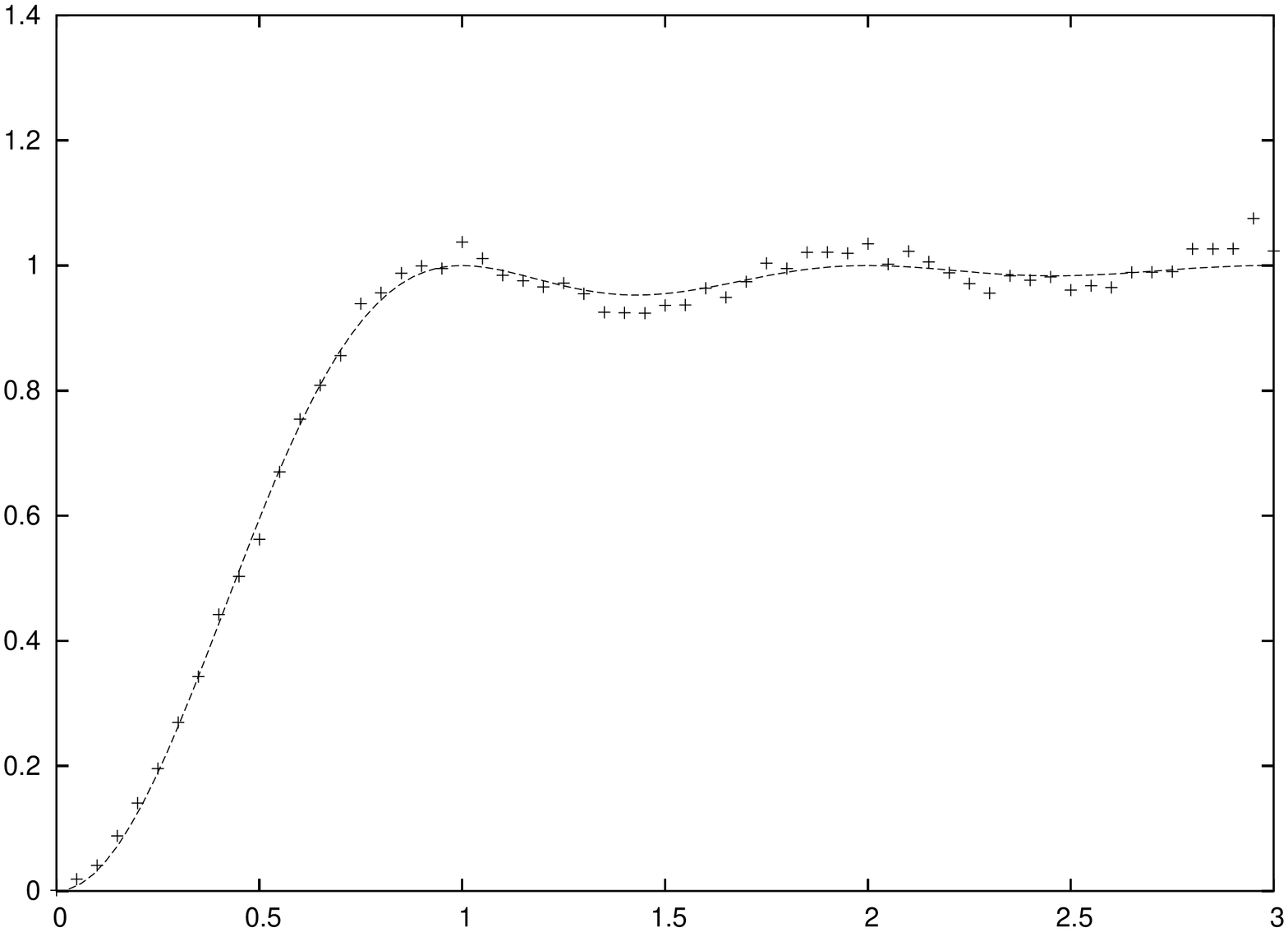,width=2in,angle=-90}
            \psfig{figure=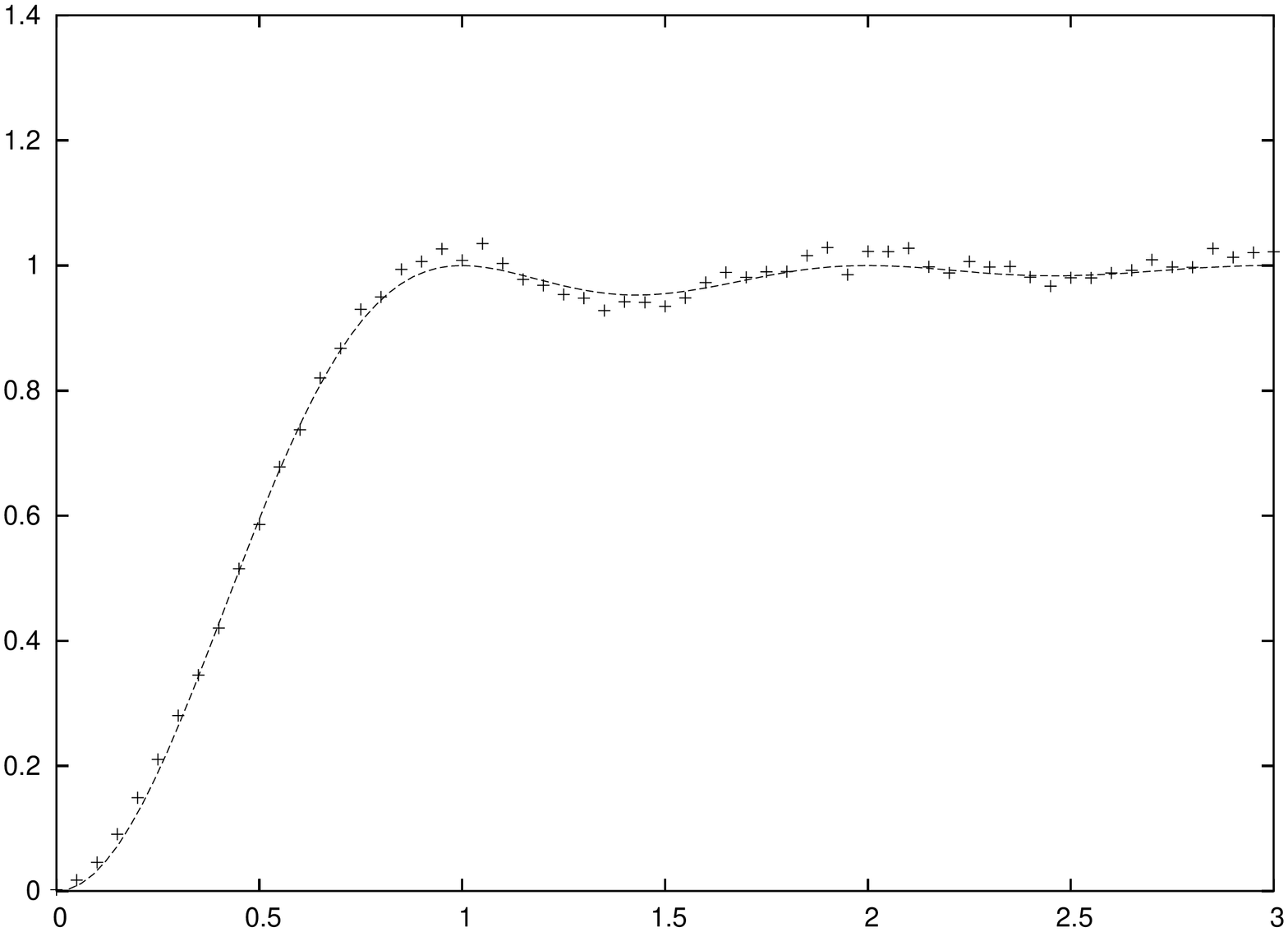,width=2in,angle=-90}
            \psfig{figure=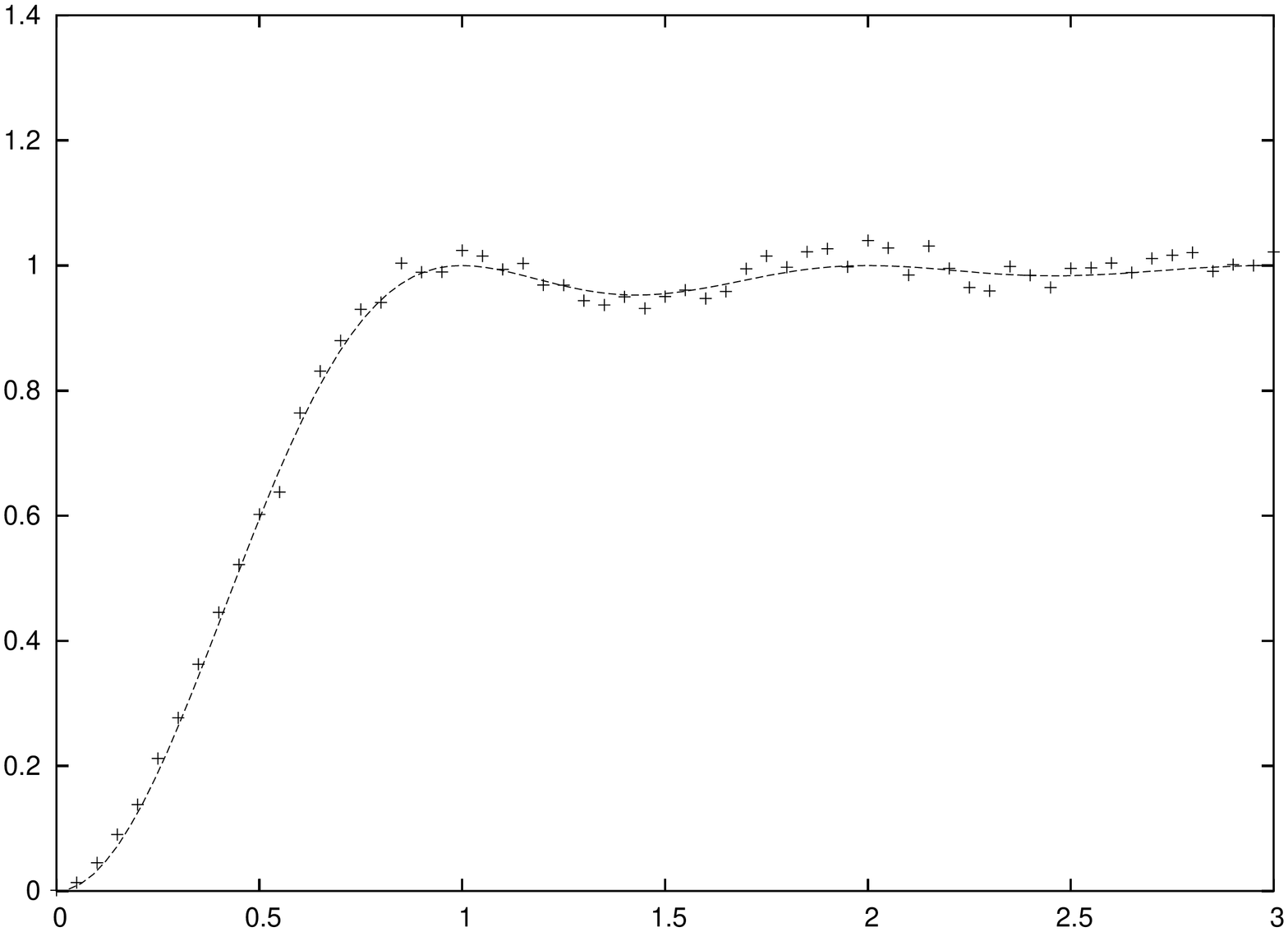,width=2in,angle=-90}
    }
    \caption
    {Pair correlation for zeros of all primitive $L(s,\chi)$, $3 \leq q \leq 17$, 
     the Ramanujan $\tau$ $L$-function, and five elliptic curve $L$-functions} 
    \label{fig:pair correlation}
\end{figure}

The author has tested the
pair correlation conjecture for a number of $L$-functions. Figure~\ref{fig:pair correlation}
depicts the 
same experiment as in Odlyzko's figure, but for 
various Dirichlet $L$-functions and $L$-functions associated to cusp forms.
Altogether there are eighteen graphs.

The first twelve graphs depict the pair
correlation for all primitive Dirichlet $L$-functions, $L(s,\chi)$ for 
conductors $q=3$, $4$, $5$, $7$, $8$, $9$, $11$, $12$, $13$, $15$, $16$, $17$. 
Each graph shows the average pair correlation
for each $q$, i.e. the pair correlation was computed individually for each 
$L(s,\chi)$, and then averaged over $\chi \mod q$. 

In the case of $q=3,4$ there is only one primitive $L$-function for 
either $q$, and approximately five million
zeros were used for each ($4,772,120$ and $5,003,411$ zeros respectively to be precise).
In the case of $q=5,7,8,9,11,12,13,15,16,17$ there are $3,5,2,4,9,1,11,3,4,15$ 
primitive $L$-functions respectively. 
For $q=5,7,8,9,11,12$ either $2,000,000$ zeros
or $1,000,000$ zeros were computed for each $L(s,\chi)$, depending on whether
$\chi$ was real or complex. In the case of $q=16,17$ half as many zeros were computed.

The last six graphs are for $L$-functions associated to cusp forms. The first of these
six shows the pair correlation of the first $284,410$ zeros of the
Ramanujan $\tau$ $L$-function, corresponding to the cusp form of
level one and weight twelve. The next five depict the pair correlation of
the first $100,000$ zeros of
the $L$-functions associated to the elliptic curves of conductors $11,14,15,17,19$.
These last six graphs use larger bins 
since data in these cases is more limited.

The quality of the fit is comparable to what one finds with 
zeros of $\zeta(s)$ up to the same height. See, for example, figures 1 and 3 
in~\cite{O3}.  It would be possible to 
extend the $L(s,\chi)$ computations and obtain data near the $10^{20}$th or
higher zero, at least for reasonably sized $q$. Using the
methods of section~\ref{sec:L algorithms} the time required to
compute $L(1/2+it,\chi)$ is $O(|qt|^{1/2})$, compared to $O(|t|^{1/2})$
for $\zeta(1/2+it)$. Adapting the Odlyzko-Sch\"{o}nhage algorithm
would allow for many evaluations of these $L$-functions at essentially
the cost of a single evaluation.
While such a computation might be manageable for Dirichlet $L$-functions,
it is hopeless for cusp form $L$-functions where the time and
also the number of Dirichlet coefficients required is $O(|N^{1/2}t|)$, i.e. linear
in $t$.  Here $N$ is the conductor of the $L$-function. Using present algorithms and 
hardware, it might be possible to extend these cusp form computations to $t=10^8$ or 
$10^9$. 

Slight care is needed to normalize these zeros correctly as the 
formula for the number of zeros of $L(s)$ depends on the degree of
the $L$-function and on its conductor.
For Dirichlet $L$-functions 
$L(s,\chi)$, $\chi \mod q$, we should normalize its zeros $1/2+i\gamma$
as follows:
\begin{equation}
     \notag
    \tilde{\gamma} = \gamma \frac{\log(|\gamma|q/(2\pi e))}{2\pi}
\end{equation}
For a cusp form $L$-function of conductor $N$, we should take the following
normalization:
\begin{equation}
     \notag
    \tilde{\gamma} = \gamma \frac{\log(|\gamma|N^{1/2}/(2\pi e))}{\pi}
\end{equation}

From a graphical point of view, it is hard to display information
concerning higher order correlations. Instead one can look at a 
statistic that involves knowing~\cite{KS} all the $n$-level correlations for
characteristic functions,
namely the nearest neighbour spacings distribution.

\begin{figure}[htp]
    \centerline{
            \psfig{figure=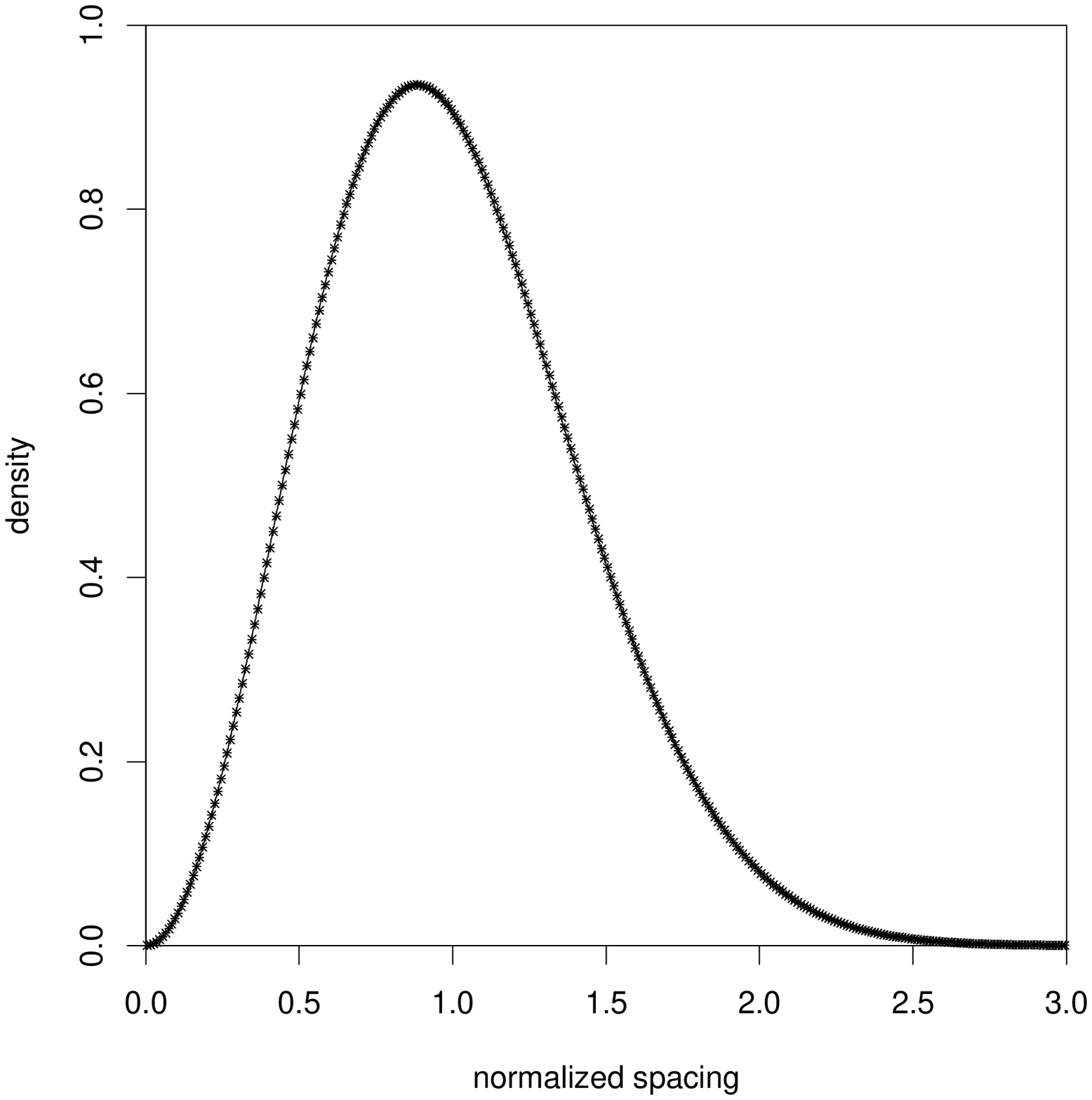,width=3in,angle=-90}
            \psfig{figure=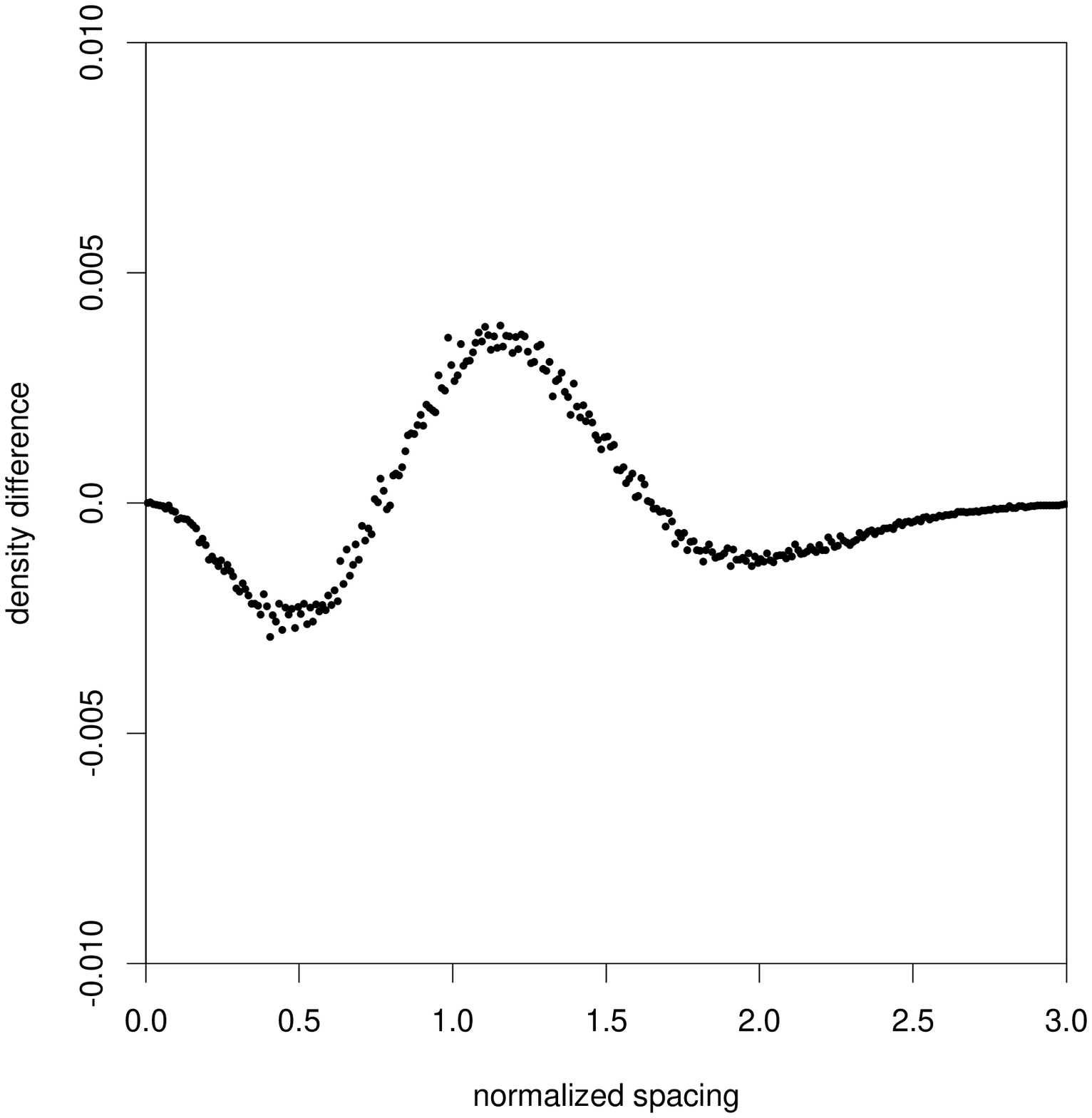,width=3in,angle=-90}
    }
    \caption
    {The first graph shows Odlyzko's nearest neighbour spacings distribution for $2\times 10^8$ zeros of
     $\zeta(s)$ near the $10^{23}$rd zero. The second graph shows the difference he computed between
     the histogram and the predicted density function.
     Recently, Bogomolny, Bohigas and Leboeuf have explained the 
     role of secondary terms in shaping the difference displayed.}
    \label{fig:odlyzko2}
\end{figure}

In Figure~\ref{fig:odlyzko2} we display Odlyzko's
picture for the distribution
of the normalized spacings $\delta_j$ for $2\times 10^8$ zeros of 
$\zeta(s)$ near the $10^{23}$rd zero.
This is computed by breaking up
the $x$-axis into small bins and counting how many $\delta_j$'s fall into
each bin, and then comparing this 
against the nearest neighbour spacings distribution of
the normalized eigenangles of matrices in $\text{U}(N)$, as $N \to \infty$,
again averaged according to Haar measure on $\text{U}(N)$. The density
function for this distribution is given~\cite{M} as
$$
    \frac{d^2}{dt^2} \prod_n (1-\lambda_n(t))
$$
where $\lambda_n(t)$ are the eigenfunctions of the integral operator
\begin{equation}
    \label{eq:integral operator}
    \lambda(t) f(x) = \int_{-1}^1 \frac{\sin(\pi t (x-y))}{\pi(x-y)}f(y) dy,
\end{equation}
sorted according to $1 \geq \lambda_0(t) \geq \lambda_1(t) \geq \ldots \geq 0$.
See~\cite{O3} for a description of how the density function can be computed.

\begin{figure}[htp]
    \centerline{
            \psfig{figure=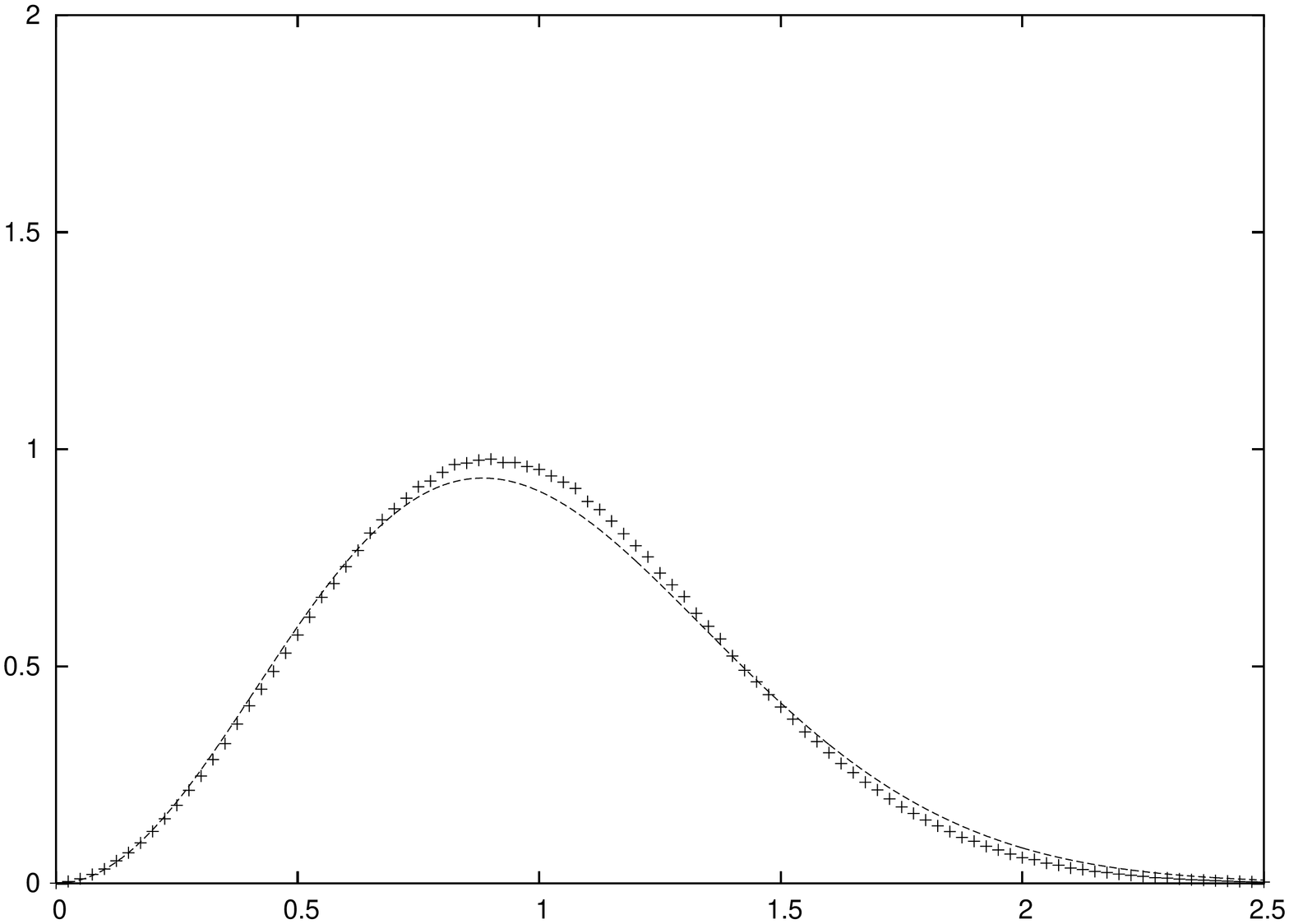,width=3in,height=1.9in,angle=-90}
            \psfig{figure=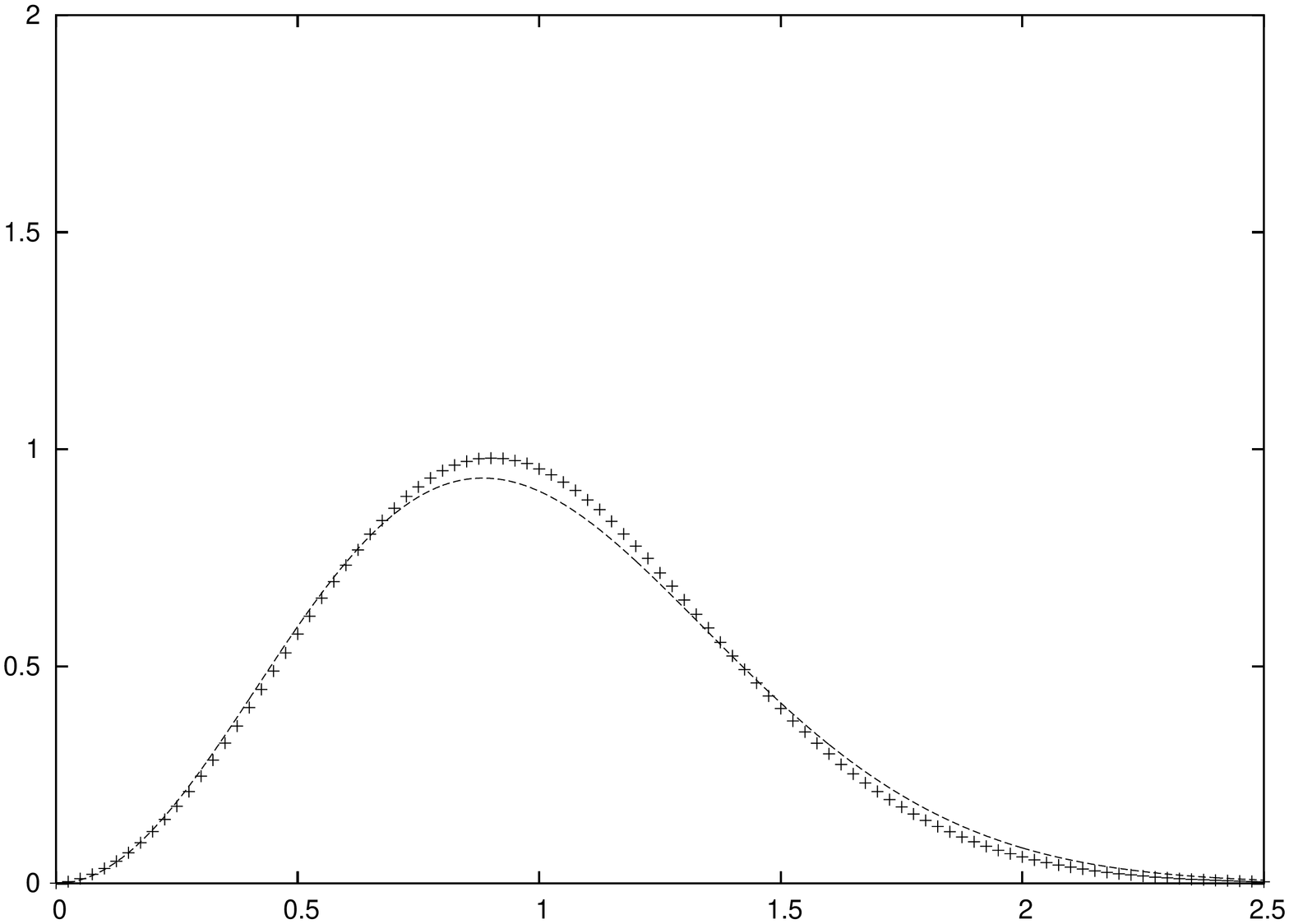,width=3in,height=1.9in,angle=-90}
    }
    \centerline{
            \psfig{figure=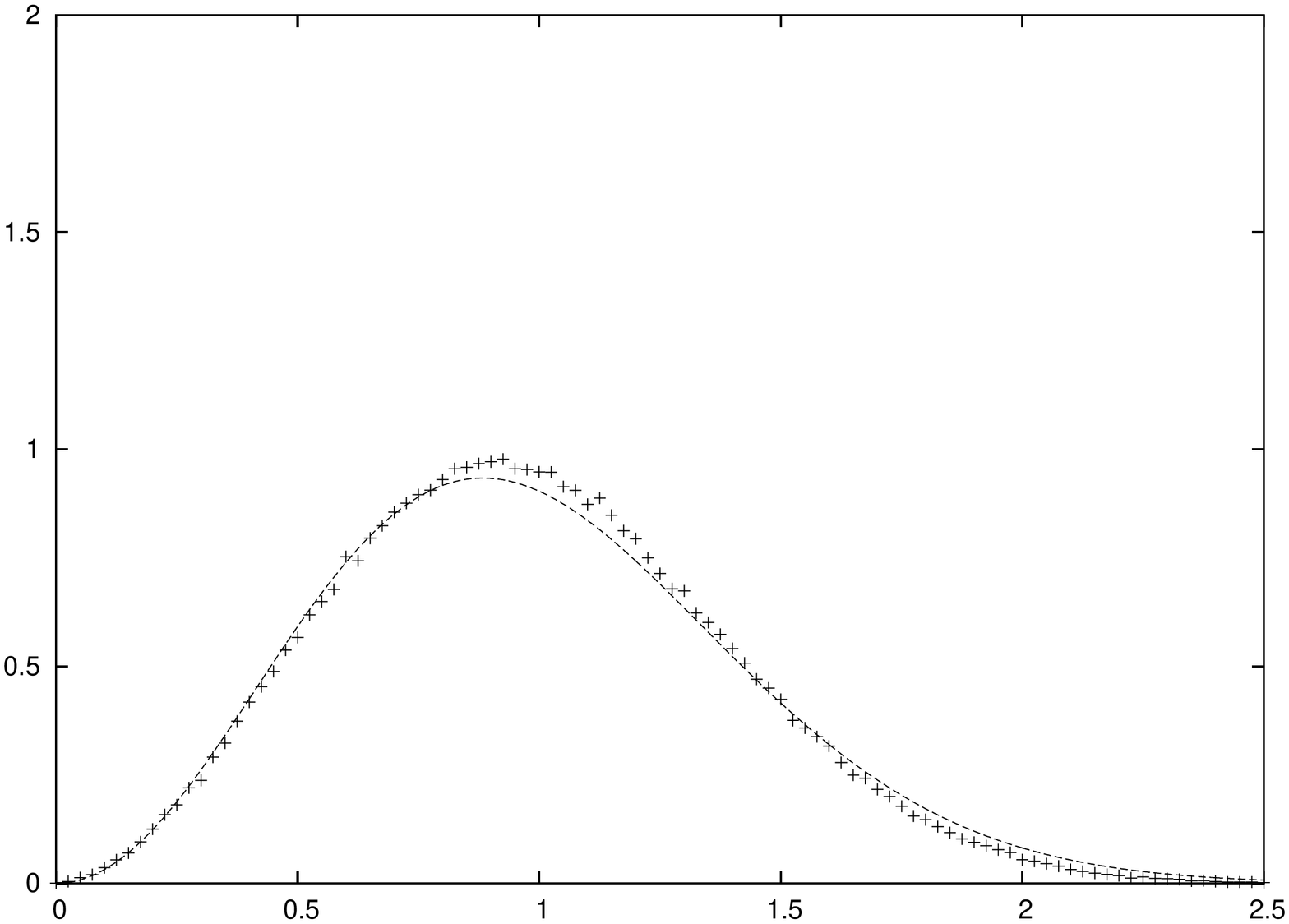,width=3in,height=1.9in,angle=-90}
            \psfig{figure=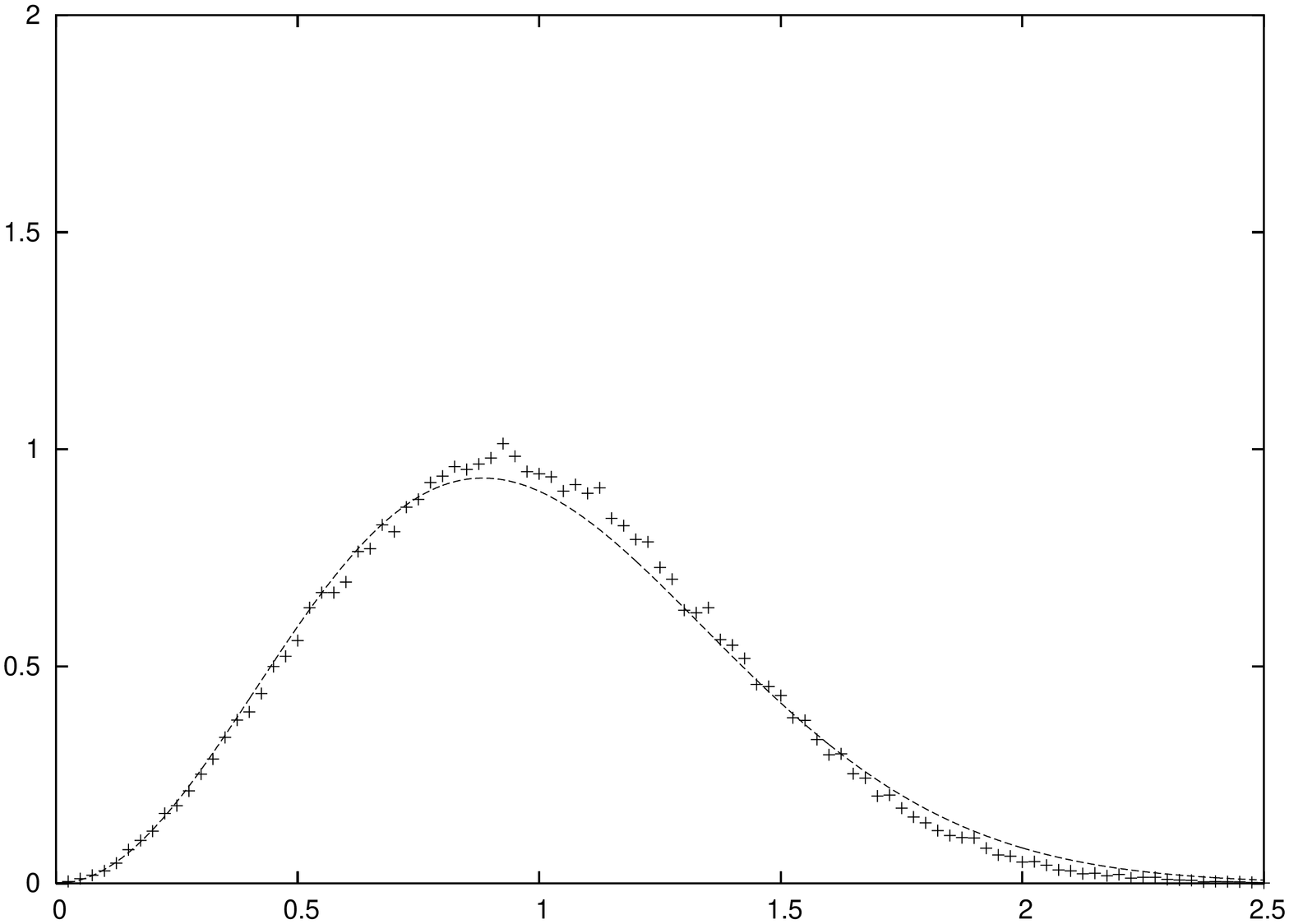,width=3in,height=1.9in,angle=-90}
    }
    \centerline{
            \psfig{figure=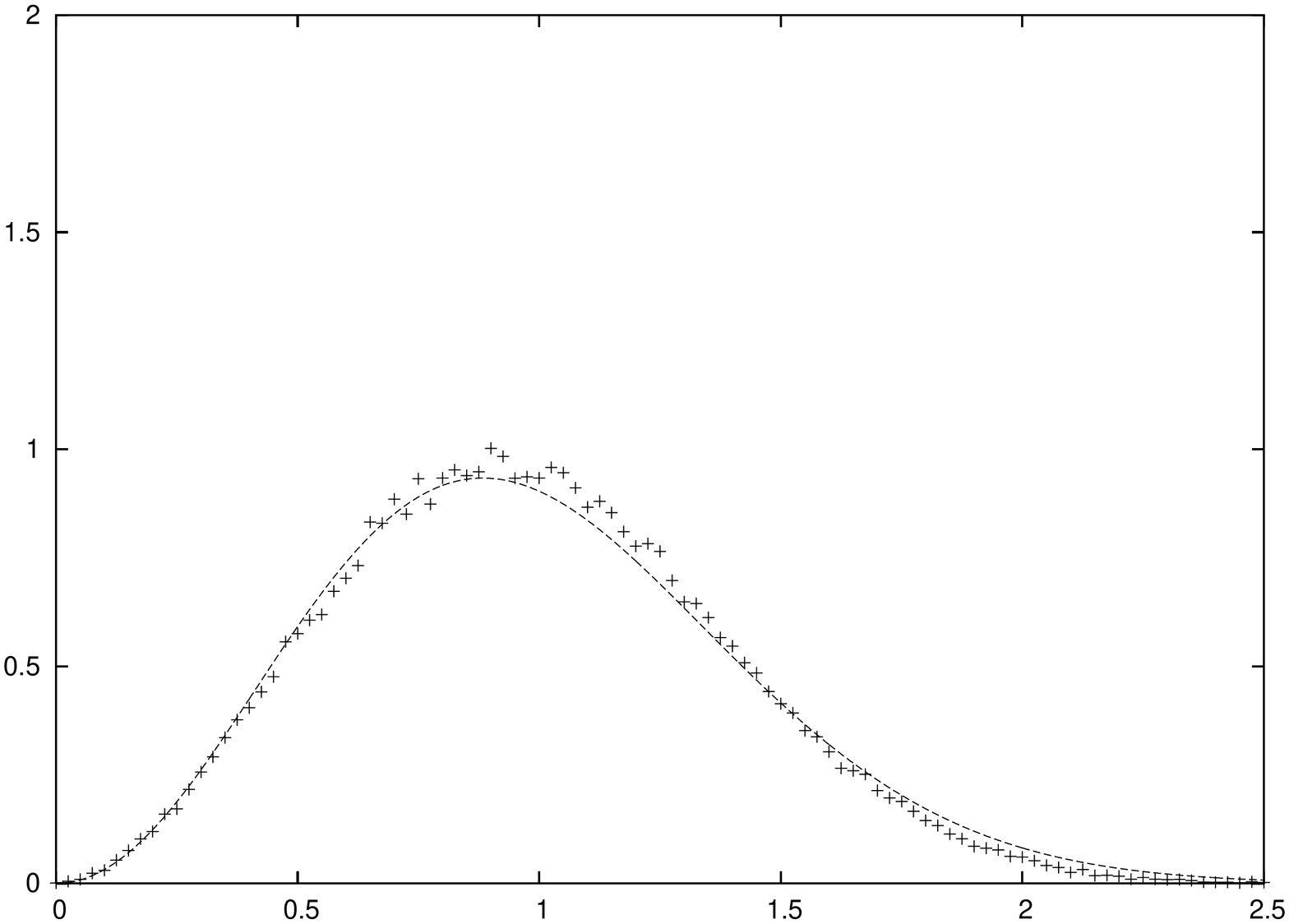,width=3in,height=1.9in,angle=-90}
            \psfig{figure=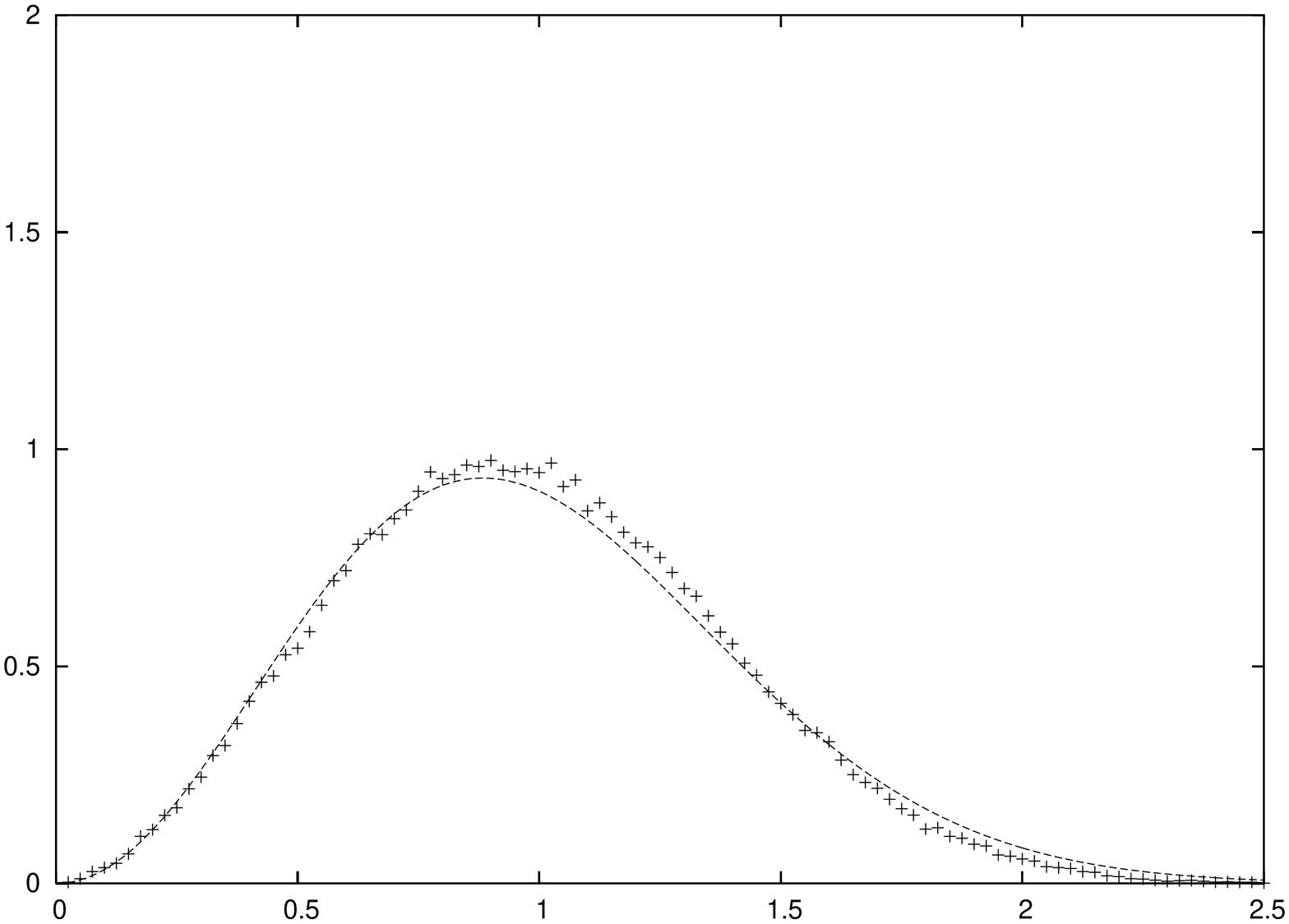,width=3in,height=1.9in,angle=-90}
    }
    \centerline{
            \psfig{figure=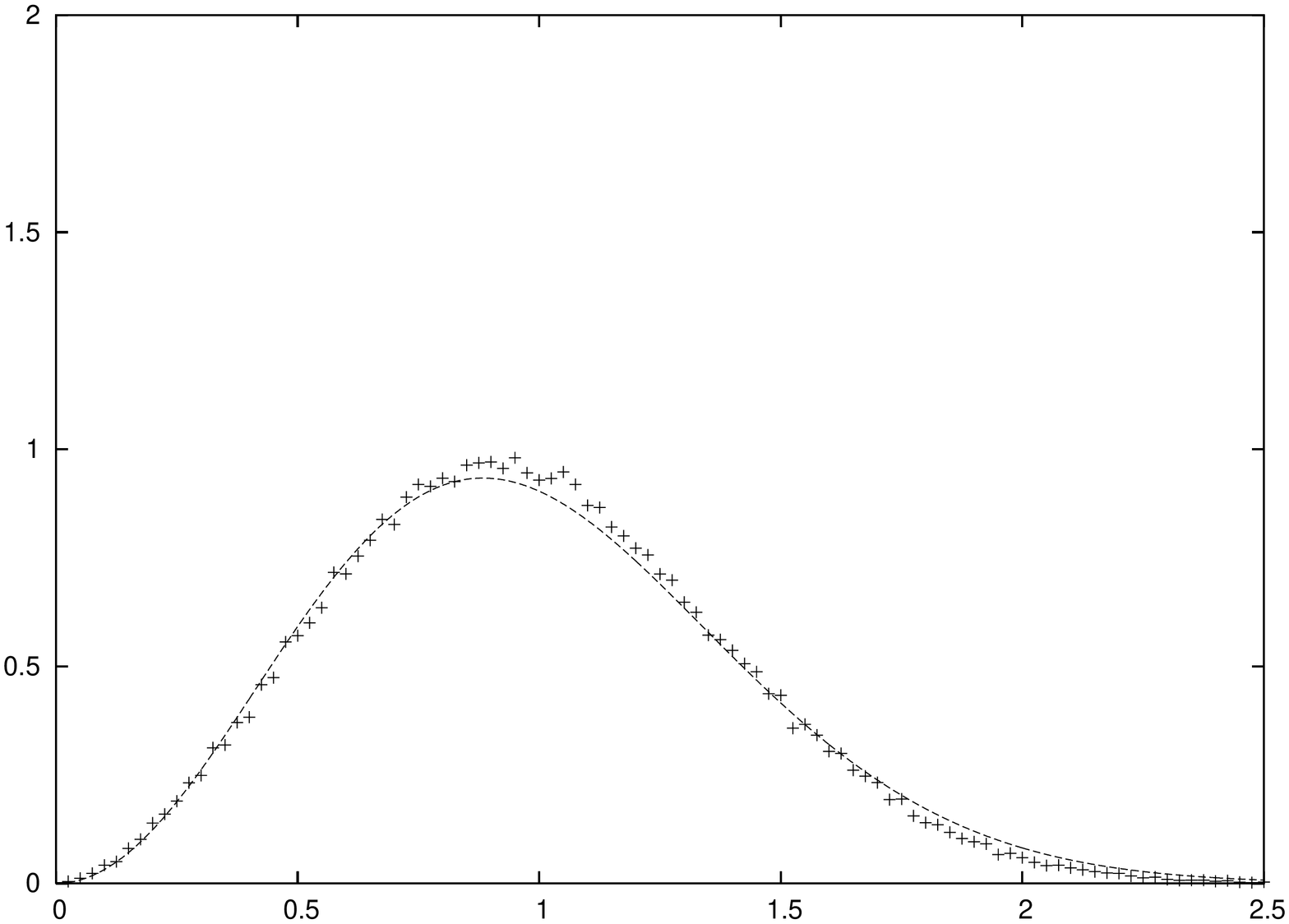,width=3in,height=1.9in,angle=-90}
            \psfig{figure=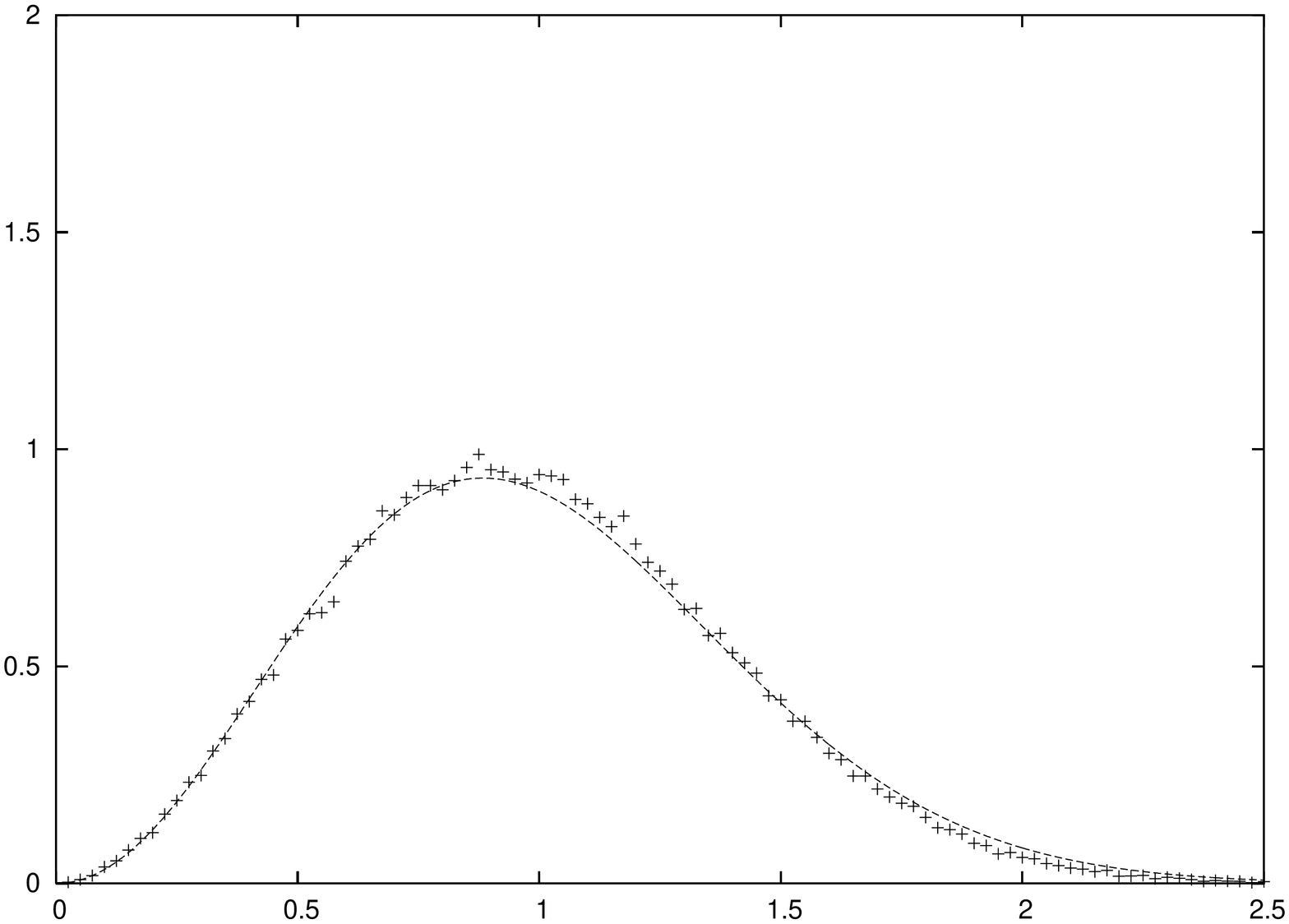,width=3in,height=1.9in,angle=-90}
    }
    \caption
    {Nearest neighbour spacings distribution for several Dirichlet and cusp form
     $L$-functions. The first is for $L(s,\chi)$, $q=3$. 
     The second is the average nearest neighbour spacing for all primitive $L(s,\chi)$,
     $3 \leq q \leq 19$. The last six are for the Ramanujan $\tau$ $L$-function, and
     five $L$-functions associated to elliptic curves.}
    \label{fig:nearest neighbour}
\end{figure}

In Figure~\ref{fig:nearest neighbour} we display the nearest neighbour spacings 
distribution for the sets of zeros described above, namely millions of zeros
of primitive $L(s,\chi)$, with conductors $3 \leq q \leq 17$, and hundreds of thousands
of zeros of six cusp form $L$-functions. We also
depict the nearest neighbour spacings for the first $500,000$ zeros of
each of the 16 primitive $L(s,\chi)$ with $\chi \mod 19$ complex, and $1,000,000$
zeros for the one primitive real $\chi \mod 19$.

Eight graphs are displayed. The first is for the $4,772,120$ zeros of
$L(s,\chi)$, $\chi \mod 3$. The second one depicts the average spacings distribution
for all 76 primitive $L(s,\chi)$, $\chi \mod q$ with $3 \leq q \leq 19$, i.e.
the spacings distribution was computed individually for each of these
$L$-functions and then averaged. The next six graphs show the spacings distribution for
the Ramanujan $\tau$ $L$-function, and
the $L$-functions associated to the elliptic curves of conductors $11,14,15,17,19$.
Again, the fit is comparable to the fit one gets with the same number of zeros
of $\zeta(s)$.

\subsection{Density of zeros}

Rather than look at statistics of a single $L$-function, we can form statistics involving
a collection of $L$-functions. This has the advantage of allowing us to study the behaviour
of our collection near the critical point where specific information about the collection
may be revealed. This idea was formulated by
Katz and Sarnak~\cite{KS}~\cite{KS2} who studied function field zeta functions
and conjectured that the various classical compact groups should be relevant
to questions about $L$-functions. 

While the eigenvalues of matrices in all the classical compact groups share,
on average, the same limiting correlations and spacings distributions, their 
characteristic polynomials do exhibit distinct behaviour near the point $z=1$.
Using the idea that the unit circle for characteristic
polynomials in the classical compact groups correponds to the critical line, with the point
$z=1$ on the unit circle corresponding to the critical point, Katz and Sarnak were
led to formulate conjectures regarding the density of zeros near the critical
point for various collections of $L$-functions. This is detailed in 
section~\ref{sec:density2} below.

\begin{figure}[htp]
    \centerline{
            \psfig{figure=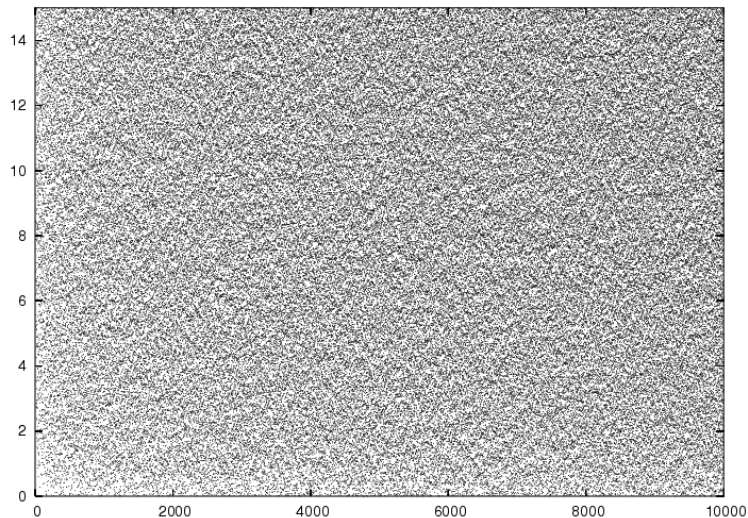,width=4in,angle=-90}
    }
    \caption
    {Zeros of $L(s,\chi)$ with $\chi$ a generic non-real primitive Dirichlet character for the
     modulus $q$, with $5 \leq q \leq 10000$. The horizontal axis is $q$ and, for each
     $L(s,\chi)$, the imaginary parts of its zeros up to height $15$ are listed. 
    }
    \label{fig:zeros complex chi}
\end{figure}

The fact that different families of $L$-funtions exhibit distinct behaviour near the 
critical point is illustrated in figure \ref{fig:zeros complex chi}.
This plot depicts the imaginary parts of the zeros 
of many $L(s,\chi)$ with $\chi$ a generic non-real primitive Dirichlet character for the 
modulus $q$, with $5 \leq q \leq 10000$. 
Other than the fact that, at a fixed height, the zeros become more dense proportionally to 
$\log{q}$, the zeros appear to be uniformly dense.

\begin{figure}[htp]
    \centerline{
            \psfig{figure=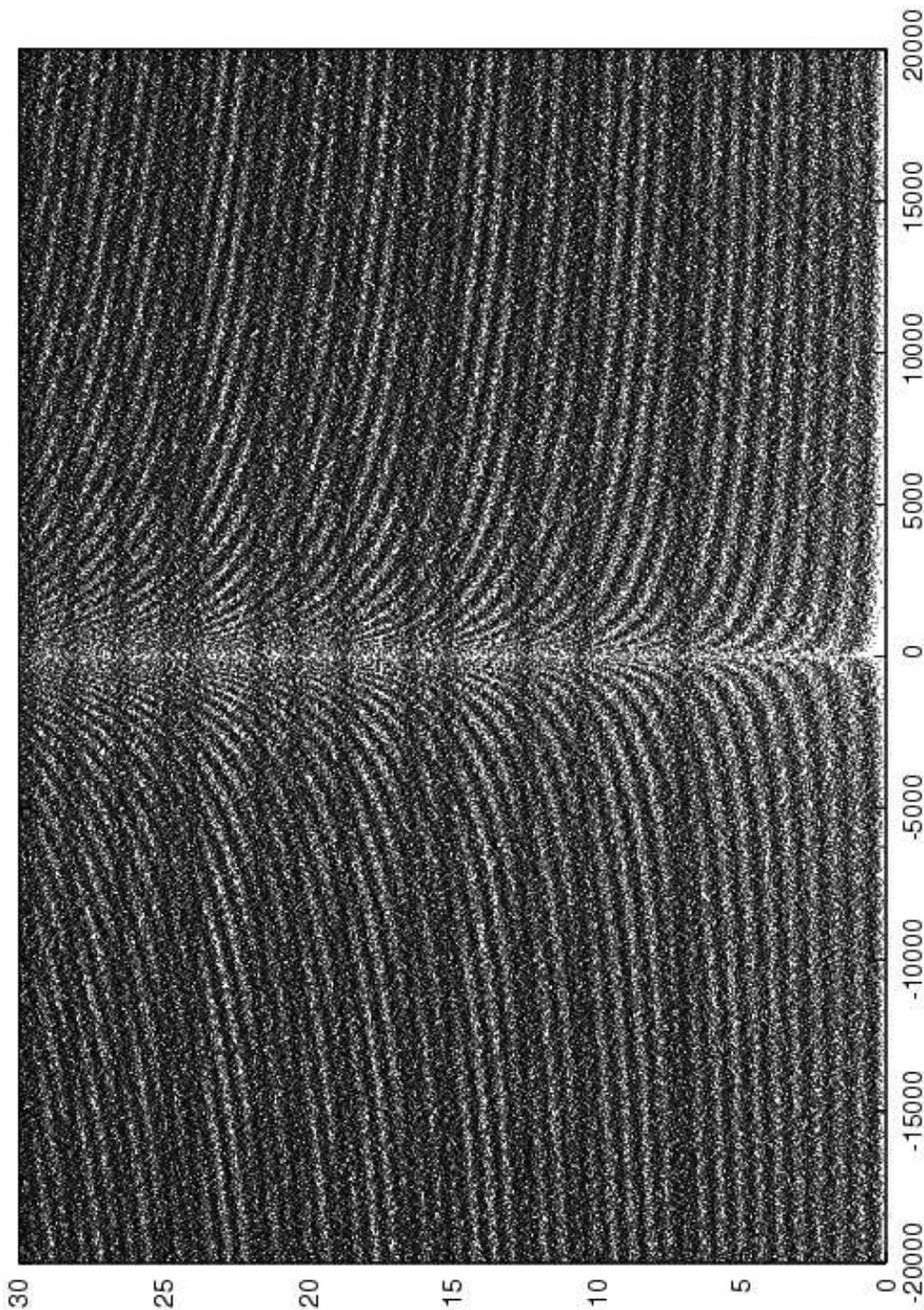,width=5.5in,angle=0}
    }
    \caption
    {Zeros of $L(s,\chi_d)$ with $\chi_d(n)=\left(\frac{d}{n} \right)$,
     the Kronecker symbol. We restrict $d$
     to fundamental discriminants $-20000< d < 20000$.
     The horizontal axis is $d$ and, for each
     $L(s,\chi_d)$, the imaginary parts of its zeros up to height $30$ are listed.
     A higher resolution image can be obtained from the author's webpage under
     `Publications'.
    }
    \label{fig:zeros real chi}
\end{figure}

This contrasts sharply with the plot in figure~\ref{fig:zeros real chi}
which depicts the zeros of
$L(s,\chi_d)$ where $\chi_d$ is a real primitive character (the Kronecker symbol), 
and $d$ ranges over 
fundamental discriminants with $-20000<d<20000$. Here we see the density of zeros 
fluctuating as one moves away from the real axis. 

Other features can be seen in the plot.
First, from the white band near the $x$-axis 
we notice that the lowest zero for each $L(s,\chi_d)$ tends to stay away from
the critical point. We can also see the effect of secondary terms on this
repulsion. The lowest zero for $d>0$ tends to be higher than the lowest zero for
$d<0$. This turns out to be related to the fact that the $\Gamma$-factor in the 
functional equation 
for $L(s,\chi_d)$ is $\Gamma(s/2)$ if $d>0$, but is $\Gamma((s+1)/2)$ when $d<0$.

\begin{figure}[htp]
    \centerline{
            \psfig{figure=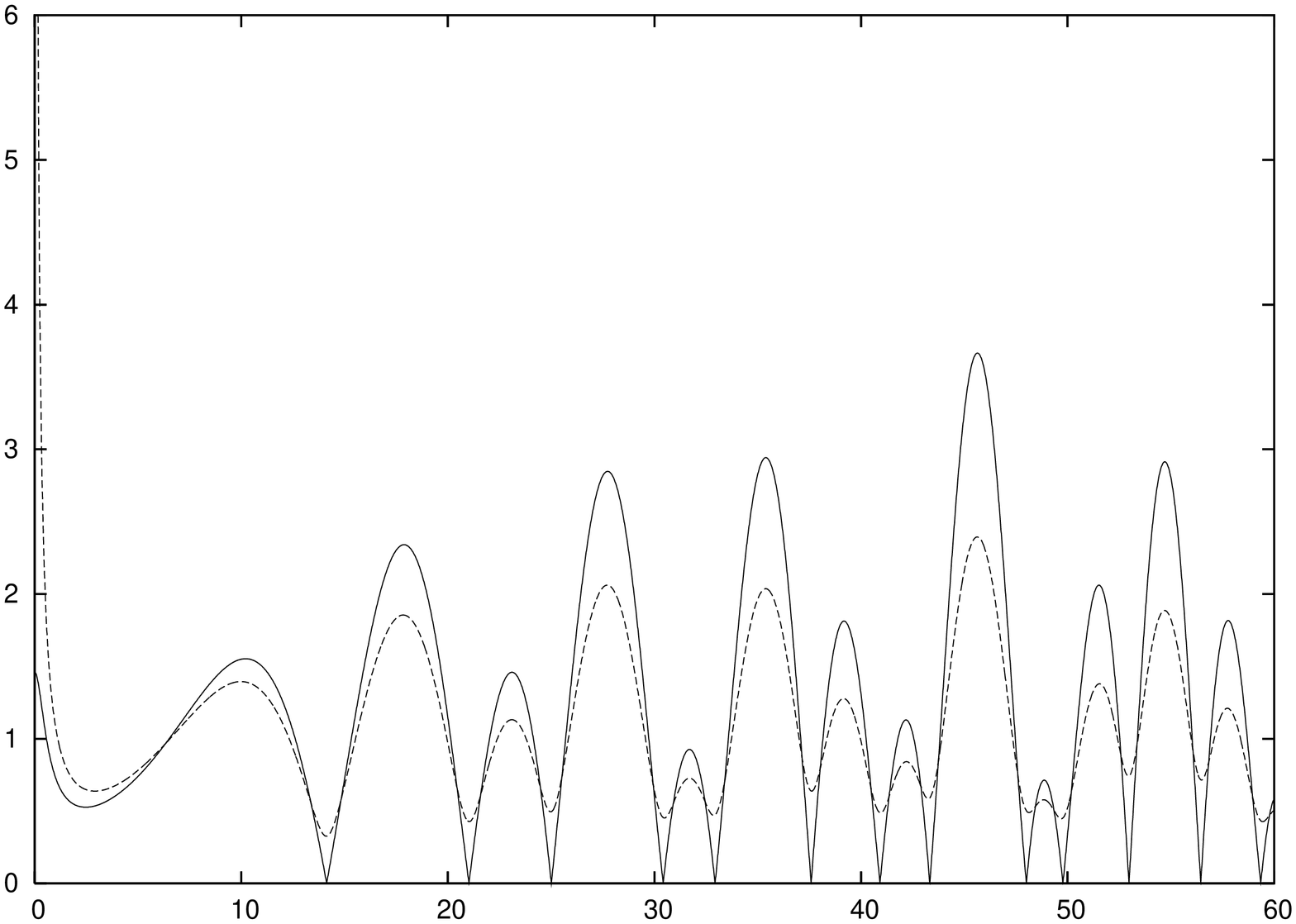,width=4in,angle=-90}
    }
    \caption
    {
       A graph illustrating that, at least initially, the minima of $|\zeta(1+iy)|$
       occur very close the zeros of $|\zeta(1/2+iy)|$. The dashed
       line is the graph of the former while the solid line is the graph of the 
       the latter.
    }
    \label{fig:zeta half and one line}
\end{figure}

We can also see slightly darker regions appearing in horizontal strips.
The first one occurs roughly at height $7.$, half the height of
the first zero of $\zeta(s)$.
These horizontal strips are due to secondary terms in the density of zeros for this 
collection of $L$-functions
which include~\cite{Ru3}~\cite{CS} a term that is proportional to
$$
    \Re\frac{\zeta'(1+2it)}{\zeta(1+2it)}.
$$ 
This is large when $\zeta(1+2it)$ is small. Surprisingly, $\zeta(1+iy)$ and
$\zeta(1/2+iy)$ track each other very closely, see figure \ref{fig:zeta half and one line}, 
and the minima of 
$|\zeta(1+iy)|$ appear close to the zeros of $\zeta(1/2+iy)$.
This is similar to a phenomenon
that occurs when we look at secondary terms in the pair correlation of the
zero of $\zeta(s)$ which also involves $\zeta(s)$ on the one 
line~\cite{BoK}~\cite{BK2}.

\subsubsection{$n$-level density}
\label{sec:density2}

The $n$-level density is used to measure the average density of the 
zeros of a family of $L$-functions or matrices. It is arranged to
be sensitive to the low lying zeros in the family, i.e. those near the
critical point if we are dealing with $L$-functions, and those near
the point $z=1$ on the unit circle if we are dealing with characteristic 
polynomials from the classical compact groups.

Let $A$ be an $N \times N$ matrix in on of the classical compact groups.
Write the eigenvalues of $A$ as $\lambda_j=e^{i \theta_j}$ with
$$
    0 \leq \theta_1 \leq \ldots \leq \theta_N < 2\pi.
$$
Let
$$
    H^{(n)}(A,f)
    =
    \sum_{1 \leq j_1,\ldots,j_n \leq N \atop \text{distinct}}
    f\br{\theta_{j_1} N/(2\pi),\ldots,\theta_{j_n} N/(2\pi)}
$$
with $f: {\mathbb{R}}^n \to {\mathbb{R}}$, bounded, Borel measurable,
and compactly supported. Because of the normalization by $N/(2\pi)$, 
and the assumption that $f$ has compact support,
$H^{(n)}(A,f)$ only depends on the small $\theta_{j}$'s.

Katz and Sarnak~\cite{KS} proved the following family dependent result:
\begin{equation}
    \label{eq:non-universal}
    \lim_{N\to\infty}
    \int_{G(N)}
    H^{(n)}(A,f) dA
    = 
    \int_0^\infty
    \ldots
    \int_0^\infty
    W_G^{(n)}(x) f(x) dx
\end{equation}
for the following families: 
\\ 
\\ 
\\
\begin{tabular}{|c|c|}
$G$                                   & $W_G^{(n)}$ \\ \hline
$\text{U}(N)$,$\text{U}_\kappa(N)$         & $\det\br{K_{0}(x_j,x_k)}_{
                                            {1 \leq j \leq n} \atop
                                            {1 \leq k \leq n}}$ \vspace{.01in}\\ \hline
$\text{USp}(N)$                       & $\det\br{K_{-1}(x_j,x_k)}_{
                                            {1 \leq j \leq n} \atop
                                            {1 \leq k \leq n}}$ \vspace{.01in}\\ \hline
$\text{SO}(2N)$                       & $\det\br{K_{1}(x_j,x_k)}_{
                                            {1 \leq j \leq n} \atop
                                            {1 \leq k \leq n}}$ \vspace{.01in}\\ \hline
$\text{SO}(2N+1)$                     & $\det\br{K_{-1}(x_j,x_k)}_{
                                            {1 \leq j \leq n} \atop
                                            {1 \leq k \leq n}}
                              + \sum_{\nu=1}^n \delta(x_\nu) \det\br{K_{-1}(x_j,x_k)}_{
                                            {1 \leq j\neq\nu \leq n} \atop
                                            {1 \leq k\neq\nu \leq n}}$ \vspace{.01in}\\ \hline
\end{tabular}
with
$$
    K_{\e}(x,y) = \frac{\sin(\pi(x-y))}{\pi(x-y)}
                  +\e \frac{\sin(\pi(x+y))}{\pi(x+y)}.
$$
Here
$$
    \text{U}_\kappa(N)  = \{A \in \text{U}(N): {\det}(A)^\kappa=1\}.
$$
The delta functions in the $\text{SO}(2N+1)$ case
are accounted for by the eigenvalue at $1$.
Removing this zero from~(\ref{eq:non-universal})
yields the same $W_G^{(n)}$ as for USp. 

Let 
$$
   D(X) = \left\{ 
       \text{$d$ a fundamental discriminant : $|d| \leq X$}
   \right\}
$$
and let $\chi_d(n)=\br{\frac{d}{n}}$ be Kronecker's symbol.
Write the non-trivial zeros of $L(s,\chi_d)$ as
$$
    1/2 + i\gam{j}, \hspace{.5in} j= \pm 1, \pm 2, \ldots
$$
sorted by increasing imaginary part, and
\begin{equation}
    \notag
    \gam{-j} = - \gam{j}.
\end{equation}
The author proved~\cite{Ru2} that 
\begin{align}
    \label{eq:Wsp}
    &\lim_{X\to\infty}
    \frac{1}{|D(X)|}
    \sum_{d \in D(X)}
    \sum_{j_i \geq 1 \atop \text{distinct}}
        f\br{l_d\gam{j_1}, l_d\gam{j_2}, \ldots, l_d\gam{j_n}} \notag \\
    &= 
    \int_0^\infty
    \ldots
    \int_0^\infty
    f(x) W_{\text{USp}}^{(n)}(x) dx,
\end{align}
where
\begin{equation}
    \notag
     l_d = \frac{\log(|d|/\pi)}{2\pi}.
\end{equation}
Here,  $f$ is assumed to be smooth, and rapidly decreasing
with $\hat{f}(u_1,\ldots,u_n)$ supported in $\sum_{i=1}^{n} \abs{u_i} <1$.
This generalized the $n=1$ case that had been achieved earlier~\cite{OzS}~\cite{KS2}.
Assuming the Riemann Hypothesis for all $L(s,\chi_d)$, the $n=1$ case has been
extended to $\hat{f}$ supported in $(-2,2)$~\cite{OzS2}~\cite{KS3}.
Chris Hughes has an alternate derivation of~(\ref{eq:Wsp}) appearing
in the notes of these proceedings.

This result confirms the connection between zeros of $L(s,\chi_d)$ and
eigenvalues of unitary symplectic matrices and explains the repulsion
away from the critical point and the 
fluctuations seen
in figure~\ref{fig:zeros real chi}, at least near the real axis,
because,
when $n=1$, the density of zeros is described by the
function $W_{\text{USp}}^{(1)}(x)$ which equals
$$
    1-\frac{\sin(2\pi x)}{2\pi x}.
$$
At height $x$, we therefore also expect, as we average over larger and larger $|d|$,
for the fluctuations to diminish proportional to $1/x$.
However, if we allow $x$ to grow with $d$ then the fluctuations
actually persist due to secondary fluctuating terms that can be large 
if $x$ is allowed to grow with $d$~\cite{Ru3}~\cite{CS}.

The above suggests that the distribution of the lowest zero,
i.e. the one with smallest imaginary part,
in this family of $L$-functions ought to be modeled by the 
distribution of the smallest eigenangle of characteristic polynomials in
$\text{USp}(N)$, with $N \to \infty$. Similary we expect that the distribution, say, of
the second lowest zero ought to fit the distribution of the second smallest
eigenangle.

The probability densities describing the distribution of the smallest and
second smallest eigenangles, normalized by $N/(2\pi)$, for
characteristic polynomials in $\text{USp}(N)$, with $N$ even and tending to 
$\infty$
are given~\cite{KS} respectively by
\begin{equation}
    \notag
    \nu_1(\text{USp})(t) = -\frac{d}{dt} E_{-,0}(t)
\end{equation}
and
\begin{equation}
    \notag
    \nu_2(\text{USp})(t) = -\frac{d}{dt} (E_{-,0}(t)+E_{-,1}(t)),
\end{equation}
where
\begin{align}
    E_{-,0}(t) &= \prod_{j=0}^\infty (1-\lambda_{2j+1}(2t)) \notag \\
    E_{-,1}(t) &= \sum_{k=1}^\infty \lambda_{2k+1}(2t)
                  \prod_{j=0 \atop j\neq k}^\infty (1-\lambda_{2j+1}(2t)). \notag
\end{align}
Here, the $\lambda_j(t)$'s are the eigenvalues of the integral equation
in~(\ref{eq:integral operator}).

\begin{figure}[htp]
    \centerline{
            \psfig{figure=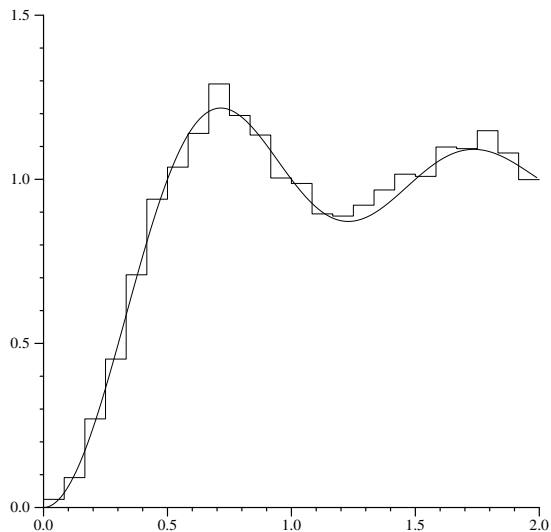,width=3in,angle=0}
    }
    \caption
    {Density of zeros of $L(s,\chi_d)$ for
     7243 prime values of $|d|$ lying in the interval $(10^{12},10^{12}+200000)$.
     Compared against the random matrix theory prediction, $1-\sin(2\pi x)/(2\pi x)$.}
    \label{fig:density zeros usp}
\end{figure}

\begin{figure}[htp]
    \centerline{
            \psfig{figure=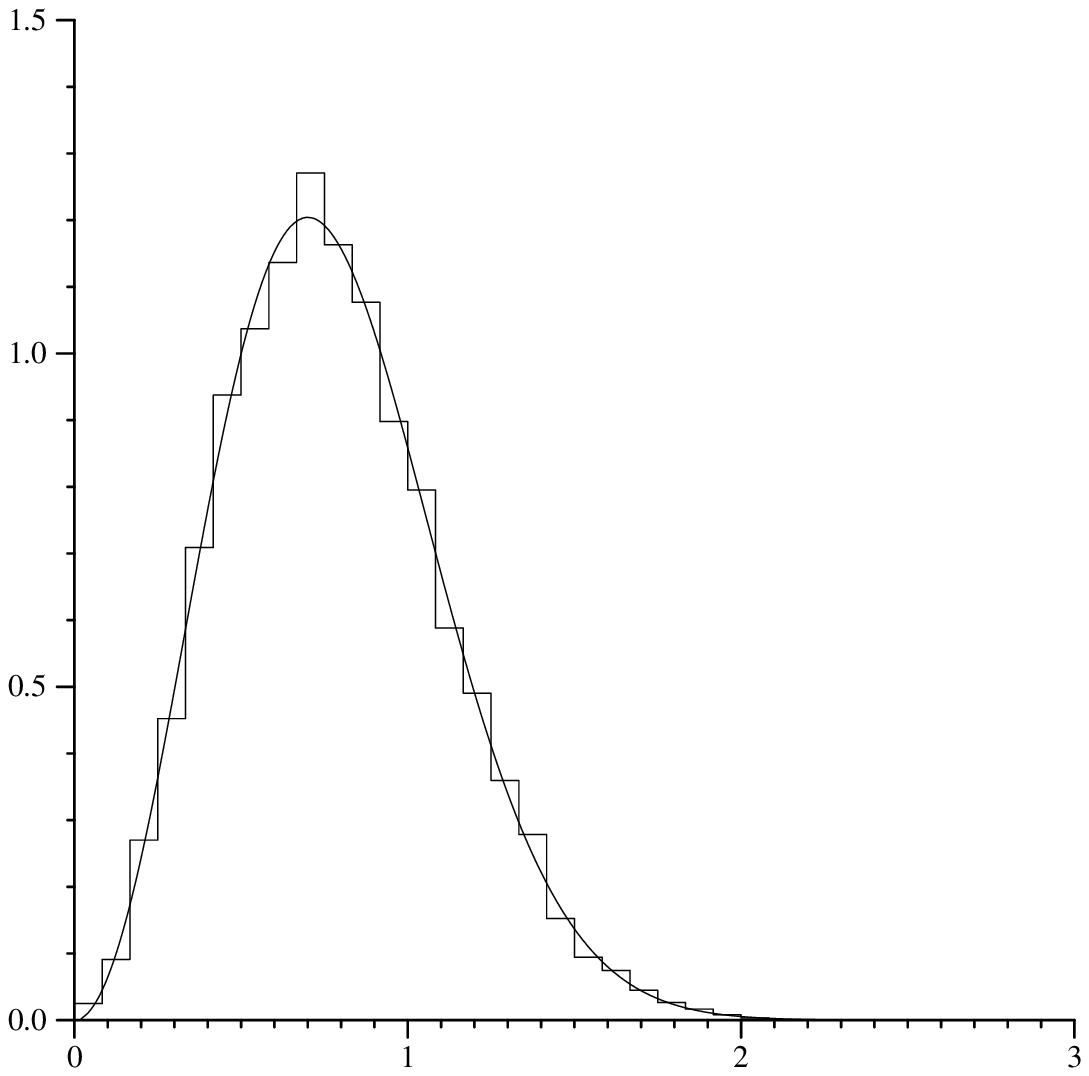,width=3in,angle=0}
            \psfig{figure=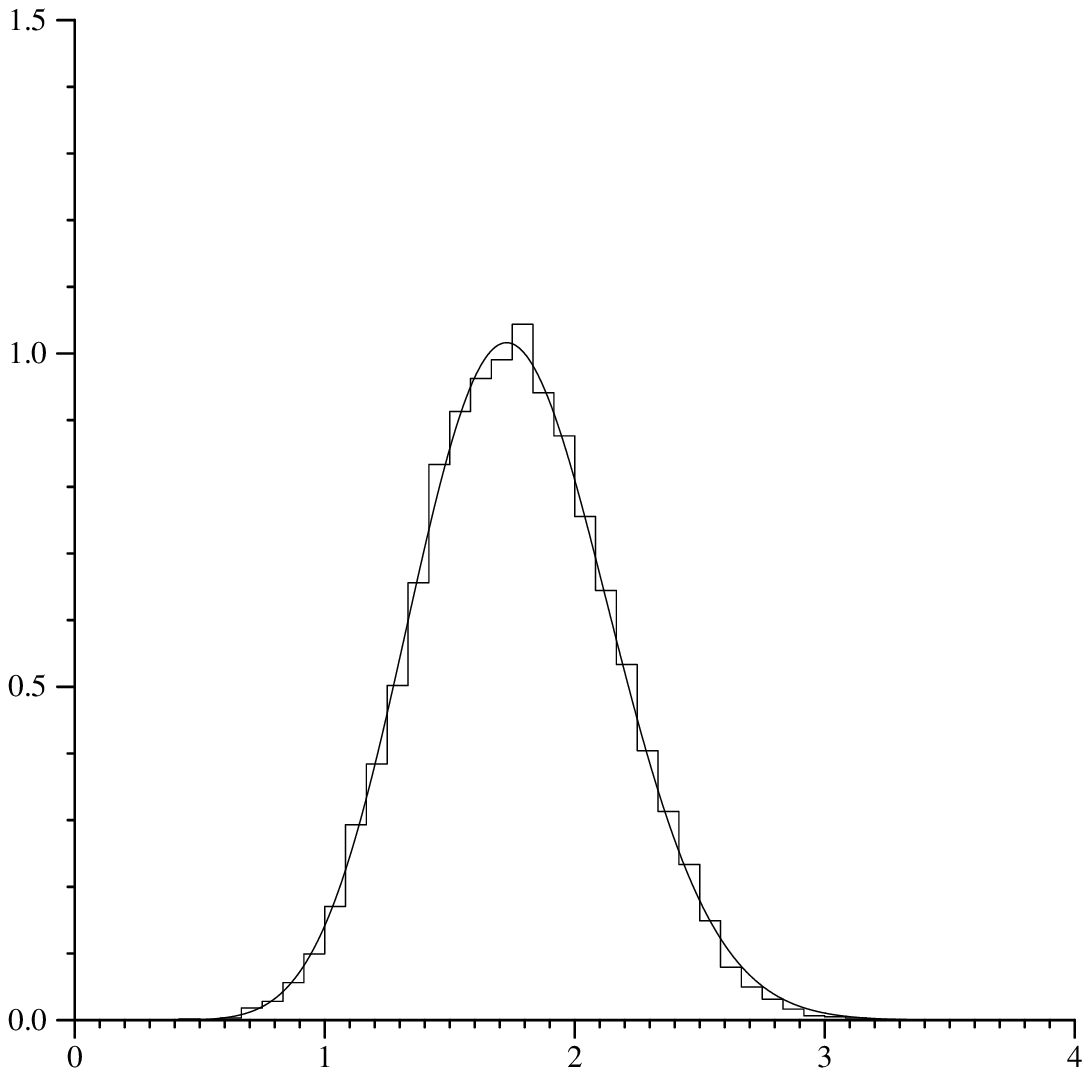,width=3.85in,height=2.83in,angle=0}
    }
    \caption
    {Distribution of the lowest and second lowest zero of $L(s,\chi_d)$ for
     7243 prime values of $|d|$ lying in the interval $(10^{12},10^{12}+200000)$.
     Compared against the random matrix theory predictions.}
    \label{fig:lowest zeros usp}
\end{figure}

This also suggests that the means of the the first and second lowest 
zeros are given by
\begin{align}
    \lim_{X\to\infty}
    \frac{1}{\abs{D(X)}} \sum_{d \in D(X)} \gam{1} l_d
    =
    \int_0^\infty t \nu_1(\text{USp})(t) dt =.78\ldots \notag \\
    \lim_{X\to\infty}
    \frac{1}{\abs{D(X)}} \sum_{d \in D(X)} \gam{2} l_d
    =
    \int_0^\infty t \nu_2(\text{USp})(t) dt = 1.76\ldots \notag
\end{align}
However, the convergence to the predicted means is logarithmically slow
due to secondary terms of size $O(1/\log(X))$. Consequently, when comparing
against the random matrix theory predictions, one gets a better fit
by making sure the lowest zero has the correct mean. This can be achieved
by rescaling the data,
further multiplying, for a set $D$ of fundamental discriminants,
$\gam{1} l_d$ by
\begin{equation}
    \label{eq:renormalize 1}
    .78 
    \left( \frac{1}{\abs{D}} \sum_{d \in D} \gam{1} l_d \right)^{-1}
\end{equation}
and $\gam{2} l_d$ by
\begin{equation}
    \label{eq:renormalize 2}
    1.76 
    \left( \frac{1}{\abs{D}} \sum_{d \in D} \gam{2} l_d \right)^{-1}
\end{equation}
In figures~\ref{fig:density zeros usp} and~\ref{fig:lowest zeros usp}, we use the normalization described above. For our data set,
the denominator in~(\ref{eq:renormalize 1}) equals $.83$, and,
in~(\ref{eq:renormalize 2}) equals $1.84$.

In figure~\ref{fig:density zeros usp} we depict the $1$-level density
of the zeros of $L(s,\chi_d)$ for $7243$ prime $|d|$ lying in the interval
$(10^{12}, 10^{12}+200000)$. These zeros were computed in 1996 as part
of the authors PhD thesis~\cite{Ru}. Here we divide the $x$-axis into 
small bins, count how many normalized zeros of $L(s,\chi_d)$ lie in 
each bin, divide that count by the number of $d$, namely 7243,
and compare that to the graph of $1-\sin(2\pi x)/(2\pi x)$.

In figure~\ref{fig:lowest zeros usp}
we depict the distribution of the lowest and second lowest
normalized zero for the set of zeros just described.
These are compared against $\nu_1$ and $\nu_2$
which were computed using the same program, obtained from Andrew Odlyzko,
that was used in~\cite{O3}. 

\begin{figure}[htp]
    \centerline{
            \psfig{figure=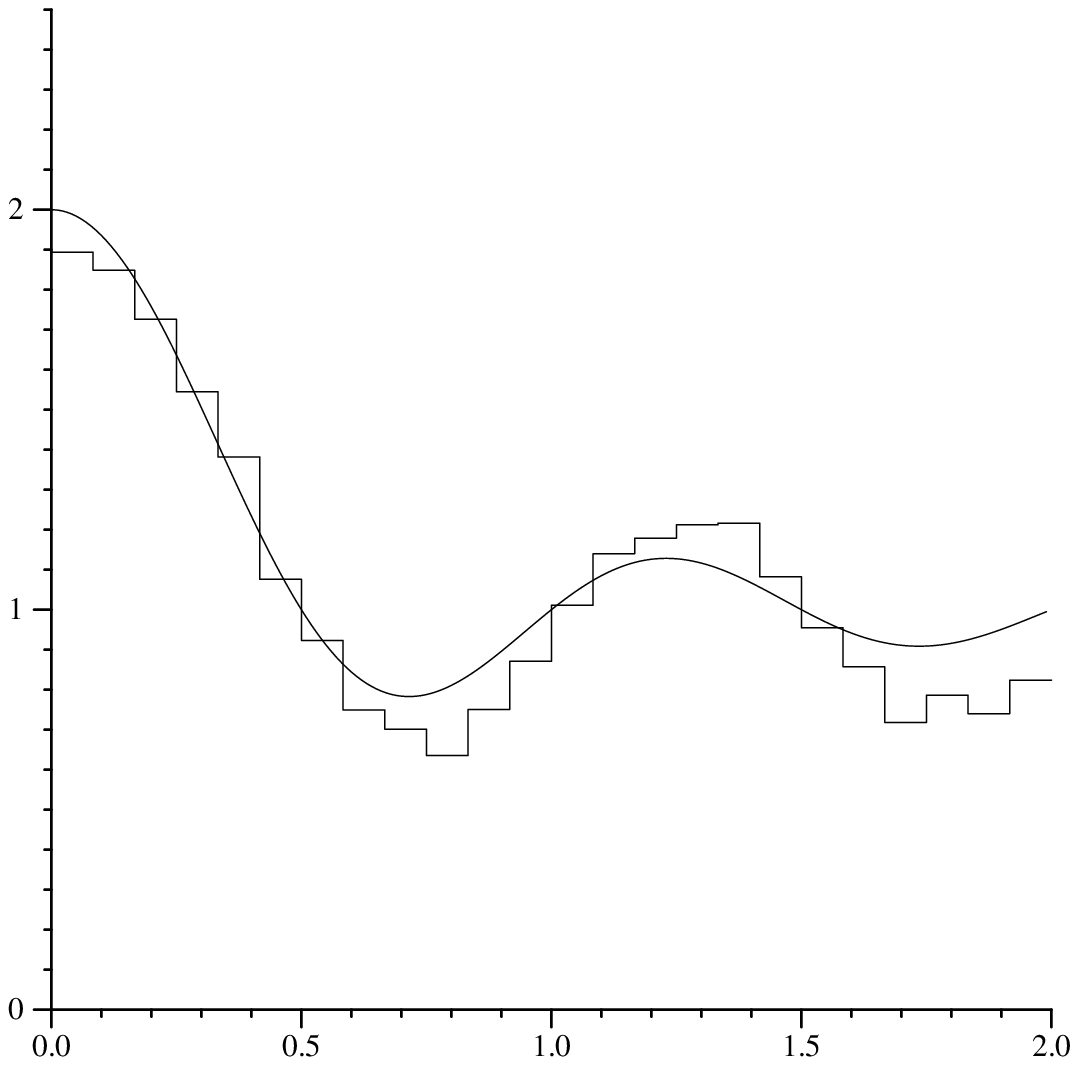,width=2.9in,angle=0}
            \psfig{figure=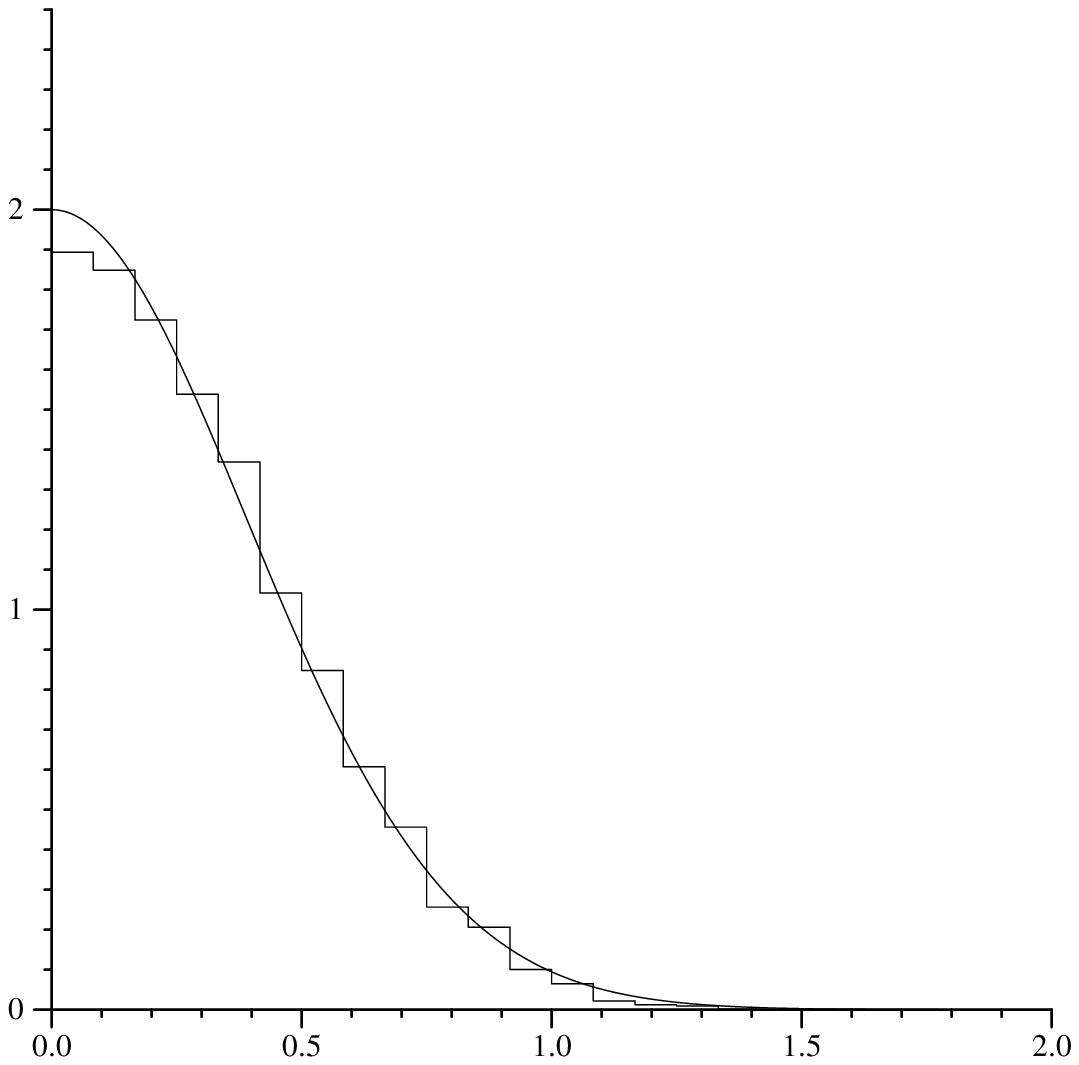,width=3in,angle=0}
    }
    \caption
    {One-level density and distribution of the lowest zero of even
     quadratic twists of the Ramanujan $\tau$ $L$-function, $L_\tau(s,\chi_d)$, for
     11464 prime values of $d>0$ lying in the interval $(350000,650000)$.
     }
    \label{fig:zeros SO}
\end{figure}

In figure~\ref{fig:zeros SO} we depict the $1$-level density and distribution of 
the lowest zeros for quadratic twists of the Ramanujan $\tau$ $L$-function,
$L_\tau(s,\chi_d)$, $d>0$. For this family of $L$-functions,
one can prove~\cite{Ru2} a result similar to~(\ref{eq:Wsp})
but with $W_{\text{USp}}$ replaced with $W_{\text{SO(\text{even})}}$, 
and the support of $\hat{f}$ reduced to $\sum_{i=1}^n \abs{u_i} < 1/2$.
The $1$-level density is therefore given by $1+\sin(2\pi x)/(2\pi x)$ and 
the probability density for the distribution of the smallest eigenangle,
normalized by $2N/(2\pi)$, for matrices in $\text{SO}(2N)$, with $N \to \infty$,
is given~\cite{KS} by 
\begin{equation}
    \notag
    \nu_1(\text{SO(\text{even})})(t) = -\frac{d}{dt} 
    \prod_{j=0}^\infty (1-\lambda_{2j}(2t)),
\end{equation}
whose mean is $.32$.
the figure uses 11464 prime values of $|d|$ lying in
$(350000,650000)$, and the zeros were normalized by $2 l_d$, and then 
rescaled so as to have mean $.32$ rather than $.29$.
The choice of using $2 l_d$ for normalizing the zeros
is the correct one up to leading term, 
but is slightly adhoc and by now a better understanding of a tighter normalization 
up to lower terms has emerged~\cite{CFKRS}~\cite{Ru3}.

\subsection{Value distribution of $L$-functions}

Keating and Snaith initiated the use of random matrix theory to study the
value distribution of $L$-functions with their important paper~\cite{KeS}
where they consider moments of characteristic polynomials
of unitary matrices and conjecture the leading-order
asymptotics for the moments of $\zeta(s)$ on the critical line.
This was followed by a second paper~\cite{KeS2} along
with a paper by Conrey and Farmer~\cite{CF} which
provide conjectures for the leading-order asymptotics
of moments of various families of $L$-functions by examining
analogous questions for characteristic polynomials of the
various classical compact groups.

Keating and Snaith's technically impressive work
also represents a philosophical breakthrough. Until
their paper appeared, one would compare, say, statistics
involving zeros of $\zeta(s)$ to similar statistics for
eigenvalues of $N\times N$ unitary matrices, with $N \to \infty$.
However, their work compares the average value of $\zeta(1/2+it)$
to the average value of $N\times N$ unitary characteristic polynomials
evaluated on the unit circle, with $N \sim \log(t/(2\pi))$. This choice
of $N$ is motivated by comparing local spacings of zeros,
for example~(\ref{eq:normalized spacings2}) 
v.s.~(\ref{eq:normalized eigenangles}).
A slightly different approach to this choice of $N$
proceeds by comparing functional equations of $L$-functions to functional
equations of characteristic polynomials~\cite{CFKRS}. 

At first sight, it seems strange to compare the Riemann zeta function which
has infinitely many zeros to characteristic polynomials of
finite size matrices. However, this suggests
that a given height $t$, the Riemann zeta function can be modeled
locally by just a small number of zeros, as well as by more global information
that incorporates the role played by primes.
Recently, Gonek, Hughes, and Keating have developed such a model~\cite{GHK}.

Below we describe three three specific examples where random matrix theory has
led to important advances in our understanding of the value
distribution of $L$-functions. These concern the families:
\begin{enumerate}
   \item $\zeta(1/2+it)$, where we average over $t$.
   \item $L(1/2,\chi_d)$, where we average over fundamental discriminants $d$.
   \item $L_E(1/2,\chi_d)$, quadratic twists of the $L$-function associated
   to an elliptic curve $E$ over $\mathbb{Q}$, where we average 
   over fundamental discriminants $d$.
\end{enumerate}
These three are examples of unitary, unitary symplectic, and even orthogonal
families respectively~\cite{KS2}~\cite{CFKRS}. Note that in the last example, we normalize
the Dirichlet coefficient of the $L$-function as in~(\ref{eq:dirichlet series E})
so that the functional equation of $L_E(1/2,\chi_d)$
brings $s$ into $1-s$ with the critical point being $s=1/2$.

We first illustrate that these three examples exhibit distinct behaviour
by contrasting their value distributions. The first-order asymptotics for
the moments of $|\zeta(1/2+it)|$ are conjectured by Keating and Snaith~\cite{KeS} to 
be given by
\begin{equation}
   \label{eq:ks moment1}
   \frac{1}{T}
   \int_0^T |\zeta(1/2+it)|^{r} dt \sim
   a_{r/2} \prod_{j=1}^N \frac{\Gamma(j) \Gamma(j+r)}{\Gamma(j+r/2)^2}, \quad \Re{r} > -1,
\end{equation}
with $N \sim \log(T)$
and $a_{r/2}$ defined by~(\ref{eq:zeta_ak}).

For quadratic Dirichlet $L$-functions Keating and Snaith~\cite{KeS2}
conjecture that
\begin{equation}
   \label{eq:ks moment2}
   \frac{1}{\abs{D(X)}} \sum_{d \in D(X)} 
   L(1/2,\chi_d)^r
   \sim
   b_r 
   2^{2Nr}
   \prod_{j=1}^{N} \frac{\Gamma(N+j+1)\Gamma(j+\frac12+r)
   }  {
      \Gamma(N+j+1+r)
   \Gamma(j+\frac12)
   }, \quad \Re{r} > -3/2,
\end{equation}
with $N \sim \log(X)/2$, 
where the sum runs over fundamental discriminants $|d| \leq X$, and,
as suggested by Conrey and Farmer~\cite{CF},
$$
    b_r=\prod_p \frac{(1-\frac 1p)^{\frac{r(r+1)}{2}}}
    {1+\frac 1p}
    \left(\frac{(1-\frac{1}{\sqrt{p}})^{-r}+
    (1+\frac{1}{\sqrt{p}})^{-r}}{2}+\frac{1}{p}\right).
$$

Next, let $q$ be the conductor of the elliptic curve $E$. 
Averaging over fundamental discriminants and restricting to 
discriminants for which $L_E(s,\chi_d)$ has an even functional
equation, the conjecture asserts~\cite{CF}~\cite{KeS2} that
\begin{equation}
   \label{eq:ks moment3}
   \frac{1}{\abs{D(X)}} \sum_{{d \in D(X) \atop (d,q)=1}\atop \text{even funct eqn}} 
   L_E(1/2,\chi_d)^r
   \sim
   c_r
    2^{2Nr} \prod_{j=1}^{N} \frac{\Gamma(N+j-1)
    \Gamma(j-\frac12+r)}
    {
    \Gamma(N+j-1+r)
    \Gamma(j-\frac12)
    }, \quad \Re{r} > -1/2,
\end{equation}
with $N \sim \log(X)$ and
\begin{equation}
     \notag
    c_r =
    \prod_p 
    \left(1-\frac{1}{p}\right)^{k(k-1)/2}
    R_{r,p} 
\end{equation}
where, for $p \nmid q$,
\begin{equation}
     \notag
    R_{k,p}  =     
    \left(1+\frac 1 p\right)^{-1}\left(\frac 1 p +\frac{1}{2}
    \left(
    \left(1-\frac{a_p}{p}+\frac{1}{p}\right)^{-k}+\left(1+\frac{a_p}{p}+\frac{1}{p}
    \right)^{-k}
    \right)\right) .
\end{equation}
In the above equation, $a_p$ stands for the $p$th coefficient of the Dirichlet series of $L_E$.

In the case of the Riemann zeta function we take absolute values,
$|\zeta(1/2+it)|$, otherwise the moments would be zero.
In the other two cases, the $L$-values
are conjectured to be non-negative real numbers, hence we
directly take their moments.

We should also observe that, while statistics such as the pair correlation
or density of zeros
involving zeros of $L$-functions
have arithmetic information appearing in the secondary terms,
moments already reveal such behaviour at the level of the main term.
This reflects the global nature of the moment statistic as compared to
the local nature of statistics of zeros that have been discussed.

\begin{figure}[htp]
    \centerline{
            \psfig{figure=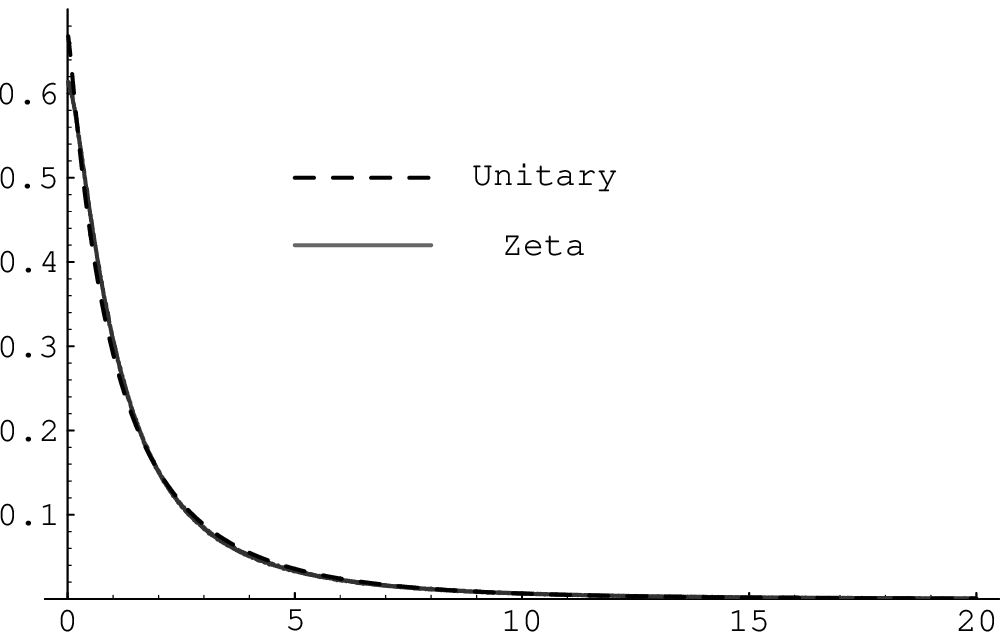,width=2.7in,angle=0}
    }
    \centerline{
            \psfig{figure=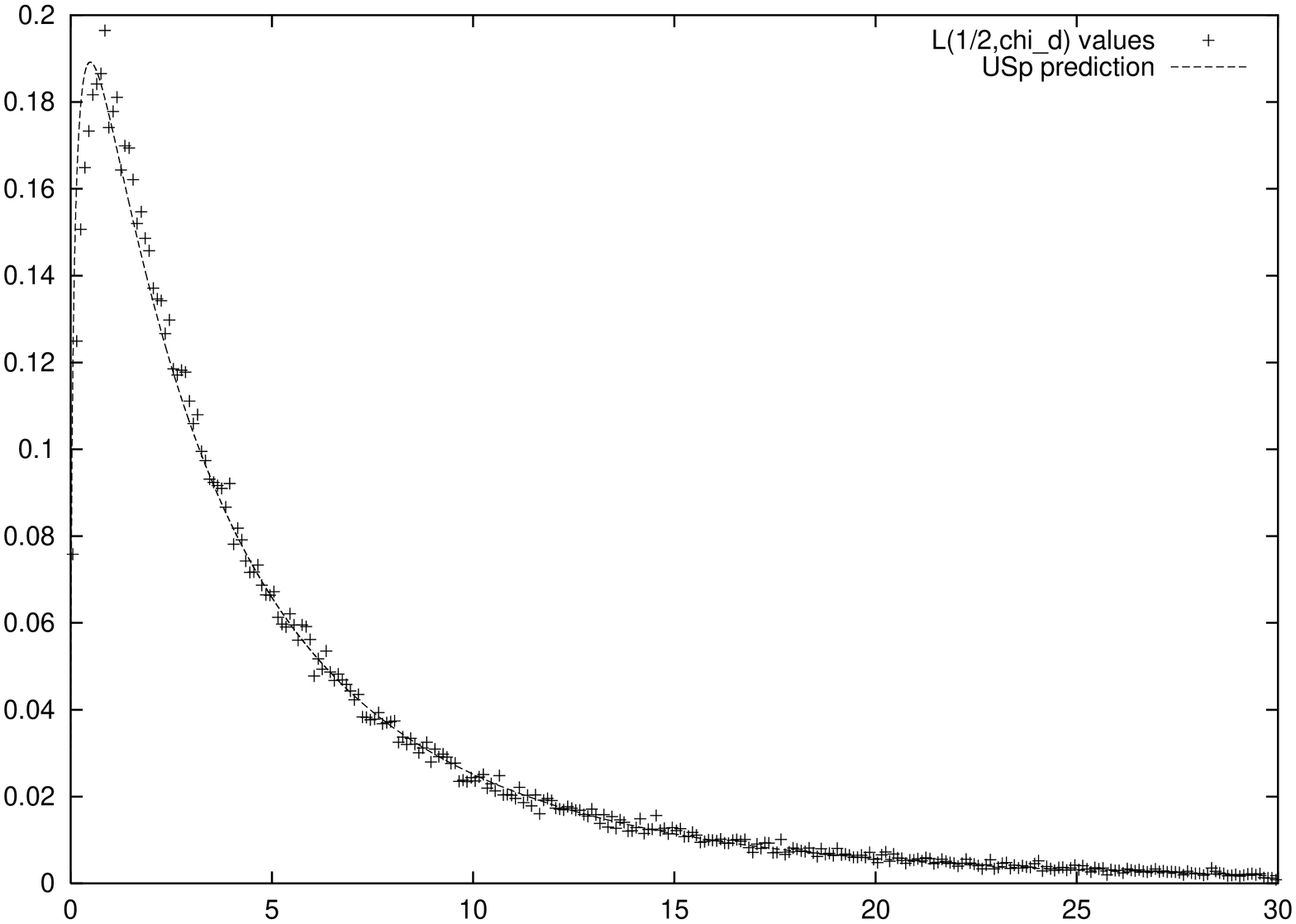,width=2.7in,angle=-90}
            \psfig{figure=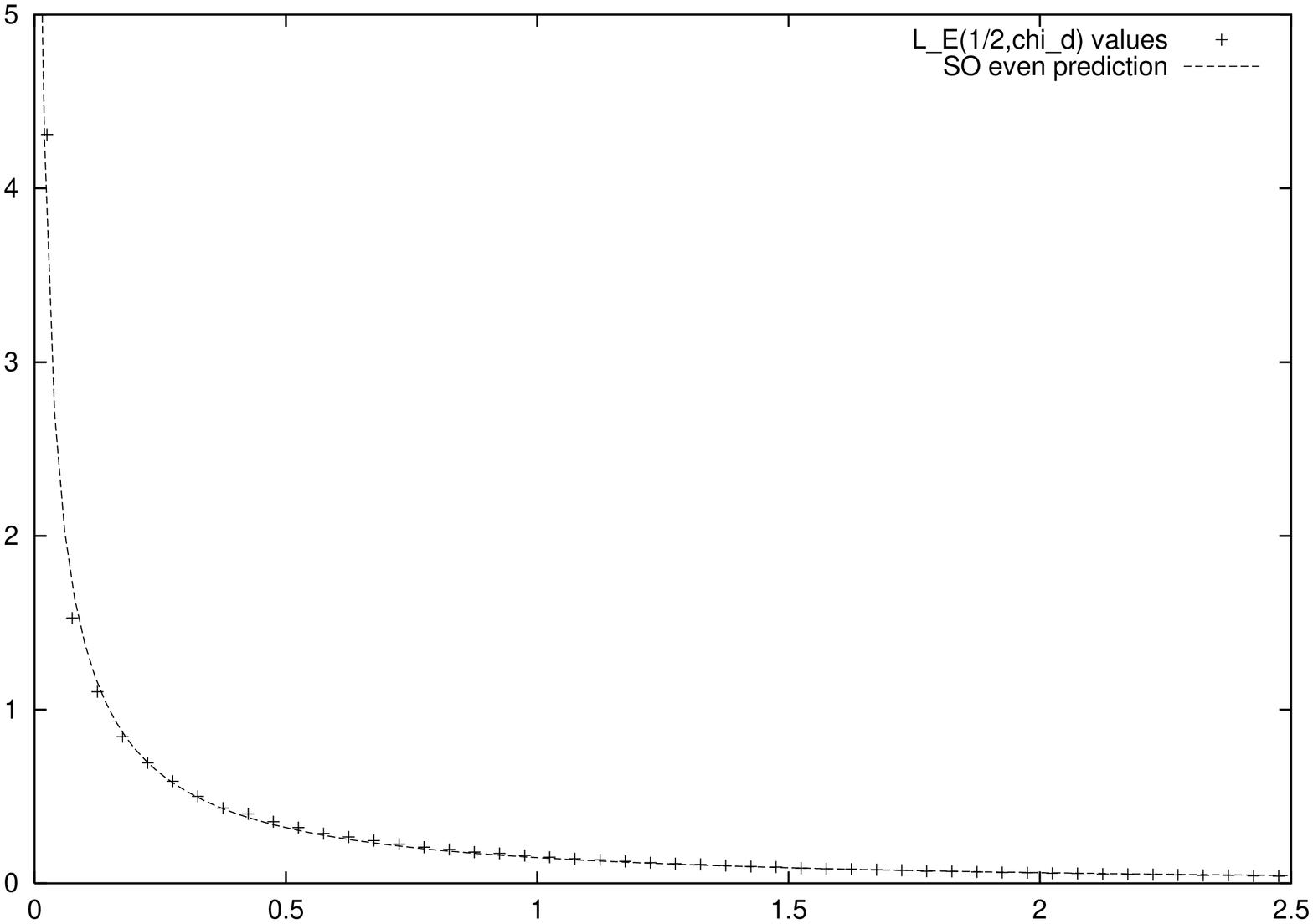,width=2.7in,angle=-90}
    }
    \caption{
        Value distribution of $L$-functions compared to the random matrix theory counterparts.
        The first picture, depicts the value distribution of $|\zeta(1/2+it)|$, with $t$ near $10^6$,
        the second of $L(1/2,\chi_d)$ with $800000< |d| < 10^6$, and the third
        of $L_{E_{11}}(1/2,\chi_d)$, with $-85000000<d<0$, $d=2,6,7,8,10 \mod 11$.}
    \label{fig:value dist three families}
\end{figure}

Using the above conjectured asymptotics we can naively plot value distributions.
Figure~\ref{fig:value dist three families} compares numerical value distributions
for data in these three examples against the counterpart densities
from random matrix theory.
Notice that these three graphs behave distinctly near the origin. 
The solid curves are computed by taking inverse Mellin transforms, as 
in~(\ref{eq:inverse mellin p(t)}), of the right hand sides
of equations~(\ref{eq:ks moment1}), ~(\ref{eq:ks moment2}), and ~(\ref{eq:ks moment3}), 
but without the arithmetic
factors $a_k$, $b_k$, $c_k$. Shifting the inverse Mellin transform line integral
to the left, the location of the first pole in each integrand dictates the
behaviour of the corresponding density functions near the origin.
The locations of these three poles are at $r=-1,-3/2,$ and $-1/2$
respectively.
Taking $t$ to be the horizontal axis,
near the origin the first density is 
proportional to a constant, the second to $t^{1/2}$, and the third to $t^{-1/2}$.
In forming these graphs one takes $N$ as described above, so that 
the proportionality constants do depend on $N$.
As $N$ grows, these graphs tend to get flatter.

The first graph is reproduced from~\cite{KeS}.
In the second and third 
graphs displayed, a slight cheat was used to get a better fit. The histograms
were rescaled linearly along both axis until the histogram matched
up nicely with the solid curves. 
We must ignore the arithmetic factors when taking inverse
Mellin transforms since these factors are known~\cite{CGo} to be functions of order two 
and cause the inverse Mellin transforms to diverge. 
To properly plot the
correct value distributions we would need to use more than just the
leading-order asymptotics.
Presently, our knowledge of the moments of various families of $L$-functions
extends beyond the first-order
asymptotics, but only for positive integer values of $r$ 
(even integer in the case of $|\zeta(1/2+it)|^r$),
however, 
one would need to apply full asymptotics for complex values of $r$.
The paper by Conrey, Farmer, Keating, the author, and Snaith~\cite{CFKRS}
conjectures the full asymptotics, for example, of the
three moment problems above, but for integer $r$,
with corresponding
theorems in random matrix theory given in~\cite{CFKRS2}.
The paper of Conrey, Farmer and Zirnbauer goes even beyond this stating conjectures for 
the full asymptotics of moments
of ratios of $L$-functions, and, using 
methods from supersymmetry, proving corresponding theorems in random matrix theory~\cite{CFZ}.
Another paper, by Conrey, Forrester, and Snaith,
uses orthognal polynomials to obtain alternative proofs 
of the random matrix theory theorems for ratios~\cite{CFS}.

\subsubsection{Moments of $|\zeta(1/2+it)|$}

Next we describe the full moment conjecture from~\cite{CFKRS}
for $|\zeta(1/2+it)|$. In that paper, the conjecture is derived heuristically
by looking at products of zetas shifted slightly away from the critical line and
then setting the shifts equal to zero. 

The formula is written in terms of contour integrals and involve
the Vandermonde:
\begin{equation}
    \notag
    \Delta(z_1,\dots,z_{m})=\prod_{1\le i < j\le m}(z_j-z_i).
\end{equation}

Suppose $g(t)=f(t/T)$ with
$f:{\mathbb{R}}^+ \to {\mathbb{R}}$ non-negative, bounded, and
integrable. 
The conjecture of~\cite{CFKRS} states that, as $T \to \infty$,
\begin{equation}
     \notag
    \int_0^\infty |\zeta(1/2 +it)|^{2k}g(t)\,dt
    \sim \int_0^\infty
    P_k\left(\log (t/(2 \pi))\right)
    g(t) \,dt ,
\end{equation}
where $P_k$ is the polynomial of degree $k^2$ given by the
$2k$-fold residue
\begin{equation}
     \notag
    P_k(x)= \frac{(-1)^k}{k!^2}\frac{1}{(2\pi i)^{2k}} \oint\cdots
    \oint \frac{G(z_1, \dots,z_{2k})\Delta^2(z_1,\dots,z_{2k})}
    {\displaystyle \prod_{j=1}^{2k} z_j^{2k}}
    e^{\tfrac x2 \sum_{j=1}^{k}z_j-z_{k+j}}\,dz_1\dots dz_{2k} ,
\end{equation}
where one integrates over small circles about $z_i=0$, with
\begin{equation}
     \notag
    G(z_1,\dots,z_{2k})= A_k(z_1,\dots,z_{2k})
    \prod_{i=1}^k\prod_{j=1}^k\zeta(1+z_i-z_{k+j}) ,
\end{equation}
and $A_k$ is the Euler product
\begin{align}
    A_k(z) &=\prod_p \prod_{i=1}^k\prod_{j=1}^k
    \left(1-\frac{1}{p^{1+z_i-z_{k+j}}}\right) \int_0^1 \prod_{j=1}^k
    \left(1-\frac{e^{2\pi i\theta}}{p^{\frac12 +z_j}}\right)^{-1}
    \left(1-\frac{e^{-2\pi i\theta}}{p^{\frac12 -z_{k+j}}}\right)^{-1}\,d\theta \notag \\
    &=
    \prod_{p}
    \sum_{m=1}^k
        \prod_{i\neq m}
           \frac{\displaystyle \prod_{j=1}^k \left(1-\frac{1}{p^{1+z_j-z_{k+i}}}\right)}
                             {1-p^{z_{k+i}-z_{k+m}}}. \notag
\end{align}

When $k=1$ or $2$, this conjecture agrees with theorems for 
the full asymptotics 
as worked out by Ingham~\cite{I} and Heath-Brown respectively~\cite{H}.
In the first case $A_1(z)=1$ and in the second case 
$A_2(z)=\zeta(2+z_1+z_2-z_3-z_4)^{-1}$,
and one can
write down the coefficients of the polynomials $P_k(x)$ in terms of
known constants. 
When $k=3$ the product over primes becomes rather complicated.
However, one can numerically evaluate~\cite{CFKRS3} the coefficients of $P_3(x)$
and the polynomial is given by:
\begin{eqnarray}
    P_3(x) =  0.000005708527034652788398376841445252313\,x^9&& \notag \\
    + 0.00040502133088411440331215332025984\,x^8&& \notag \\
    +  0.011072455215246998350410400826667\,x^7&& \notag \\
    + 0.14840073080150272680851401518774\,x^6&& \notag \\
    +  1.0459251779054883439385323798059\,x^5&& \notag \\
    + 3.984385094823534724747964073429\,x^4&& \notag \\
    +  8.60731914578120675614834763629\,x^3&& \notag \\
    + 10.274330830703446134183009522\,x^2&& \notag \\
    +  6.59391302064975810465713392\,x&& \notag \\
    + 0.9165155076378930590178543.&&
     \notag
\end{eqnarray}
In the $k=3$ case the moments of $|\zeta(1/2+it)|$ have not
been proven, and it makes sense to test the moment conjecture numerically.
Table~\ref{tab:sixthmoment}, reproduced from~\cite{CFKRS}, depicts
\begin{equation}
\int_C^D |\zeta(1/2 +it)|^6 dt \label{eqn:zeta6smooth}
\end{equation}
as compared to
\begin{equation}
\int_C^D P_3(\log(t/2\pi)) dt,
 \label{eqn:p3smooth}
\end{equation}
along with their ratio, for various blocks $[C,D]$ of length
50000, as well as a larger block of length 2,350,000. 

\begin{table}[h!tb]
\centerline{
\begin{tabular}{|c|c|c|c|}
\hline
$[C,D]$ & conjecture (\ref{eqn:p3smooth}) & reality (\ref{eqn:zeta6smooth})& ratio \cr \hline
[0,50000] & 7236872972.7 & 7231005642.3    & .999189\cr
[50000,100000] & 15696470555.3 & 15723919113.6 & 1.001749\cr
[100000,150000]& 21568672884.1&21536840937.9   &     .998524\cr
[150000,200000]& 26381397608.2&26246250354.1   &     .994877\cr
[200000,250000]& 30556177136.5&30692229217.8   &    1.004453\cr
[250000,300000]& 34290291841.0&34414329738.9   &    1.003617\cr
[300000,350000]& 37695829854.3&37683495193.0   &     .999673\cr
[350000,400000]& 40843941365.7&40566252008.5   &     .993201\cr
[400000,450000]& 43783216365.2&43907511751.1   &    1.002839\cr
[450000,500000]& 46548617846.7&46531247056.9   &     .999627\cr
[500000,550000]& 49166313161.9&49136264678.2  &      .999389\cr
[550000,600000]& 51656498739.2&51744796875.0  &     1.001709\cr
[600000,650000]& 54035153255.1&53962410634.2  &      .998654\cr
[650000,700000]& 56315178564.8&56541799179.3  &     1.004024\cr
[700000,750000]& 58507171421.6&58365383245.2  &      .997577\cr
[750000,800000]& 60619962488.2&60870809317.1  &     1.004138\cr
[800000,850000]& 62661003164.6&62765220708.6  &     1.001663\cr
[850000,900000]& 64636649728.0&64227164326.1  &      .993665\cr
[900000,950000]& 66552376294.2&65994874052.2  &      .991623\cr
[950000,1000000]& 68412937271.4&68961125079.8  &     1.008013\cr
[1000000,1050000]& 70222493232.7&70233393177.0  &
1.000155\cr [1050000,1100000]& 71984709805.4&72919426905.7  &
1.012985\cr [1100000,1150000]& 73702836332.4&72567024812.4  &
.984589\cr [1150000,1200000]& 75379769148.4&76267763314.7  &
1.011780\cr [1200000,1250000]& 77018102997.5&76750297112.6  &
.996523\cr [1250000,1300000]& 78620173202.6&78315210623.9  &
.996121\cr [1300000,1350000]& 80188090542.5&80320710380.9  &
1.001654\cr [1350000,1400000]& 81723770322.2&80767881132.6  &
.988303\cr [1400000,1450000]& 83228956776.3&83782957374.3  &
1.006656\cr [0,2350000]& 3317437762612.4&3317496016044.9  &
1.000017 \cr \hline
\end{tabular}
}
\caption{ Sixth moment of $\zeta$
versus the conjecture. The `reality' column, i.e. integrals
involving $\zeta$, were computed using Mathematica.
}\label{tab:sixthmoment}\end{table}

\subsubsection{Moments of $L(1/2,\chi_d)$}

Another conjecture listed in~\cite{CFKRS} concerns the full asymptotics for
the moments of $L(1/2,\chi_d)$. We quote the conjecture here:

Suppose $g(t)=f(t/T)$ with
$f:{\mathbb{R}}^+ \to {\mathbb{R}}$ non-negative, bounded, and
integrable. 
Let $X_d(s)=|d|^{\frac12-s}X(s,a)$ where
$a=0$ if $d>0$   and  $a=1$ if $d<0$, and
\begin{equation}
    \notag
    X(s,a) = \pi^{s-\frac12}
    \Gamma\left(\frac{1+a-s}{2}\right)/\Gamma\left(\frac{s+a}{2}\right) .
\end{equation}
That is, $X_d(s)$ is the factor in the functional equation
$L(s,\chi_d)=X_d(s)L(1-s,\chi_d)$.
Summing over negative fundamental discriminants $d$ we have, as $T \to \infty$,
\begin{equation}
    \notag
    \sum_{\!d<0}L(1/2 ,\chi_d)^kg(|d|) \sim \sum_{\!d<0}\,Q_k(\log
    {|d|}) g(|d|)
\end{equation}
where $Q_k$ is the polynomial of degree $k(k+1)/2$ given by
the $k$-fold residue
\begin{equation}
    \notag
    Q_k(x)= \frac{(-1)^{k(k-1)/2}2^k}{k!}
    \frac{1}{(2\pi i)^{k}}
    \oint \cdots \oint
    \frac{G(z_1, \dots,z_{k})\Delta(z_1^2,\dots,z_{k}^2)^2}
    {\displaystyle \prod_{j=1}^{k} z_j^{2k-1}}
    e^{\tfrac x2 \sum_{j=1}^{k}z_j}\,dz_1\dots dz_{k} ,
\end{equation}
where
\begin{equation}
    \notag
    G(z_1,\dots,z_k)=B_k(z_1,\dots,z_k)
    \prod_{j=1}^k X(1/2+z_j,1)^{-\frac12}
    \prod_{1\le i\le j\le k}\zeta(1+z_i+z_j),
\end{equation}
and $B_k$ is the Euler product, absolutely convergent for
$|\Re z_j|<\frac12 $, defined by
\begin{align}
    B_k(z_1,\dots,z_k) & =  \prod_p \prod_{1\le i \le j \le k}
    \left(1-\frac{1}{p^{1+z_i+z_j}}\right) \notag \\ &\times \left(\frac
    12 \left(\prod_{j=1}^k\left( 1-\frac{1}{p^{\frac 12+z_j}}\right)^{-1} +
    \prod_{j=1}^k\left(1+\frac{1}{p^{\frac12+z_j}}\right)^{-1}
    \right)+\frac 1p \right) \left( 1+ \frac{1}{p}\right)^{-1}. \notag
\end{align}
We can also sum over $d>0$ but then need to replace $X(1/2+z_j,1)$
with $X(1/2+z_j,0)$.

This conjecture agrees with theorems in the case of $k=1,2,3$~\cite{J}~\cite{S}
(only the leading term has been checked in the case of $k=3$, but in principle
the lower terms could be verified).

Figure \ref{fig:conj vs reality, even chi}, reproduced from~\cite{CFKRS},
depicts, for $k=1,\ldots,8$ and $X=10000,20000,\ldots,10^7$,
\begin{equation}
    \notag
    \sum_{0<d\leq X} L(1/2 ,\chi_d)^k
\end{equation}
divided by
\begin{equation}
   \notag
   \sum_{0<d\leq X}\,Q_k(\log{d}).
\end{equation}

\begin{figure}[htp]
    \centerline{
            \psfig{figure=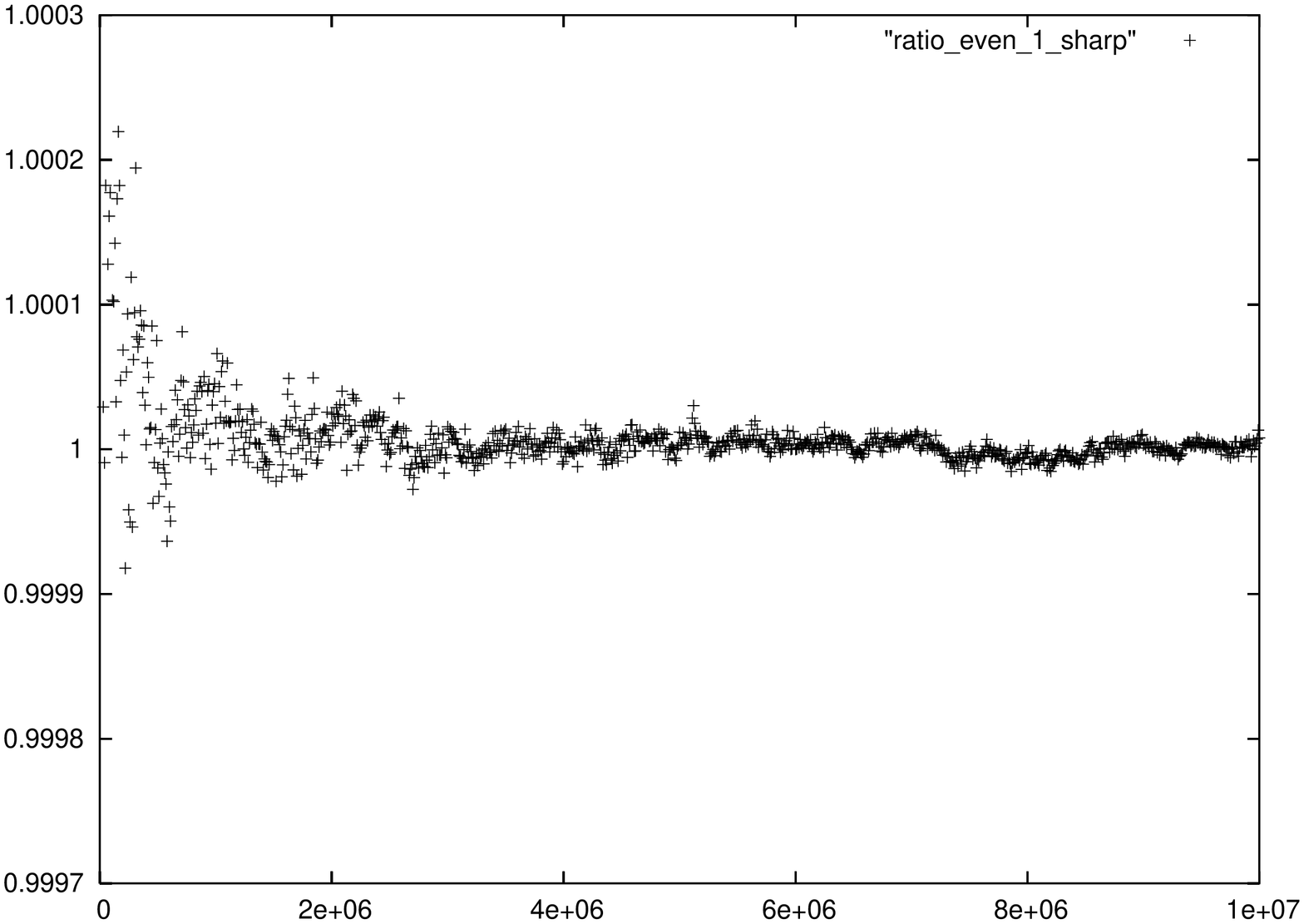,width=3in,height=1.9in,angle=-90}
            \psfig{figure=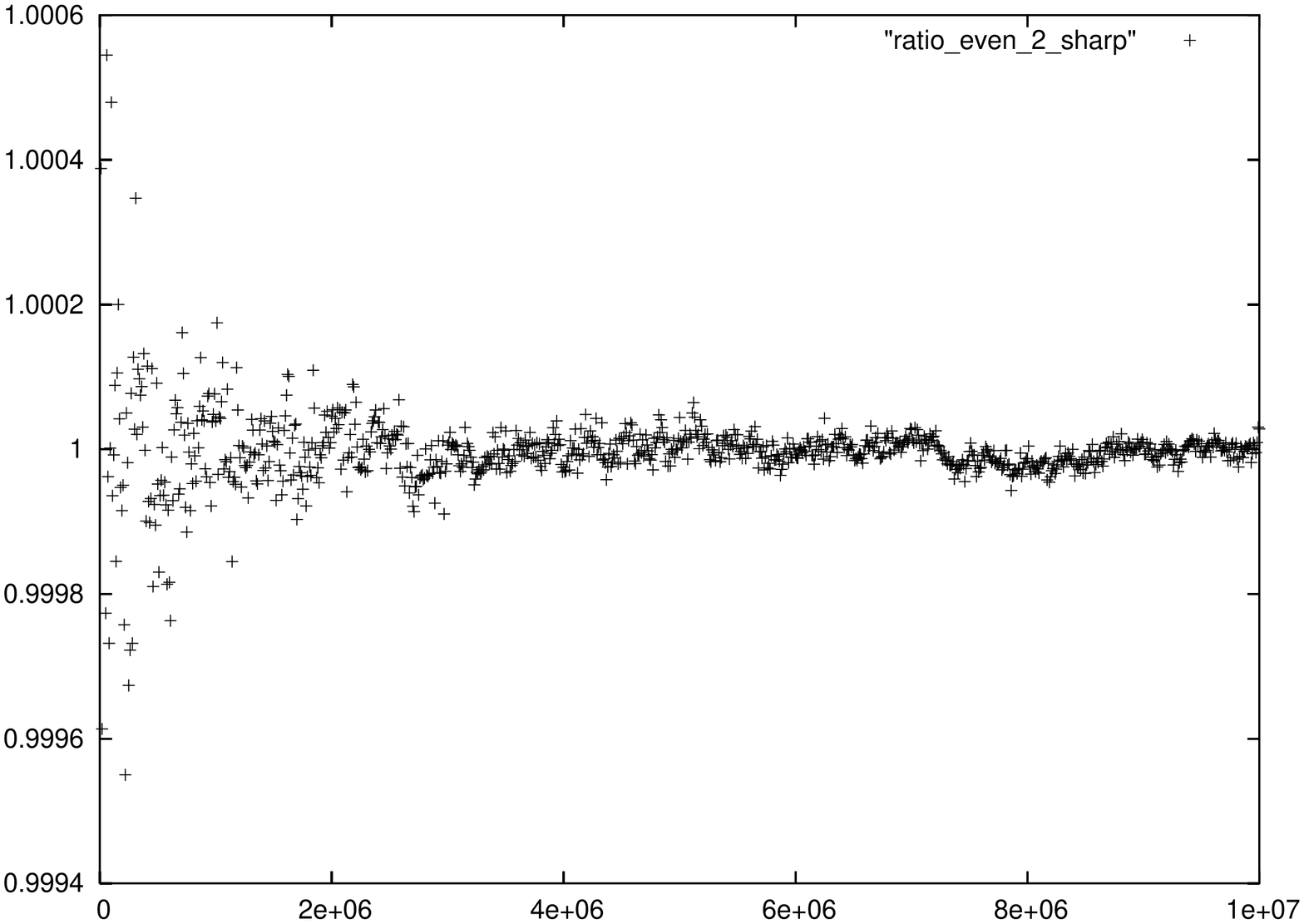,width=3in,height=1.9in,angle=-90}
    }
    \centerline{
            \psfig{figure=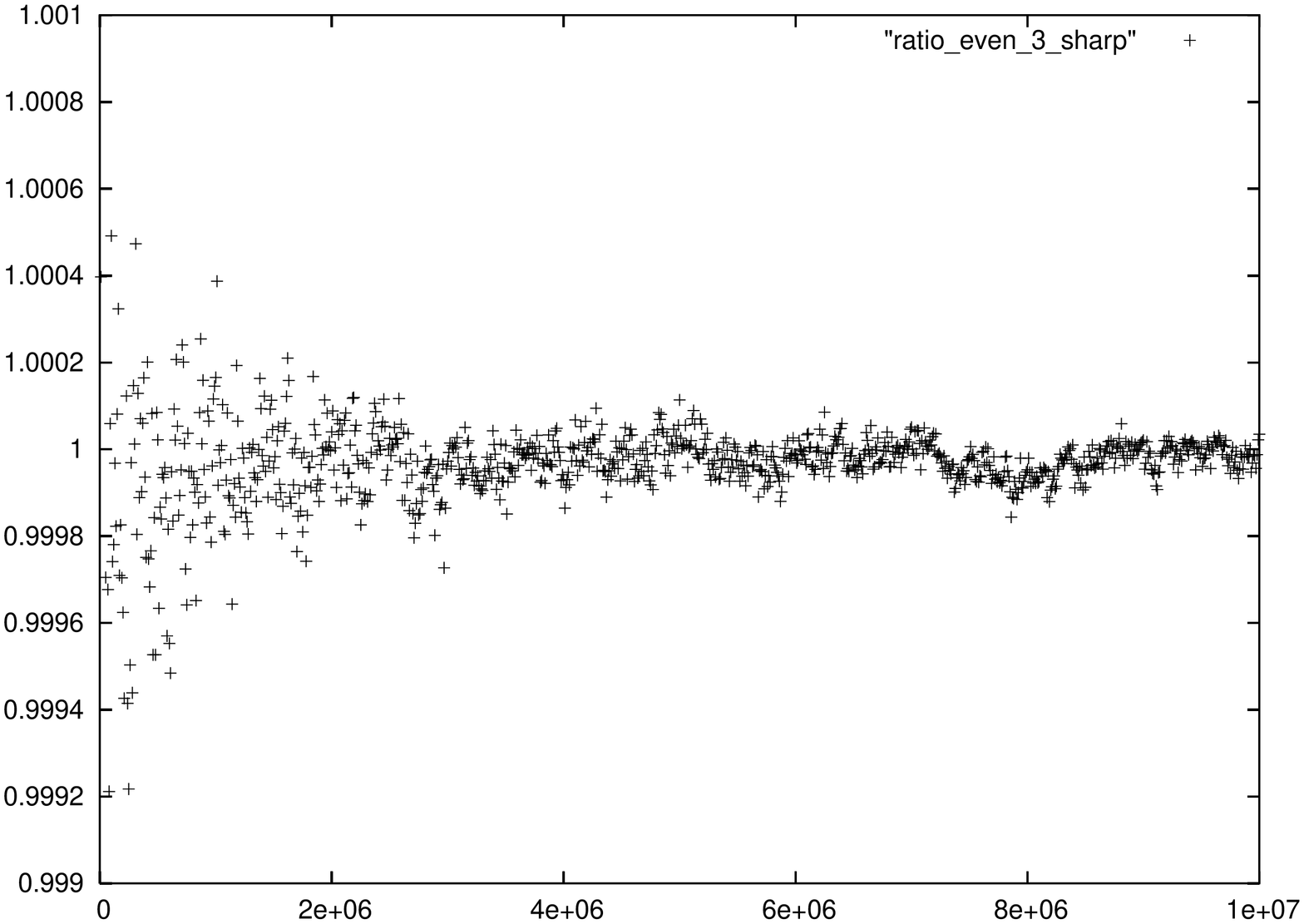,width=3in,height=1.9in,angle=-90}
            \psfig{figure=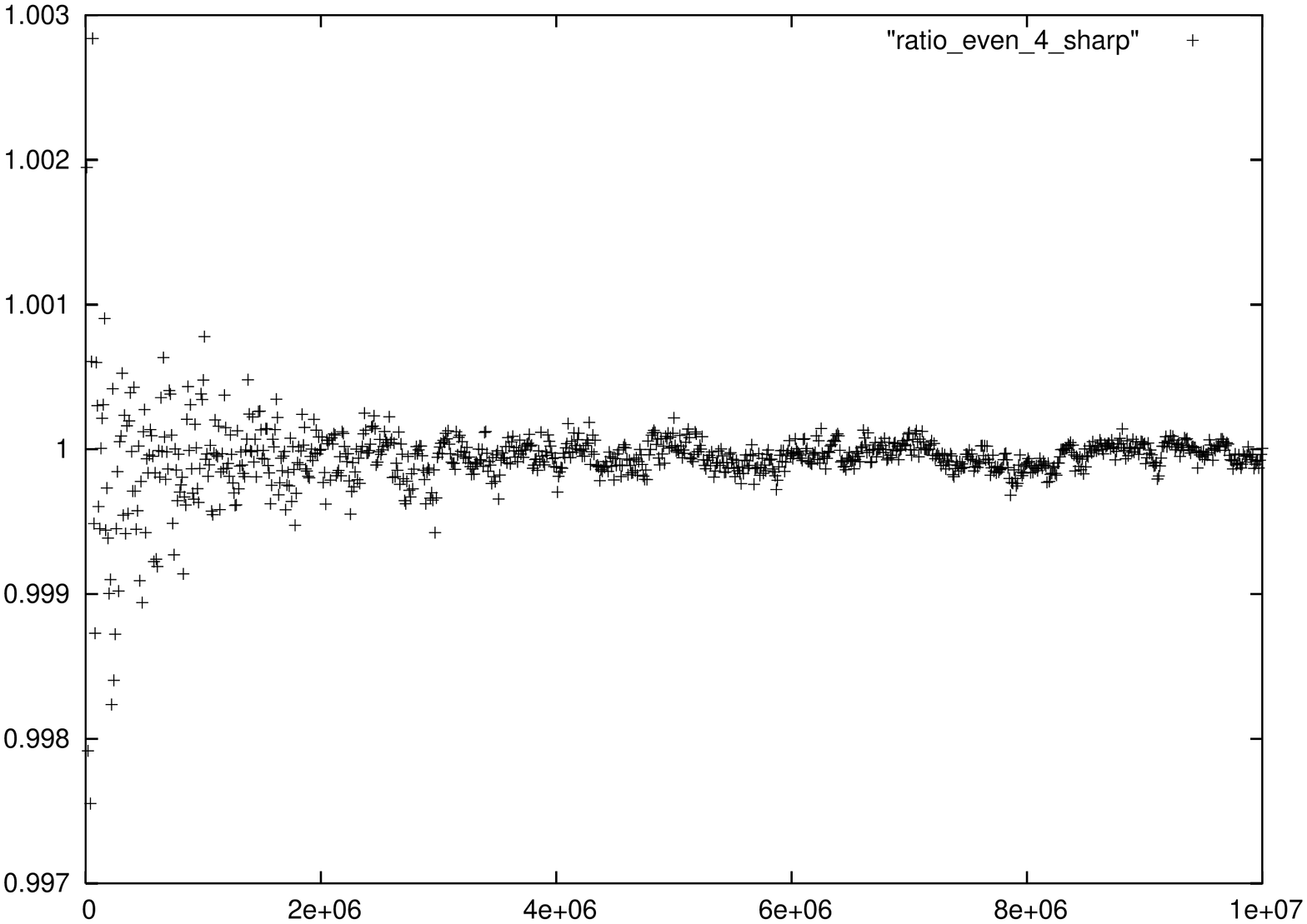,width=3in,height=1.9in,angle=-90}
    }
    \centerline{
            \psfig{figure=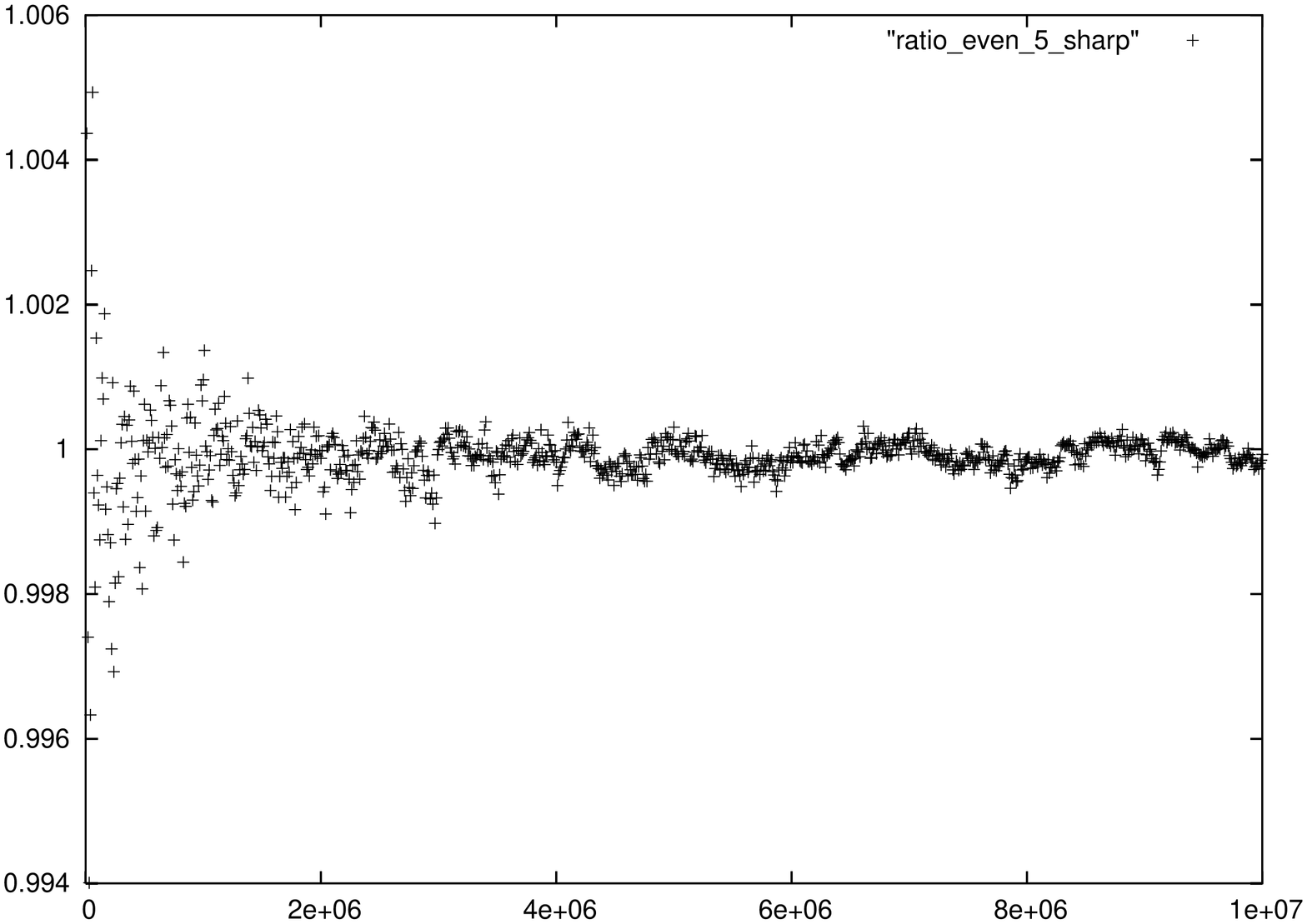,width=3in,height=1.9in,angle=-90}
            \psfig{figure=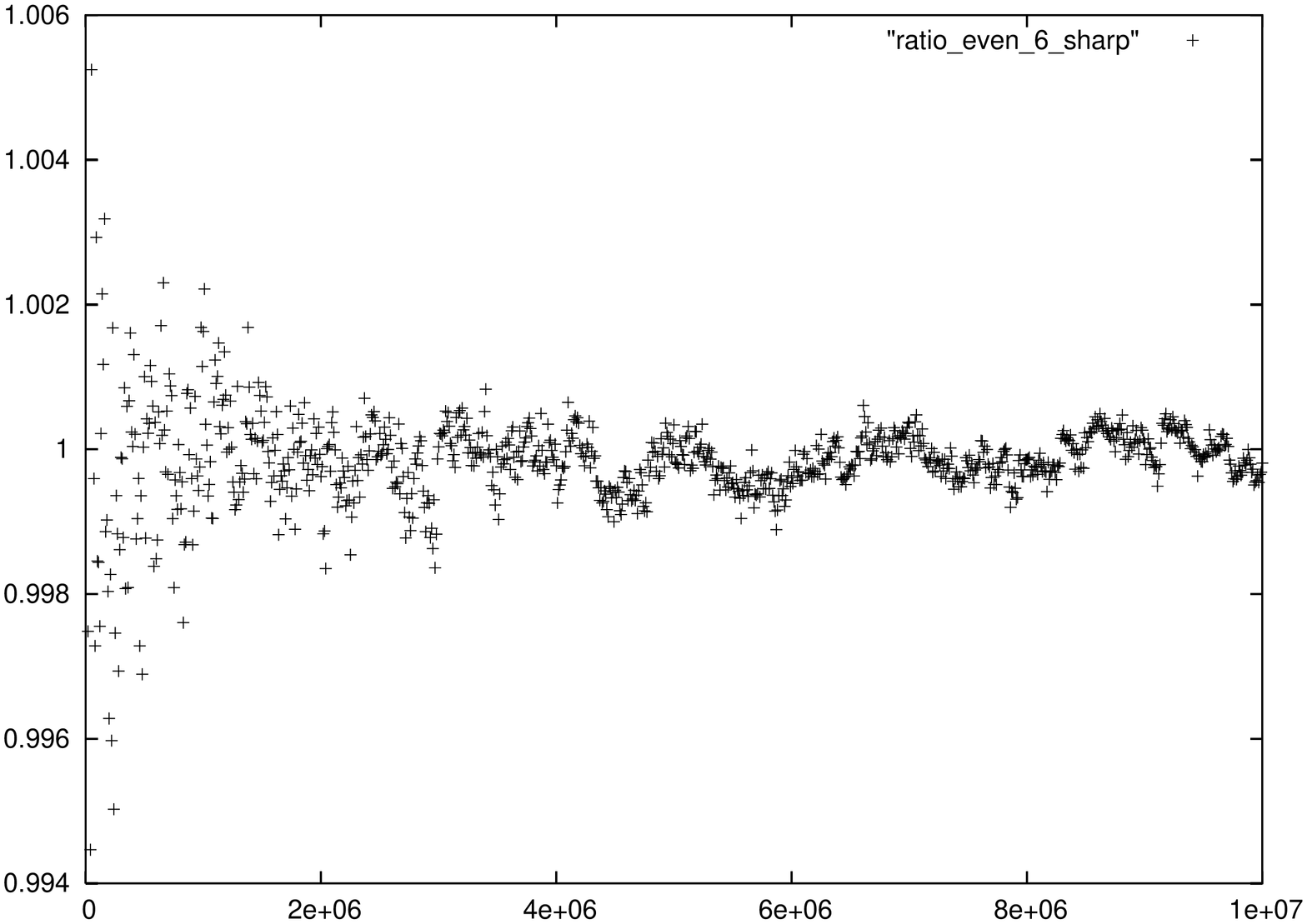,width=3in,height=1.9in,angle=-90}
    }
    \centerline{
            \psfig{figure=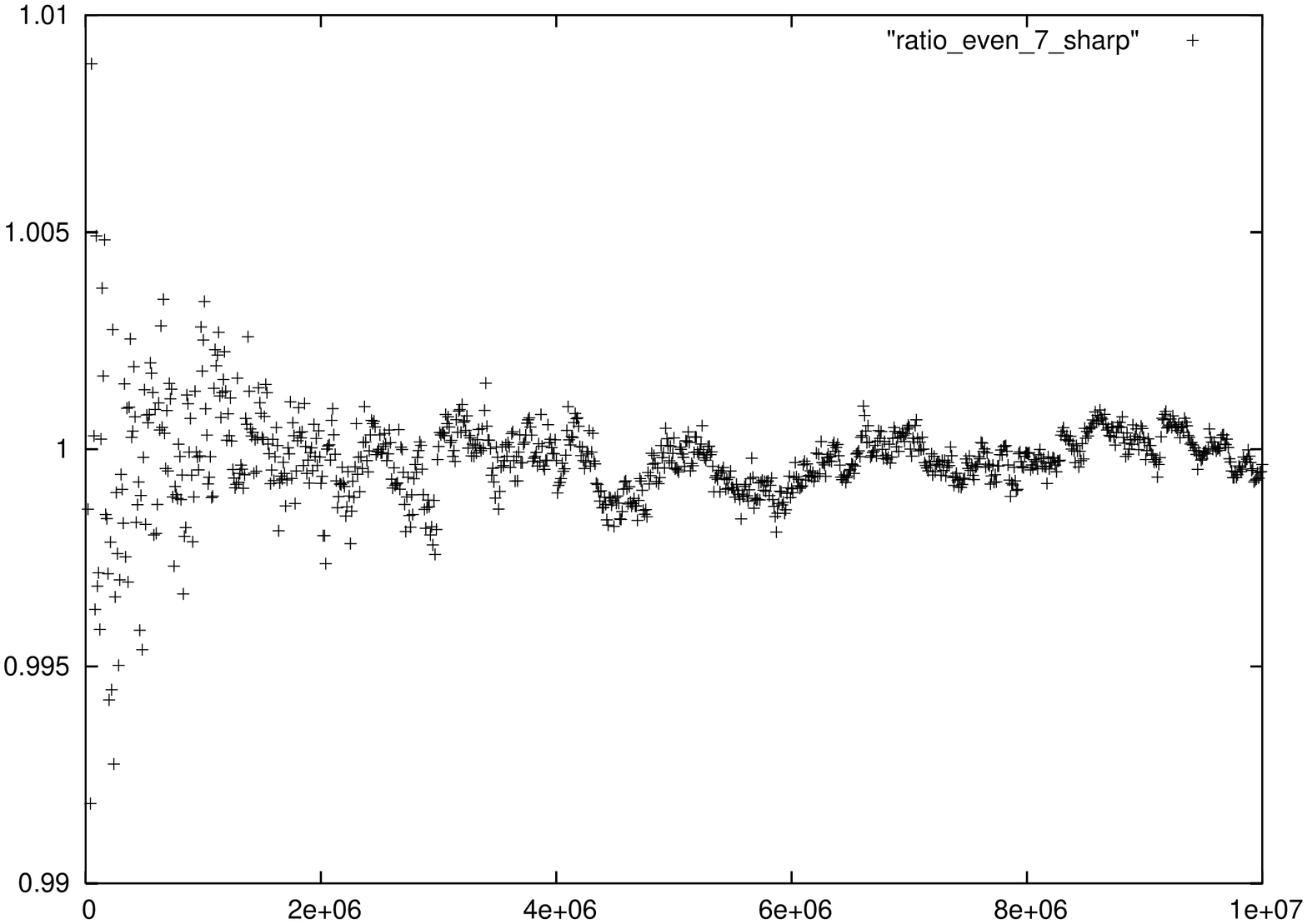,width=3in,height=1.9in,angle=-90}
            \psfig{figure=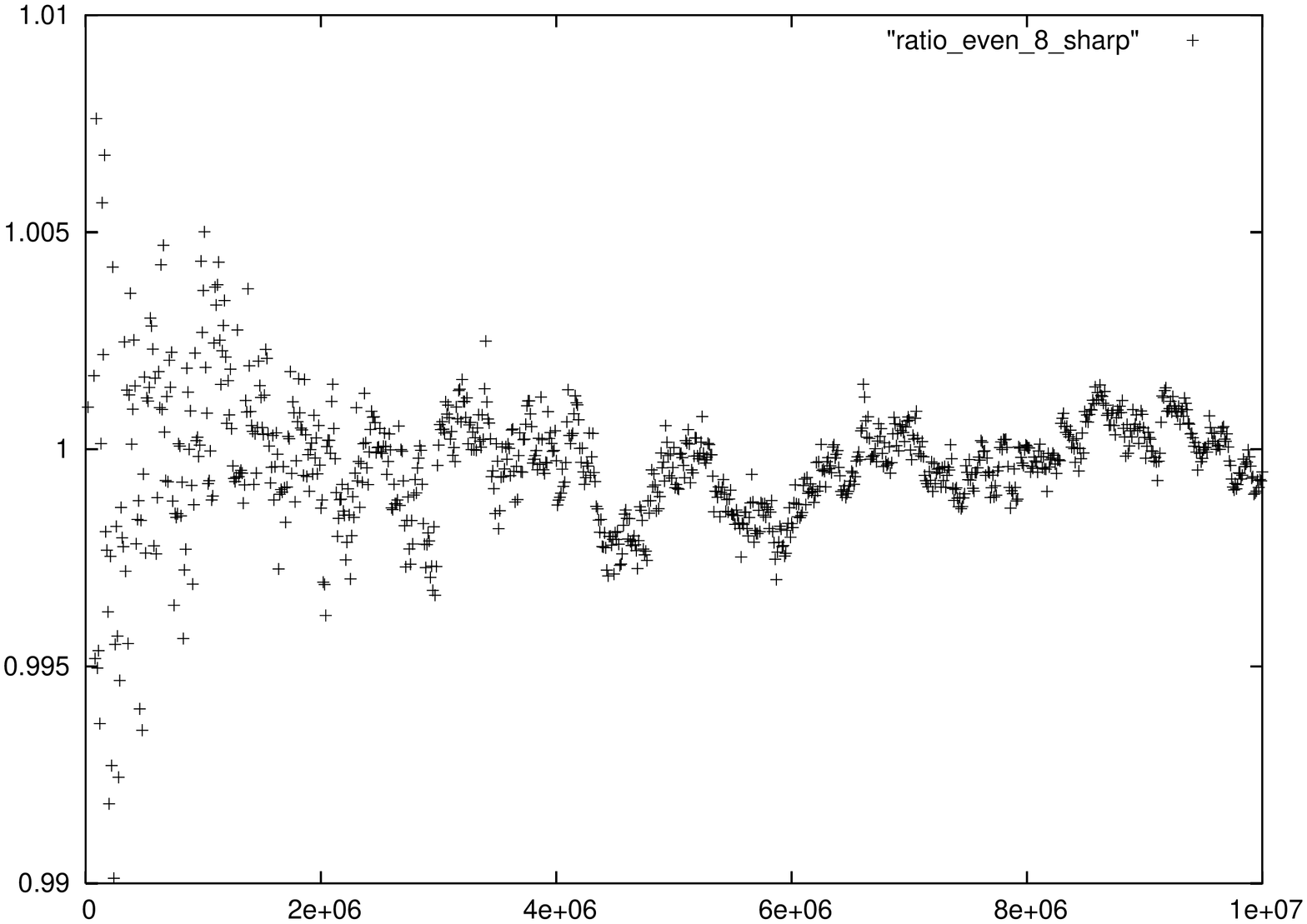,width=3in,height=1.9in,angle=-90}
    }
    \caption
    {Horizontal axis in each graph is $X$. These graphs depict the first eight
     moments, sharp cutoff, of $L(1/2,\chi_d)$, $0< d \leq X$ divided by the conjectured value,
     sampled at $X=10000,20000,\ldots,10^7$.
     We see the graphs fluctuating above and below one. Notice that the vertical
     scale varies from graph to graph}.
    \label{fig:conj vs reality, even chi}
\end{figure}

\subsubsection{Vanishing of $L_E(1/2,\chi_d)$}
\label{sec:L_E values}

In~\cite{CKRS}, Conrey, Keating, the author, and Snaith apply the moment conjecture
~(\ref{eq:ks moment3}) to the problem of predicting asymptotically the
number of vanishings of $L_E(1/2,\chi_d)$.
Using the fact that these
$L$-values are discretized, for example via the Birch and Swinnerton-Dyer conjecture
or the theorem of Kohnen-Zagier~\cite{KZ},
and by studying, up to leading term and for small values,
the density function predicted by~(\ref{eq:ks moment3})
they conjectured that
\begin{equation}
    \notag
    \sum_{{{d \in D(X) \atop (d,q)=1}\atop \text{even funct eqn}}\atop L_E(1/2,\chi_d)=0} 1
    \sim \alpha_E X^{3/4} \log(X)^{\beta_E}.
\end{equation}
The power on the logarithm depends on the underlying curve $E$ because,
in the Birch Swinnerton Dyer conjecture, the Tamagawa factors can contribute
powers of $2$ depending on the prime factors of $d$ and on $E$ 
and this affects the discretization. 
The constant $\alpha_E$ depends on $c_{-1/2}$ and the real period of $E$, but also
on some extra subtle 
arithmetic information that seems to be related to Delaunay's heuristics for
Tate-Shafarevich groups~\cite{De} and is not yet fully understood.
Numerical evidence in favour of this conjecture is presented in~\cite{CKRS2}.
One can skirt these delicate issues, the power on the logarithm and the constant
$\alpha_E$, as follows.

Let $p \nmid q$ be prime.
Sort the $d$'s for which $L_E(1/2,\chi_d)=0$ by residue classes
mod $p$, according to whether $\chi_d(p)=1$ or $-1$, and consider the
ratio
\begin{equation}
    \notag
    R_p(X)=
    \frac{
        \sum_{{{{d \in D(X) \atop (d,q)=1}\atop \text{even funct eqn}}\atop L_E(1/2,\chi_d)=0}\atop\chi_d(p)=1} 1
    }
    {
        \sum_{{{{d \in D(X) \atop (d,q)=1}\atop \text{even funct eqn}}\atop L_E(1/2,\chi_d)=0}\atop\chi_d(p)=-1} 1
    }.
\end{equation}
One can formulate~\cite{CKRS2}~\cite{CFKRS} conjectures for the moments in these 
two subfamilies and the moments agree except for a factor that depends on $p$. 
By considering this ratio, the 
powers of $X$, of $\log X$, and the constant $\alpha_E$ should all cancel out,
except for a single factor that depends on $p$.
This leads to a conjecture~\cite{CKRS} for $R_p(X)$:
\begin{equation}
    \notag
    R_p=\lim_{X\rightarrow \infty}
    R_p(X)=\sqrt{\frac{p+1-a_p}{p+1+a_p}},
\end{equation}
where $a_p$ denotes the $p$th coefficient of the Dirichlet series for $L_E$.
The square root in this conjecture is a consequence of the moments having
a pole at $r=-1/2$. 

\begin{figure}[htp]
    \centerline{
            \psfig{figure=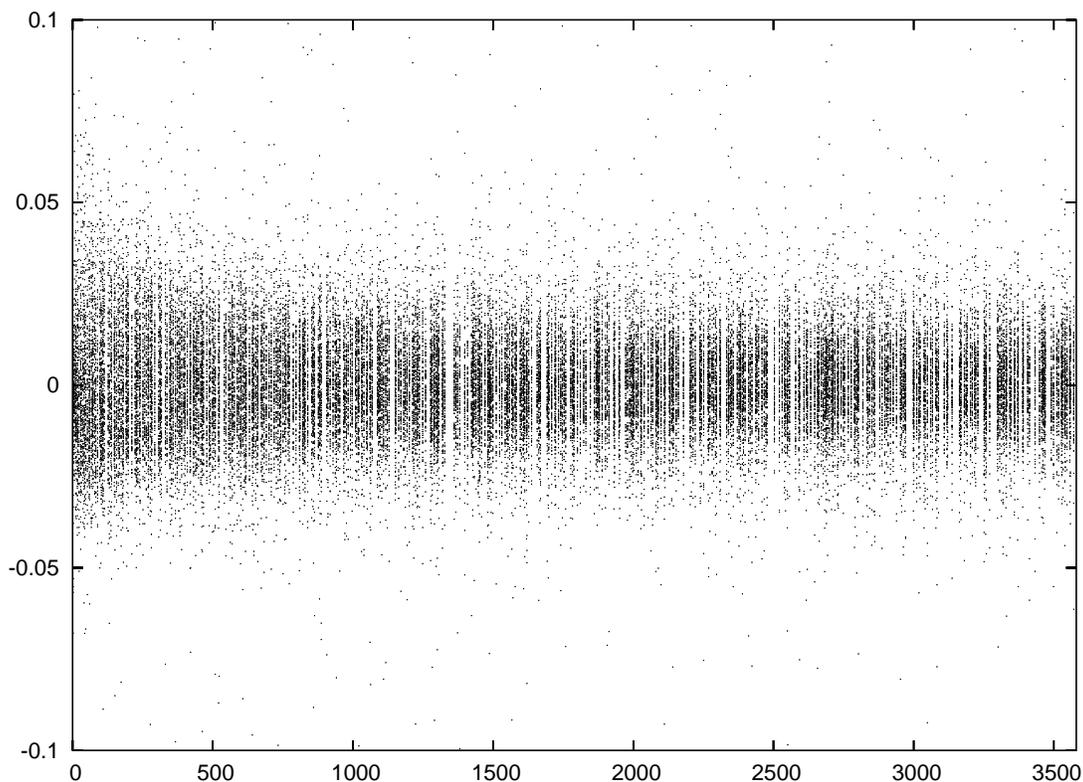,width=6in,angle=-90}
    }
    \caption
    {A plot for one hundered elliptic curves of $R_p(X)-R_p$ for $2\leq p< 3500$, $X=10^8$.}
    \label{fig:1st}
\end{figure}

We end this paper with a plot that substantiates this conjecture.
Figure~\ref{fig:1st} compares, for one hundred elliptic curves $E$,
the predicted value of $R_p$ to the actual value $R_p(X)$, with $X=10^8$
and the set of $d$'s restricted to certain residue classes depending on $E$
as described in ~\cite{CKRS2}. The $L$-values were computed in this special
case by exploiting their connection to the coefficients of certain
weight three halves modular forms and using a table of Rodriguez-Villegas and
Tornaria~\cite{RT}.

The horizontal axis is $p$. For each $p$, and each of the one hundred 
elliptic curves $E$ we plot $R_p(X)-R_p$. We see the values fluctuating
about zero, most of the time agreeing to within about two percent.
The convergence in $X$ is predicted from secondary terms to be 
logarithmically slow and one gets a better fit by including more terms~\cite{CKRS2}.

\subsubsection{Acknowledgements}

The author wishes to thank Andrew Odlyzko for providing 
him with figures~\ref{fig:odlyzko1} and~\ref{fig:odlyzko2} and
the solid curves used 
in figures~\ref{fig:pair correlation} and~\ref{fig:nearest neighbour}.
Nina Snaith supplied the first graph in figure~\ref{fig:value dist three families}.
Atul Pokharel assisted in the preperation of figure~\ref{fig:1st}.
He thanks Fernando Rodriguez-Villegas and Gonzalo Tornaria for giving him 
a list of ternary quadratic forms that he used in computing the $L$-values 
for section~\ref{sec:L_E values}. He is grateful for the kind use of Andrew
Granville's and William Stein's computer clusters on which
some of the computations described were run. Brian Conrey,
David Farmer and Ralph Furmaniak provided feedback on the manuscript.
He also wishes to thank the organizers of the Random Matrix Approaches
in Number Theory program for inviting him to the Newton Institute to participate.


\begin{thebibliography}{R}

\bibitem[AT]{AT}
\newblock
S. Akiyama\ and\ Y. Tanigawa, \emph{Calculation of values of $L$-functions associated to elliptic curves},
Math. Comp. {\bf 68} (1999), no.~227, 1201--1231

\bibitem[B]{MR90b:01039}
\newblock
B. Berndt, \emph{Ramanujan's notebooks. {P}art {I}{I}}, Springer-Verlag,
  New York, 1989.

\bibitem[BK]{BK}
\newblock
M. Berry and J. Keating, \emph{A new asymptotic representation for
  $\zeta(\frac12+it)$ and quantum spectral determinants}, Proc. Roy. Soc.
  London Ser. A \textbf{437} (1992), no.~1899, 151--173.

\bibitem[BK2]{BK2}
\newblock
M. Berry and J. Keating, \emph{The Riemann zeros and eigenvalue asymptotics},
Siam Review {\bf 41} (1999), no.~2, 236--266.

\bibitem[BoK]{BoK}
\newblock
E. Bogomolny and J. Keating, \emph{Gutzwiller's trace formula and spectral
statistics: beyond the diagonal approximation}, Physical Review Letters {\bf 77}
(1996) no.8, 1472--1475.

\bibitem[BoK2]{BoK2}
\newblock
E. Bogomolny and J. Keating, \emph{Random matrix theory and the Riemann zeros II:
$n$-point correlation}, Nonlinearity {\bf 9},
(1996), 911--935.


\bibitem[Bu]{MR97k:11080}
\newblock
D. Bump, \emph{Automorphic forms and representations}, Cambridge Studies in
  Advanced Mathematics, vol.~55, Cambridge University Press, Cambridge, 1997.

\bibitem[BCDT]{BCDT}
\newblock
C. Breuil, B. Conrad, F. Diamond, and R. Taylor, 
J. Amer. Math. Soc. {\bf 14} (2001), no.~4, 843--939.

\bibitem[C]{C}
\newblock 
H. Cohen, {\it High precision computation of Hardy-Littlewood constants}, draft.
Available at {\tt www.math.u-bordeaux.fr/$\sim$cohen}. 

\bibitem[CF]{CF}
\newblock 
B. Conrey and D. Farmer,  {\it Mean values
of $L$-functions and symmetry}, \jour  Internat. Math. Res. Notices
\yr 2000 {\bf 17},  883--908.

\bibitem[CFKRS]{CFKRS}
\newblock   B. Conrey, D. Farmer, J. Keating, M. Rubinstein, and N. Snaith,
{\it Integral moments of $\zeta$ and ${L}$-functions}, 
\jour  Proceedings of the London Mathematical Society, 
to appear.


\bibitem[CFKRS2]{CFKRS2}
\newblock   B. Conrey, D. Farmer, J. Keating, M. Rubinstein, and N. Snaith,
{\it Autocorrelation of random matrix polynomials}, \jour Commun. Math. Phys
{\bf 237} (2003) 3, pp. 365-395.

\bibitem[CFKRS3]{CFKRS3}
\newblock   B. Conrey, D. Farmer, J. Keating, M. Rubinstein, and N. Snaith,
{\it Lower order terms in the moments of $L$-functions}, \jour preprint.

\bibitem[CFZ]{CFZ}
\newblock 
B. Conrey, D. Farmer, and M. Zirnbauer {\it Autocorrelations of ratios
of characteristic polynomials of $L$-functions}, preprint.

\bibitem[CFS]{CFS}
\newblock 
B. Conrey, P. Forrester, and N. Snaith {\it Averages of ratios of characteristic polynomials
for the compact classical groups}, preprint.

\bibitem[CG]{CG}
\newblock
B. Conrey and A. Ghosh,  {\it Mean values of the Riemann zeta-function}, \jour Mathematika
{\bf 31} \yr 1984
159--161.

\bibitem[CG2]{CG2}
\newblock
BJ. Conrey and A. Ghosh,  {\it A conjecture for the sixth power moment of the Riemann zeta-function}, \jour
 Int. Math. Res. Not.
 {\bf 15} \yr 1998
 pp. 775--780.

\bibitem[CGo]{CGo}
\newblock
B. Conrey and S.  Gonek,  {\it High moments of the Riemann zeta-function}, \jour Duke Math. Jour.
\yr 2001
 {\bf 107}  pp. 577--604.

\bibitem[CKRS]{CKRS}    
B. Conrey, J. Keating, M. Rubinstein, and N. Snaith,
{\it On the frequency of vanishing of quadratic twists of modular ${L}$-functions},
in {\it Number Theory for the Millennium I: Proceedings of the
Millennial Conference on Number Theory}; editor, M.A.~Bennett et~al., pages
301--315. A K Peters, Ltd, Natick, 2002.

\bibitem[CKRS2]{CKRS2}    
B. Conrey, J. Keating, M. Rubinstein, and N. Snaith,
{\it Random Matrix Theory and the Fourier Coefficients of
Half-Integral Weight Forms}, ar{X}iv:math.nt/0412083

\bibitem[CS]{CS}
\newblock 
B. Conrey and N.Snaith, {\it Applications of the $L$-functions ratios conjectures},
preprint.

\bibitem[D]{D}
\newblock
H. Davenport,  Multiplicative Number Theory, \publ GTM 74
\publaddr  Springer-Verlag, New York, NY
\yr 2000.

\bibitem[DH]{DH}
\newblock
D.~Davies and C. Haselgrove, \emph{The evaluation of {D}irichlet
  ${L}$-functions}, Proc. Roy. Soc. Ser. A \textbf{264} (1961), 122--132.

\bibitem[De]{De}
\newblock
C. Delaunay, {\it Heuristics on {T}ate-{S}hafarevitch groups of elliptic curves defined over {$\mathbb{Q}$}},
   Experiment. Math.,
    {\bf 10} (2001), 2 , 191--196.

\bibitem[Del]{MR49:5013}
\newblock
P. Deligne, \emph{La conjecture de {W}eil. {I}}, Inst. Hautes \'Etudes Sci.
  Publ. Math. (1974), no.~43, 273--307.

\bibitem[Do]{Do}
\newblock
T. Dokchister,
{\it Computing Special Values of Motivic $L$-Function}, ar{X}iv:math.NT/0207280

\bibitem[E]{E}
\newblock
H. Edwards,
{\it Riemann's Zeta Function}, 
\publ Academic Press 
\yr 1974

\bibitem[EMOT]{MR84h:33001b}
\newblock
A. Erd{\'e}lyi, W. Magnus, F. Oberhettinger, and F. Tricomi, 
\emph{Higher transcendental functions. {V}ol. {I}{I}}, Robert E.
  Krieger Publishing Co. Inc., Melbourne, Fla., 1981, Based on notes left by
  Harry Bateman, Reprint of the 1953 original.

\bibitem[FKL]{FKL}
\newblock
D. Farmer, W. Kranec, and S. Lemurell,
{\it Maass forms on $\Gamma_0(11)$}, draft.

\bibitem[F]{F}
\newblock
S. Fermigier, \emph{Z\'eros des fonctions ${L}$ de courbes
  elliptiques}, Experiment. Math. \textbf{1} (1992), no.~2, 167--173.

\bibitem[Fr]{Fr}
\newblock
E. Friedman, 
{\it Hecke's integral formula},
S\'eminaire de Th\'eorie des Nombres, 1987--1988,
Exp.\ No.\ 5, 23,
Univ. Bordeaux I.

\bibitem[G]{G}
\newblock
W. Gabcke, 
{\it Neue Herleitung und explicite Restabschatzung der Riemann-SiegelFormel}, 
\jour Ph.D. Dissertation, 
\publaddr Gottingen 
\yr 1979.

\bibitem[GHK]{GHK}
\newblock
S. Gonek, C. Hughes, J, Keating, {\it A New Statistical Model for the Riemann Zeta Function}, preprint.

\bibitem[H]{H}
\newblock
R. Heath-Brown,  {\it The fourth power moment of the Riemann zeta-function}, \jour Proc. London Math. Soc.
 (3)
\yr 1979      
 {\bf 38}  pp. 385 -- 422.

\bibitem[I]{I}
\newblock
 A. E. Ingham,  {\it Mean-value theorems in the theory of the Riemann zeta-function}, \jour Proceedings of the
 London Mathematical Society  (92)
\yr 1926
 {\bf 27}  pp. 273--300.

\bibitem[J]{J}
\newblock     
 M. Jutila,  {\it On the mean value of $L({\frac12},\,\chi )$ for real characters}, \jour Analysis
 {\bf 1} \yr 1981
\pages149--161.

\bibitem[KS]{KS}
\newblock
N. Katz and P. Sarnak, {\it Random matrices, Frobenius eigenvalues, and monodromy}, 
Amer. Math. Soc., Providence, RI \yr 1999.

\bibitem[KS2]{KS2}
\newblock
N. Katz and P. Sarnak, {\it Zeroes of zeta functions and symmetry},
Bull. Amer. Math. Soc. (N.S.) {\bf 36} (1999), no.~1, 1--26.

\bibitem[KS3]{KS3}
\newblock
N. Katz and P. Sarnak, {\it Zeros of zeta functions, their spacings
and their spectral nature}, 1997 preprint of KS2.


\bibitem[K]{K}
\newblock
J. Keating,  {\it Periodic orbits, spectral statistics, and the Riemann zeros},
in Supersymmetry and Trace Formulae: Chaos and Disorder, 
J. Keating, D. Khmelnitskii, and I. Lerner, eds., Plenum, New York, 1998, 1--15.

\bibitem[KeS]{KeS}
\newblock
J. Keating and N. Snaith,  {\it Random matrix theory and $\zeta(\frac12 +it)$}, \jour Comm. Math. Phys. {\bf 214} \yr 2000
57--89.

\bibitem[KeS2]{KeS2}
\newblock   
 J. P. Keating and N. C. Snaith,  
{\it Random matrix theory and $L$-functions at $s=\frac12 $}, \jour Comm. Math. Phys.
{\bf 214} \yr 2000
 pp.  91--110. 

\bibitem[Ke]{Ke}
\newblock
J. Keiper, \emph{On the zeros of the {R}amanujan $\tau$-{D}irichlet series 
  in the critical strip}, Math. Comp. \textbf{65} (1996), no.~216, 1613--1619.

\bibitem[KZ]{KZ}
W.~Kohnen and D.~Zagier,
{\it Values of ${L}$-series of modular forms at the center of the critical strip},
Invent. Math., {\bf 64} (1981), 175--198.    

\bibitem[Kn]{MR93j:11032}
\newblock
A. Knapp, \emph{Elliptic curves}, Mathematical Notes, vol.~40,
  Princeton University Press, Princeton, NJ, 1992.

\bibitem[LO]{MR80g:12010}
\newblock
J. Lagarias and A. Odlyzko, \emph{On computing {A}rtin ${L}$-functions in
  the critical strip}, Math. Comp. \textbf{33} (1979), no.~147, 1081--1095.

\bibitem[L]{L}
\newblock
A. Lavrik, {\it Approximate functional equation for Dirichlet Functions}, 
Izv. Akad. Nauk SSSR {\bf 32} (1968), 134--185.

\bibitem[Le]{Le}
\newblock
R.~Lehman, {\em On the distribution of the zeros of the Riemann zeta
function}, Proc. London Math. Soc. (3) {\bf 20} (1970), 303-320.
\MR{41:3414}

\bibitem[LRW]{LRW}
\newblock
J.~van~de Lune, H. te~Riele, and D. Winter, \emph{On the zeros of the
  {R}iemann zeta function in the critical strip. {I}{V}}, Math. Comp.      
  \textbf{46} (1986), no.~174, 667--681.

\bibitem[M]{M}
\newblock
M. Mehta, {\it Random Matrices}, 2nd edition, Academic Press, 1991.

\bibitem[Mo]{Mo}
\newblock
H. Montgomery, \emph{The pair correlation of zeros of the zeta function},
Analytic number theory (Proc. Sympos. Pure Math., Vol. XXIV, St. Louis Univ.,
St. Louis, Mo., 1972) (Providence, R.I.), Amer. Math. Soc., 1973,    
pp.~181--193.

\bibitem[O]{O}
\newblock
A. Odlyzko, 
{\it The $10^{20}$-th zero of the Riemann zeta function and 175 million of its neighbors}, unpublished.
{\tt www.dtc.umn.edu/$\sim$odlyzko}

\bibitem[O2]{O2}
\newblock
A. Odlyzko, 
{\it The $10^{22}$-nd zero of the Riemann zeta function}, 
Dynamical, Spectral, and Arithmetic Zeta Functions, M. van Frankenhuysen and M. L. Lapidus, eds., 
Amer. Math. Soc., Contemporary Math. series, {\bf 290}, 2001, 139--144.

\bibitem[O3]{O3}
\newblock
A. Odlyzko, 
{\it On the distribution of the spacings between zeros of the zeta function}, 
Math. Comp., {\bf 48} (1987), 273--308.

\bibitem[O4]{O4}
\newblock
A. Odlyzko, private communication.

\bibitem[OS]{OS}
\newblock
A. Odlyzko and A. Sch\"{o}nhage, 
{\it Fast algorithms for multiple evaluations of the Riemann zeta function}
\jour  Trans. Am. Math. Soc., {\bf 309}
\yr 1988, 797--809.

\bibitem[Og]{MR41:1648}
\newblock
A. Ogg, \emph{Modular forms and {D}irichlet series}, W. A. Benjamin, Inc.,
  New York-Amsterdam, 1969.

\bibitem[Ol]{MR55:8655}
\newblock
F. Olver, \emph{Asymptotics and special functions}, Academic Press [A
  subsidiary of Harcourt Brace Jovanovich, Publishers], New York-London, 1974,
  Computer Science and Applied Mathematics.

\bibitem[OzS]{OzS} 
\newblock
A. \"Ozl\"uk\ and\ C. Snyder, 
\emph{Small Zeroes of Quadratic
$L$-Functions}, Bull. Aust. Math. Soc. {\bf 47} (1993), 307--319.

\bibitem[OzS2]{OzS2} 
\newblock
A. \"Ozl\"uk\ and\ C. Snyder, 
\emph{On the distribution of the nontrivial zeros of 
quadratic $L$-functions close to the real axis},
Acta Arith. {\bf 91} (1999), no.~3, 209--228.

\bibitem[P]{P}
\newblock
R. Paris, \emph{An asymptotic representation for the {R}iemann zeta function
  on the critical line}, Proc. Roy. Soc. London Ser. A \textbf{446} (1994),
  no.~1928, 565--587.

\bibitem[R]{MR91j:01070b}
\newblock
Bernhard Riemann, \emph{Gesammelte mathematische {W}erke, wissenschaftlicher
  {N}achlass und {N}achtr\"age}, Springer-Verlag, Berlin, 1990, Based on the
  edition by Heinrich Weber and Richard Dedekind, Edited and with a preface by
  Raghavan Narasimhan.

\bibitem[RT]{RT}
\newblock
F. Rodriguez-Villegas and G. Tornaria, private communication.

\bibitem[Ru]{Ru}
\newblock
M. Rubinstein,
{\it Evidence for a spectral interpretation of the zeros of
      $L$-functions}. Princeton Ph.D. Dissertation, 1998.

\bibitem[Ru2]{Ru2}
\newblock
M. Rubinstein,
{\it Low lying zeros of $L$-functions and random matrix
theory}. Duke Mathematical Journal {\bf 109} (2001), no.~1,
147--181.

\bibitem[Ru3]{Ru3}
\newblock
M. Rubinstein,
{\it Lower terms in the density of zeros of quadratic Dirichlet $L$-functions},
preprint.

\bibitem[Ru4]{Ru4}
\newblock
M. Rubinstein,
{\it The $L$-function class library and command line interface},
{\tt www.math.uwaterloo.ca/$\sim$mrubinst/L\_function/L.html}.


\bibitem[Rud]{MR88k:00002}
\newblock
W. Rudin, {\it Real and complex analysis}, third ed., McGraw-Hill Book
Co., New York, 1987.

\bibitem[Rum]{Rum}
\newblock
R.~Rumely, {\em Numerical computations concerning the ERH}, Math. Comp. {\bf
61} (1993), 415--440.

\bibitem[RS]{RS}
\newblock
M. Rubinstein and P. Sarnak, {\it Chebyshev's Bias},
\jour Experimental Mathematics {\bf 3} \yr 1994, no.~3, 173--197.

\bibitem[RudS]{RudS}
\newblock
Z. Rudnick and P. Sarnak, {\it Zeros of principal $L$-functions and random matrix theory},
Duke Mathematical Journal (2) \textbf{81} (1996), 269--322.

\bibitem[S]{S}
\newblock   K. Soundararajan,  {\it Non-vanishing of
quadratic Dirichlet $L$-functions at $s=\frac12$ }, \jour  Ann. of Math.
(2)  {\bf 152} \yr 2000  pp. 447--488.

\bibitem[Sp]{Sp}
\newblock
R. Spira, \emph{Calculation of the {R}amanujan $\tau $-{D}irichlet series},
  Math. Comp. \textbf{27} (1973), 379--385.

\bibitem[St]{St}
\newblock
A. Strombergsson, {\it On the zeros of $L$-functions associated to Maass waveforms},
IMRN (1999), No. 15.

\bibitem[TW]{MR96d:11072}
\newblock
R. Taylor and A. Wiles, \emph{Ring-theoretic properties of certain
  {H}ecke algebras}, Ann. of Math. (2) \textbf{141} (1995), no.~3, 553--572.

\bibitem[T]{T}
\newblock
N. Temme, {\it The asymptotic expansions of the incomplete gamma
functions}, SIAM J. Math. Anal. {\bf 10} (1979), 757--766.

\bibitem[Ti]{MR88c:11049}
\newblock
E. Titchmarsh, \emph{The theory of the {R}iemann zeta-function}, second ed.,
  The Clarendon Press Oxford University Press, New York, 1986, Edited and with
  a preface by D. R. Heath-Brown.

\bibitem[To]{To}
\newblock
E. Tollis, \emph{Zeros of Dedekind zeta functions in the critical
  strip}, Math. Comp. \textbf{66} (1997), no.~219, 1295--1321.

\bibitem[Tu]{Tu}
\newblock
A. Turing, {\em Some calculations of the Riemann zeta function}, Proc.
London Math. Soc. (3) {\bf 3} (1953), 99--117. 


\bibitem[W]{W}
\newblock
S. Wedeniwski, {\it Verification of the Riemann Hypothesis}, 
{\tt www.zetagrid.net}.

\bibitem[Wi]{MR96d:11071}
\newblock
A. Wiles, \emph{Modular elliptic curves and {F}ermat's last theorem}, Ann.
  of Math. (2) \textbf{141} (1995), no.~3, 443--551.

\bibitem[Y]{Y}
\newblock
H. Yoshida, \emph{On calculations of zeros of ${L}$-functions related
  with {R}amanujan's discriminant function on the critical line}, J. Ramanujan
  Math. Soc. \textbf{3} (1988), no.~1, 87--95.

\bibitem[Y2]{Y2}
\newblock
H. Yoshida, \emph{On calculations of zeros of various ${L}$-functions}, J. Math. 
  Kyoto Univ. \textbf{35} (1995), no.~4, 663--696.

\end{thebibliography}
\end{document}